\documentclass[12pt]{article}
\usepackage[russian]{babel}
\usepackage{amsmath,amssymb,amsthm,amsfonts,amscd}
\usepackage{epsfig}
\input epsf

 \def\R{{\mathbb R}}\def\C{{\mathbb C}}  \def\N{{\mathbb N}}

\def\ctg{\mathop{\fam0 ctg}}
\def\ch{\mathop{\fam0 ch}}
\def\det{\mathop{\fam0 det}}
\def\diam{\mathop{\fam0 diam}}
\def\exp{\mathop{\fam0 exp}}
\def\grad{\mathop{\fam0 grad}}
\def\pr{\mathop{\fam0 pr}}
\def\SO{\mathop{\fam0 SO}}
\def\tg{\mathop{\fam0 tg}}
\def\tr{\mathop{\fam0 tr}}
\def\id{\mathop{\fam0 id}}
\def\Alt{\mathop{\fam0 Alt}}
\long\def\comment#1\endcomment{}

\begin{document}

\centerline{\uppercase{\bf Basic differential geometry}}
\smallskip
\centerline{\uppercase{\bf as a sequence of interesting problems}}
\smallskip
\centerline{A. Skopenkov}

\bigskip
\small
{\bf Abstract.}
This book is expository and is in Russian (sample English translation of two pages is given).
It is shown how in the course of solution of interesting geometric problems
(close to applications) naturally appear different notions of
{\it curvature}, which distinguish given geometry from the `ordinary' one.
Direct elementary definitions of these notions are presented.
The book is accessible for students familiar with analysis of several
variables, and could be an interesting easy reading for professional ma\-the\-ma\-ti\-cians.
The material is presented as a sequence of problems, which is peculiar not only
to Zen monasteries but also to serious mathematical education.

\normalsize
\bigskip
\hfill{\it The modern world is full of theories which are proliferating}

\hfill{\it at a wrong level of generality, we're so} good {\it at theorizing,}

\hfill{\it and one theory spawns another, there's a whole industry}

\hfill{\it of abstract activity which people mistake for thinking.}

\hfill{\it I. Murdoch, The Good Apprentice.}

\bigskip
\centerline{\bf Main curvatures}

\smallskip
{\bf 1.} Let $A_1,\dots,A_s$ be points in $\R^3$ with masses $m_1,\dots,m_s$.
The {\it inertia momentum} of this system w.r.t. a line $l$ is the number
$I(l):=m_1|A_1l|^2+\dots+m_s|A_sl|^2$, where $|A_il|$ is the distance from
$A_i$ to $l$.

(a) Let $I_+$ and $I_-$ be the maximal and the minimal value of the inertia
momenta w.r.t. lines in the plane pasing through a fixed point $O$
(posibly, $I_+=I_-$).
Take one of the lines $l_+$ such that $I(l_+)=I_+$.
Then
$I(l)=I_+\cos^2\varphi+I_-\sin^2\varphi$, где $\varphi=\angle(ll_+)$.

(b)* For each $O\in\R^3$ there are lines $l_1,l_2,l_3$ passing through $O$ and
such that for each line $l$ pasng through $O$ we have
$I(l)=I(l_1)\cos^2(ll_1)+I(l_2)\cos^2(ll_2)+I(l_3)\cos^2(ll_3)$.

\smallskip
\centerline{\epsffile{p.20}}
\centerline{\it Figure: a normal and a `skew' sections of a surface.}

\smallskip
A {\it coorientation} of a surface $\Pi\subset\R^3$ is a field of unit lengh
vectors $n(P)$ normal to $\Pi$; the field should continuously depend on
$P\in\Pi$.

A {\it curvature of a (non-parametrized) curve on a cooriented surface}
is the normal projection of the acceleration of the unit length velocity
motion on this curve.

{\it Main curvatures} $\lambda_+$ and $\lambda_-$ of a cooriented surface
$\Pi\subset\R^3$ at $P\in\Pi$ are the maximal and the minimal values of
the curvatures at $P$ of (non-parametrized) curves on $\Pi$ that are the
intersection of $\Pi$ and planes passing through $n(P)$.
The corresponding planes are denoted by $\alpha_+$ and $\alpha_-$.


\smallskip
{\bf 2.} How do main curvatures change under

(a) change of the coorientation to the opposite? \qquad

(b) dilatation of the space?

\smallskip
{\bf 3.} (a) {\it The Euler formula.}
Let $\Pi\subset\R^3$ be a cooriented surface and $P\in\Pi$.
Denote by $k(\varphi)$ the curvature at $P$ of the curve on $\Pi$ that is the
intersection of $\Pi$ and the plane passing through $n(P)$ and making the angle
$\varphi$ to the plane $\alpha_+$ (for which such a curvature is maximal).
Then
$k(\varphi)=\lambda_+\cos^2\varphi+\lambda_-\sin^2\varphi$.

(b) If $\lambda_+\ne\lambda_-$, then $\alpha_+\perp\alpha_-$.

\smallskip
{\bf 4.}
(a) What is the ratio of curvatures for the intersection curves of a surface
with two planes ($\alpha$ and $\alpha_n$ in Fig. 2) containing $P$ and
intersecting the tangent plane by the same line, if $\alpha_n$ passes through
the normal vector at $P$ and $\alpha$ makes angle $\theta$ to the normal
vector?

(b) The normal projection at $P=\gamma(0)$ of the acceleration of a parametrized curve $\gamma$ on a surface $\Pi$ depends only on the velocity $\gamma'(0)$ of this curve at $P$.

(c) State and prove the analogue of problems 3.a and 4.a for a cooriented 3-dimensional
surface $\Pi\subset\R^4$.

\smallskip
{\bf Hints.}

{\bf 1a and 3a.} Follow because both inertia moment and the curvature
have the form
$f(\varphi)=A\cos^2\varphi+2B\cos\varphi\sin\varphi+C\sin^2\varphi$.

{\bf 1b and 4c.} Follow because both inertia moment and the curvature are
quadratic forms (prove this!).

{\bf 4.} (a) {\it The Meunieur Theorem.} $k\cos\theta=k_n$.
Follows by (b).

(b) Denote $n=n(\gamma(t))$.
Then
$$n\cdot\gamma'=0\quad\Rightarrow\quad n'\cdot\gamma'+n\cdot\gamma''=0
\quad\Rightarrow\quad n\cdot\gamma''=-\gamma'\cdot\partial n/\partial\gamma'.$$

\bigskip
\centerline{\bf A simple geometric definition of the scalar and the Ricci curvatures}


\smallskip
{\bf 1.*} Let $\Pi\subset\R^3$ be a surface of revolution.
Denote by $L(R)$ the length of the circle in
$\Pi$ of radius $R$ centered at $P\in\Pi$.
Prove that

(a) $\lim\limits_{R\to0}\dfrac{2\pi R-L(R)}{R^2}=0.$

(b) there exists $\tau:=6\lim\limits_{R\to0}\dfrac{2\pi R-L(R)}{R^3}$.
This is the {\bf scalar curvature} of $\Pi$ at $P$.

\smallskip
Let $\Pi\subset\R^m$ be an $n$-dimensional surface in $\R^m$ (a reader not
familiar with surfaces can consider surfaces of revolution, which case is
interesting enough).

If $n=2$, for $P\in\Pi$ denote by $L_{\Pi,P}(R)$ the length of the circle in
$\Pi$ of radius $R$ centered at $P\in\Pi$.
The {\bf scalar curvature} of $\Pi$ at $P$ is the number
$$\tau=\tau_{\Pi,P}:=6\lim\limits_{R\to0}\frac{2\pi R-L_{\Pi,P}(R)}{\pi R^3}.$$
One can prove that this limit indeed exists.

One can prove that $\tau=2\lambda_+\lambda_-$ for 2-surfaces in $\R^3$.

A (non-parametrized) curve $\Gamma\subset\Pi$ is called a {\bf geodesic} on
$\Pi$, if $\Gamma$ is {\it locally shortest}, i.e., if each point $x\in\Gamma$
has a neighborhood $U\subset\Pi$ such that the distance (along $\Pi$)
between any two points $y_1,y_2\in U\cap\Gamma$ is equal to the length
of a segment of $\Gamma$ from $y_1$ to $y_2$.

For $P\in\Pi$ denote by $T_P$ the plane tangent to $\Pi$ at the point $P$.
Define a {\bf (geodesic) exponential map}
$$\exp=\exp\phantom{}_P:T_P\to\Pi\quad\mbox{by}
\quad\exp(u):=\gamma_{P,u}(1),$$
where $\gamma_{P,u}:[-1,1]\to\Pi$ is the geodesic for which $\gamma_{P,u}(0)=P$ and
$\gamma'_{P,u}(0)=u$.

{\bf The Ricci curvarure} of $\Pi$ at $P\in\Pi$ is a bilinear form
$\rho=\rho_P:T_P\times T_P\to R$ such that for each unit
$n$-dimensional cube $A\subset T_P$ with vertex at $P$ we have
$$V(\exp\phantom{}_P(hA))=h^n-\frac{h^{n+2}}6\int_A\rho(u,u)du+O(h^{n+3})
\quad\mbox{when}\quad h\to0.$$

One can prove that
\quad

$\bullet$ such a quadratic form indeed exists and is unique.

$\bullet$ analogous formula with $h^n$ replaced by $h^nV(A)$ ($h^{n+2}$ and $h^{n+3}$
do not change) holds for each measurable subset $A\subset T_P$.

$\bullet$ $\tau=\tr\widetilde\rho$, where linear operator
$\widetilde\rho:T_P\to T_P$ is well-defined by
$\widetilde\rho(u)\cdot v=\rho(u,v)$.

$\bullet$ $2\rho(u,u)=\sum_i\tau_i$, where $|u|=1$, $e_1,\dots,e_{n-1}$ is an
orthonormal base in the orthogonal complement to $u$ in $T_P$
and $\tau_i$ is the scalar curvature at $P$ of the
2-dimensional surface $\exp_P\left<u,e_i\right>$.

\bigskip
[1] L. Bessieres, G. Besson, M. Boileau, La preuve de la conjecture de
Poincar\'e d'apres G. Perelman, Images des Mathematiques, 2006.
http://www.math.cnrs.fr/imagesdesmaths/IdM2006.htm



\newpage
\normalsize
\centerline{\uppercase{\bf Основы дифференциальной геометрии}}
\bigskip
\centerline{\uppercase{\bf в интересных задачах}
\footnote{Прошлая версия опубликована в издательстве МЦНМО, Москва, 2009 и 2010.
Данная обновленная версия выложена на http://arxiv.org/abs/0801.1568 с разрешения издательства.}
}
\bigskip
\centerline{А. Б. Скопенков
\footnote{Московский Физико-Технический Институт,
Независимый Московский Университет,
Инфо: www.mccme.ru/\~{ }skopenko.}
}

\bigskip
{\bf Аннотация}

Показано, как при решении интересных геометрических проблем,
близких к приложениям, естественно возникают различные понятия
{\it кривизны}, отличающей изучаемую геометрию от 'обычной'.
Приведены прямые элементарные определения этих понятий.
Брошюра предназначена студентам, аспирантам, работникам науки и образования,
изучающим и применяющим дифференциальную геометрию.
Для ее изучения достаточно владения основами анализа функций нескольких
переменных (а во многих местах не нужно даже этого).
Материал преподнесен в виде циклов задач.
Это характерно не только для дзенских монастырей, но и для серьезного изучения математики.


 \bigskip
\hfill{\it The modern world is full of theories which are proliferating}

\hfill{\it at a wrong level of generality, we're so} good {\it at theorizing,}

\hfill{\it and one theory spawns another, there's a whole industry}

\hfill{\it of abstract activity which people mistake for thinking.}

\hfill{\it I. Murdoch, The Good Apprentice.}

\tableofcontents

\newpage
\section*{Введение}
\addcontentsline{toc}{section}{Введение}

\subsection{Зачем}

Приводимые задачи подобраны так, что в процессе их решения (и обсуждения)
читатель увидит, как при решении интересных геометрических проблем,
близких к приложениям, естественно возникают различные понятия
{\it кривизны}, отличающей изучаемую геометрию от 'обычной'.
\footnote{Тем самым он освоит основы дифференциальной геометрии (в частности,
б\'ольшую часть курсов, изучаемых на факультете инноваций и высоких технологий Московского Физико-технического Института и на механико-математическом факультете
Московского Государственного Университета им. М. В. Ломоносова --- кроме
интегрирования дифференциальных форм и основ топологии).
Дальнейшие знания читатель сможет почерпнуть в книгах из списка литературы.}

Особенность этого текста --- возможность познакомиться с некоторыми {\it мотивировками} и {\it идеями}
дифференциальной геометрии при сведении к необходимому минимуму ее {\it языка}.
Я старался давать определения так, чтобы сразу было ясно, что определяемый объект интересен.
А методы вычисления уже интересных (по самому их определению) объектов
формулировать в виде теорем.
(Часто изучение материала затрудняется тем, что {\it вычислительные формулы}
преподносятся в виде {\it определений}, которые становятся немотивированными.)
Вместо {\it абстрактных общих понятий} (например, тензора и ковариантного
дифференцирования) рассматриваются их {\it конкретные используемые в курсе
частные случаи}, а обобщение остается в виде задач, которые естественны и легки
для читателя, разобравшегося с частными случаями.
\footnote{Изучение 'от общего к частному' часто приводит к абсурдному эффекту: сдающие
курс воспроизводят громоздкое определение, но не могут по этому определению
привести ни одного содержательного примера определяемого объекта.}

Простейшие кривизны --- числовые поля, более сложные --- поля квадратичных
форм, а тензор кривизны Римана ('это маленькое чудовище полилинейной алгебры'
по словам М. Громова) --- поле четырехлинейных форм.
В этом курсе даются прямые геометрические определения сначала первых, затем
вторых и потом третьего.
Конечно, простейшие кривизны выражаются через более сложные (и такие выражения
часто удобны для {\it вычисления} простейших кривизн), но {\it определение}
простых понятий через более сложные затрудняет изучение материала.

Ввиду прозрачной геометрической мотивированности изучаемых понятий изложение
в основном синтетично и бескоординатно.
Несмотря на стремление к ясности и ориентированность на приложения (а точнее,
как раз в силу такого стремления), я старался
поддержать достаточно высокий уровень строгости.
Например,
различаются параметризованные и непараметризованные кривые и поверхности
(отсутствие их четкого различения мешает начинающим, хотя допустимо и удобно для специалистов; см. также сноску в
\S\ref{0curves}).

Принятый стиль изложения отвечает духу К. Ф. Гаусса (и других
первооткрывателей), много занимавшегося приложениями и превратившего один
из разделов географии в данный раздел математики.
Изложение 'от простого к сложному' и в форме, близкой к форме рождения
материала, продолжает устную традицию, восходящую к Лао Цзы и Платону,
а в современном преподавании
математики представленную, например, книгами Д. Пойа и журналом 'Квант'.

Мне кажется, принятый стиль изложения не только сделает материал более
доступным, но позволит сильным студентам (для которых доступно даже абстрактное
изложение) приобрести математический вкус с тем, чтобы разумно выбирать
проблемы для исследования, а также ясно излагать собственные открытия, не
скрывая ошибок (или известности полученного результата) за чрезмерным
формализмом.
К сожалению, такое (бессознательное) сокрытие ошибки часто происходит с
молодыми математиками, воспитанными на чрезмерно формальных курсах (происходило
и с автором этих строк; к счастью, почти все мои ошибки исправлялись {\it
перед} публикациями).

Чтение этого текста и решение задач потребуют от большинства читателей усилий
(впрочем, некоторые читатели данного текста жаловались, что в нем нет серьезных
задач, а есть лишь тривиальные упражнения).
Однако эти усилия будут сполна оправданы тем, что вслед за великими
математиками в процессе изучения интересных геометрических проблем читатель
откроет некоторые основные понятия и теоремы дифференциальной геометрии.
Надеюсь, это поможет читателю совершить собственные настолько же полезные
открытия (не обязательно в математике)!

Этот текст основан на лекциях и семинарах, которые автор вел на мехмате МГУ
в 2004-2007 годах, в летней школе 'Современная Математика' в 2007 году и на ФИВТ МФТИ в 2013 году.
Его фрагменты были представлены на семинаре кафедры дифференциальной геометрии и приложений мехмата МГУ
(рук. А.Т. Фоменко) и на семинаре по геометрии в МЦНМО (рук. В.Ю. Протасов).
Благодарю А. Иванова, А. Клименко, С. Маркелова, А. Ошемкова, А. Пляшечника, В. Прасолова,
А. Толченникова, Ю. Торхова, Г. Челнокова и всех слушателей (точнее, решателей)
курсов за полезные замечания и обсуждения.
Благодарю А. Ошемкова и В. Прасолова за предоставление тех-файлов их материалов.

\subsection{Советы и соглашения}

Для понимания условий и для решения задач достаточно уверенного владения
основами анализа функций нескольких переменных (и, чем дальше, тем больше, линейной алгебры).
Все необходимые {\it новые} определения приводятся здесь.
Кое-где требуется также теорема о существовании и единственности решения дифференциального уравнения.

Важные факты
выделены словом `теорема' или `следствие'.
Иногда подсказками являются соседние задачи; указания даются в конце каждой темы.
Факты, для доказательства которых читателю может понадобиться литература (или консультация специалиста),
приводятся со ссылками.
Если условие задачи является формулировкой утверждения, то в задаче требуется это утверждение доказать.

Рассматриваемые понятия и факты интересны, полезны и нетривиальны даже для поверхностей вращения и графиков функций
(в основном в трехмерном пространстве), а также для поверхностей многогранников.
Например, инвариант Дена, с помощью которого была решена 3-я проблема
Гильберта, тесно связан со средней кривизной поверхности многогранника.
Поэтому не приводится примеров более сложных поверхностей (кроме плоскости Лобачевского в самом конце).
Однако для хорошего понимания материала читателю будет полезно изучить такие примеры [Ra03, MF04].

Определения даются в предположении, что определяемые объекты (в частности, пределы, производные и интегралы) существуют.
Задачи также формулируются в предположении, что заданные в их условиях объекты существуют.

Заданные в условиях функции предполагаются бесконечно дифференцируемыми, если не оговорено противное.
Векторы и вектор-функции обозначаются обычными (не жирными) буквами без стрелочек.
Читатель легко отличит их от скаляров по контексту.

Через $\cdot$, $\times$ и $\wedge$ обозначаются скалярное, векторное и
смешанное (для векторов) или внешнее (для дифференциальных форм) произведения, соответственно.

Приводимые в \S\ref{0cur}, \S\ref{0poly} определения кривизн {\it независимы друг от друга}.
Поэтому после изучения поверхностей можно сразу изучать {\it любую} из вводимых
здесь кривизн (для скалярной, средней и гауссовой кривизн необходимо еще
понятие площади, для секционной и римановой --- параллельного переноса, а для
риччиевой --- геодезических и экспоненциального отображения).
При этом, естественно, задачи о связи изучаемой кривизны с еще не изученными придется отложить на потом.

{\bf Студентам, изучающим курс.}
Задачи составлены так, что их решение (вместе с изучением лекционного курса) поможет не только успешно сдать зачет
и экзамен, но также освоить основы дифференциальной геометрии и вообще стать разумным человеком.
Некоторые задачи --- на лекционный материал.
Как и для других задач, нужно быть готовым рассказать  у доски  решения ---
включая детали доказательств, которые могли не разбираться на лекции.
Поэтому, если не оговорено противное, то теоремами, доказанными на лекциях, пользоваться без доказательства
в решениях нельзя (при этом иногда проще не повторить доказательство лекционной теоремы и использовать ее для
решения задачи, а повторить необходимый фрагмент доказательства на примере решения задачи).

\subsection{Литература}


[BBB06] L. Bessieres, G. Besson, M. Boileau, La preuve de la conjecture de
Poincar\'e d'apres G.Perelman, Images des Mathematiques, 2006,
\linebreak
http://www.math.cnrs.fr/imagesdesmaths/IdM2006.htm.
Есть рус. перевод.

[Ca28] E. Cartan, G\'eom\'etrie des espaces de Riemann, Paris, 1928.
Рус. перевод: Э. Картан, Геометрия римановых пространств, Ленинград, 1936.

[Gr90] A. Gray, Tubes. Addison-Wesley, 1990.
Рус. перевод: А. Грей, Трубки, Наука, Москва, 1997.

[Gr94] M. Gromov, Sign and geometric meaning of curvature, Rend. Sem. Mat. Fis.
Milano 61 (1991), 9-123 (1994).
Рус. перевод: М. Громов, Знак и геометрический смысл кривизны, НИЦ 'Регулярная
и хаотическая динамика', Ижевск, 2000.

[GV00] Н.Б. Васильев, В.Л. Гутенмахер, Прямые и кривые, Москва, МЦНМО, 2000,
http://www.math.ru/prkr

[MSF04] А. С. Мищенко, Ю. П. Соловьев и А. Т. Фоменко, Сборник задач по
дифференциальной геометрии и топологии, Москва, Физматлит, 2004.

[MF04] А. С. Мищенко и А. Т. Фоменко, Краткий курс дифференциальной геометрии и
топологии, Москва, Физматлит, 2004.



[Pr] В. В. Прасолов, Курс дифференциальной геометрии, готовится к публикации.

[PS] В. Прасолов и А. Скопенков, Размышления о признании геометрии Лобачевского,
представлено к публикации, http://arxiv.org/abs/1307.4902

[Ra03] П. К. Рашевский, Курс дифференциальной геометрии, Москва, УРСС, 2003.

[Ra04] П. К. Рашевский, Риманова геометрия и тензорный анализ, Москва, УРСС,
2004.

[Sk10] А. Скопенков, Объемлемая однородность, Москва, МЦНМО, 2012,
\linebreak
http://arxiv.org/abs/1003.5278



[Ta89] С. Л. Табачников, О кривизне, Квант, 1989, N5;
Дифференциальная геометрия вокруг нас, Квант, 1989, N11.
http://kvant.mirror0.mccme.ru/

[To82] Дж. Торп, Начальные главы дифференциальной геометрии, Москва, Мир, 1982.

[Vi92] Н. Я. Виленкин, О кривизне, Квант, 1992, N4.
http://kvant.mirror0.mccme.ru/


\newpage
\section{Кривизны кривых}

\subsection{Кривые}\label{0curves}

Будем обозначать точкой производную по $t$, а штрихом производную по
{\it натуральному параметру} (когда это понятие появится).

\smallskip
{\bf 0.} Нарисуйте приближенно следующие траектории и кривые на плоскости или в пространстве.
Найдите их уравнения $r(t)=(x(t),y(t))$ или $r(t)=(r(t),\varphi(t))$ в декартовых или полярных координатах на плоскости; $r(t)=(x(t),y(t),z(t))$ или $r(t)=(r(t),\varphi(t),z(t))$ или $r(t)=(r(t),\varphi(t),\theta(t))$ в декартовых, цилиндрических или сферических координатах в пространстве.
Систему координат выберите сами.
Все скорости в этой задаче предполагаются ненулевыми.

(0) {\it Парабола} --- траектория мячика, брошенного со скоростью, параллельной земле, движущийся под действием силы тяжести (без учета сопротивления воздуха etc.),
или множество точек плоскости, равноудаленных от данной прямой ({\it директрисы}) и данной точки ({\it фокуса}).

(a) Луч $OA$ равномерно вращается вокруг своего неподвижного начала $O$ с угловой скоростью $\omega$.
Точка $M$ равномерно движется по лучу $OA$, начиная из точки $O$, со скоростью $v$.
Описываемая точкой $M$ траектория называется {\it спиралью Архимеда}.

(b) {\it Винтовая линия} --- траектория конца стержня длины $2r$,
равномерно со скоростью $v$ падающего на землю, остающегося
параллельным поверхности земли
и одновременно вращающегося в горизонтальной плоскости вокруг
своей середины равномерно с угловой скоростью $\omega$.

(c) Колесо радиуса $R$ катится равномерно без проскальзывания по прямой.
Описываемая точкой на ободе колеса траектория называется {\it циклоидой}.

(d) {\it Эллипс} --- множество точек плоскости, сумма расстояний от которых до
двух данных точек ({\it фокусов}) равна фиксированной величине $d$, большей расстояния $f$ между фокусами.

(e) По какой траектории движется электрон в постоянном магнитном поле, если начальная скорость электрона
не параллельна и не перпендикулярна напряженности $H$, где $H$ --- постоянный вектор?
(Закон Био-Савара-Лапласа движения электрона утверждает, что $\ddot r=\dot r\times H$.)

(f) {\it Логарифмическая спираль} (или изогональная спираль) --- плоская кривая, касательная в каждой точке $X$ которой образует с вектором $OX$ один и тот же угол $\alpha$.
Здесь $O$ и $\alpha$ --- наперед заданные точка плоскости и угол.

Интересно, что этот вид спирали часто встречается в природе.
Например, в растущих формах, подобных раковинам моллюсков, шляпкам подсолнечников, спиралям циклонов и галактик.
По этой кривой бабочка в ночи движется к лампе.

(g)* {\it Цепная линия} --- кривая, форму которой под действием
силы тяжести принимает нерастяжимая нить с закрепленными концами.

(h) {\it Кривая Вивиани} --- пересечение сферы радиуса $R$ и прямого кругового
цилиндра диаметра $R$, одна из образующих которого проходит через центр сферы.

(i)* {\it Астроида} --- кривая, для которой длина отрезка касательной
в произвольной точке, заключенного между осями координат, постоянна и равна $a$.

(j) Окружность радиуса $R$ катится без проскальзывания снаружи по окружности того же радиуса $R$.
Траектория, описываемая точкой на внешней окружности, называется {\it кардиоидой}.

\smallskip
{\bf 1.} Дайте `определение' отображения.
Дайте определение равных отображений, образа множества при отображении, прообраза множества при отображении, инъекции
(т.е. взаимно однозначного отображения), сюръекции (т.е. отображения на) и биекции (т.е. взаимно однозначного отображения на, или взаимно однозначного соответствия).

\smallskip
Пусть $D=[a,b]$ или $D=[a,+\infty)$ или $D=\R$.

{\bf Параметризованной гладкой регулярной кривой} на плоскости называется бесконечно дифференцируемое отображение $r:D\to\R^2$ (или, что то же самое, упорядоченная пара отображений $x,y:D\to\R$), производная (скорость)
$\dot r(t)$ которого ненулевая в любой точке $t\in [a,b]$.

{\bf Непараметризованной гладкой кривой} на плоскости называется
подмножество $\Pi\subset\R^2$, для любой точки $P\in\Pi$ которого существует такая ее замкнутая окрестность $OP$ в $\R^3$, что $\Pi\cap OP$ является образом $r[0,1]$ некоторой инъективной параметризованной гладкой регулярной кривой $r:[0,1]\to\R^2$.
\footnote{Таким образом, для непараметризованных кривых свойство `регулярности' включается в понятие `гладкости'.
Это мотивировано тем, что непараметризованные кривые --- частный случай подмногообразий, а не отображений.
Частным случаем отображений являются параметризованные кривые.}

Далее прилагательные `гладкая регулярная' и `гладкая' опускаются.
(Иногда непараметризованную кривую называют {\it кривой} или {\it траекторией}.)

Непараметризованные кривые интересны с точки зрения геометрии.
С этой точки зрения параметризованные кривые --- средство изучения непараметризованных.
Кроме того, параметризованные кривые интересны с точки зрения физики.

{\it Параметрическим уравнением} или {\it параметризацией} непараметризованной кривой $\Gamma\subset\R^2$ называется
такая параметризованная кривая $r:D\to\R^2$, что $\Gamma=r(D)$.

\smallskip
{\bf 2.} (a) Приведите пример не взаимно однозначной параметризации окружности.

(b) Приведите пример двух разных взаимно однозначных параметризаций одной непараметризованной кривой.

(d)* Постройте бесконечно дифференцируемое отображение $r:[-1,1]\to\R^2$, образом которого является
прямой угол (т.е. объединение отрезков $0\times[0,1]$ и $[0,1]\times0$).

\smallskip
{\bf Длиной} параметризованной кривой $r:[a,b]\to\R^2$ называется число
$$L(r):=\sup\{|A_0A_1|+|A_1A_2|+\dots+|A_{m-1}A_m|\ :\ m\in\N,\ a=a_0\le a_1\le\dots\le a_m=b,\ A_k:=r(a_k)\}.$$

{\bf 3.} (a) Длина графика дважды дифференцируемой функции $y:[a,b]\to\R$ равна
$\int_a^b\sqrt{1+\dot y(t)^2}dt$.

(b) {\bf Теорема.}
Длина параметризованной кривой $r=(x,y):[a,b]\to\R^2$ равна
$$L(r)=\int
_a^b\sqrt{\dot x(t)^2+\dot y(t)^2}dt=\int
_a^b|\dot r(t)|dt.$$

{\bf 4.} Вычислите длины параметризованных кривых от $r(a)$ до $r(b)$
для выбранных Вами параметризаций (укажите параметризацию явно!)
\quad

(a) винтовой линии; \quad (b) параболы;  \quad
(c) спирали Архимеда;  \quad (d) циклоиды.

\smallskip
{\bf 5.} Найдите формулу для длины параметризованной кривой
$r:[a,b]\to\R^2$ в
\quad

(a) полярных; \quad (b) сферических; \quad (c) цилиндрических

координатах.


\smallskip
{\bf Длиной} непараметризованной кривой называется длина некоторой ее взаимно однозначной параметризации.

\smallskip
{\bf 6.} (a) Если $r:\R\to\R^2$ --- взаимно однозначная параметризация параболы $y=x^2$ и
$r(0)=(0,0)$, $r(2)=(2,4)$, то $0<r^{-1}(1,1)<2$.

(b) Если $r_1,r_2:[a,b]\to\R^2$ --- взаимно однозначные параметризации одной непараметризованной кривой, то отображение $r_1r_2^{-1}:[a,b]\to[a,b]$ взаимно однозначно и монотонно.

(с) {\bf Теорема.} Приведенное определение длины корректно, т.е. если $r_1$ и $r_2$ --- две взаимно однозначные параметризации одной непараметризованной кривой, то $L(r_1)=L(r_2)$.

{\it Предостережение:} если при доказательстве Вы используете формулу для длины параметризации, то не забудьте доказать, что отображение $r_1r_2^{-1}:[a,b]\to[a,b]$ имеет в каждой точке положительную производную.
Впрочем,  проще доказывать {\it по определению} длины параметризованной кривой и используя (b).

\smallskip
Аналогично определяются кривые и их длины в пространстве.
В вышеприведенных определениях и теоремах нужно заменить $\R^2$ на $\R^3$ (или на $\R^n$).

\subsection{Кривизна кривых}\label{0curture}

В задачах с явным физическим содержанием следует пренебрегать размерами
движущихся объектов (т.е. считать их материальными точками) и т. п.

\smallskip
{\bf 0.} Если велосипедист движется по криволинейной дороге с постоянной по
модулю скоростью, то в любой точке его ускорение перпендикулярно его скорости
(т.е. их скалярное произведение равно нулю).

\smallskip
{\bf 1.} (a) Велотрек имеет форму винтовой линии c параметрами $r=10$м, $v=100$м/с и $\omega=2\pi c^{-1}$,
ось которой параллельна поверхности Земли (см. задачу \ref{0curves}.0.b).
Велосипедист едет по ней с постоянной по модулю скоростью и во время движения находится `со стороны' оси
винтовой линии от велотрека.
При каком значении модуля скорости он не оторвется от велотрека в `верхней' точке велотрека
(иными словами, модуль ускорения велосипедиста в `верхней' точке велотрека равен $g$)?

{\it Предостережение:} уравнение винтовой линии не обязательно будет уравнением движения велосипедиста.

{\it Предостережение:} эта задача не претендует на адекватность математической модели и реальной ситуации.

(b) Американская горка имеет форму циклоиды, находящейся в вертикальной плоскости.
Вагончик движется по ней со скоростью 1 м/с.
При какой высоте циклоиды в ее верхней точке клиент будет чувствовать невесомость?
Почему нереалистичен полученный Вами ответ?

(c) Автомобиль едет по отрезку спирали Архимеда $r=\varphi\cdot1$м со скоростью
1 м/с, не пересекая реку, т.е. луч $\varphi=0$.
С какой угловой скоростью вращается берег реки в системе отсчета автомобиля,
когда тот находится в точке с $\varphi=\pi/2$?

(d) Дан эллипс с параметрами $d=2$ и $f=1$.
Найдите радиус соприкасающейся окружности в одной из двух точек эллипса, равноудаленных от фокусов.

Окружность называется {\it соприкасающейся} с кривой в точке $A$ кривой, если
эта окружность лучше всех окружностей приближает кривую, т.е. если эта
окружность проходит через точку $A$, имеет общую с кривой касательную  в точке
$A$ и общее с кривой ускорение в точке $A$ при натуральной (см. далее) параметризации обеих кривых.

(Если в некоторой декартовой системе координат дуга окружности, содержащая
точку $A$, имеет уравнение $y=g(x)$, и дуга кривой, содержащая точку $A$, имеет
уравнение $y=f(x)$, то условие соприкасания равносильно тому, что для
некоторого $x_0\in\R$ выполнено
$A=(x_0,f(x_0))=(x_0,g(x_0))$, $f'(x_0)=g'(x_0)$, $f''(x_0)=g''(x_0)$.)


\smallskip
{\it Скоростью} (или {\it производной}) параметризованной кривой $r=(x,y):[a,b]\to\R^2$ в точке $t\in[a,b]$
называется вектор $\dot r(t)=(\dot x(t),\dot y(t))$.

Параметризованная кривая называется {\it натуральной}, если модуль ее скорости равен 1 в любой точке.
По задаче \ref{0curves}.3.a это эквивалентно тому, что ее длина от 0 до $s$ равна $s$ для любого $s$.
Натуральную параметризованную кривую будем обозначать через $\sigma$ (а не $r$).

{\it Ускорением} (или {\it второй производной}) параметризованной кривой $r=(x,y):[a,b]\to\R^2$ в точке
$t\in[a,b]$ называется вектор $\ddot r(t)=(\ddot x(t),\ddot y(t))$.

\smallskip
{\bf 2.} (a) Для натуральной кривой ускорение перпендикулярно скорости в любой точке.

(b) Найдите натуральную параметризацию винтовой линии в виде $\sigma(s)=(x(s),y(s),z(s))$.

(c) То же для циклоиды $(t-\cos t,\sin t)$, $t\in[-\pi/2,\pi/2]$, в виде $\sigma(s)=(x(s),y(s))$.

(c') Существует ли натуральная параметризация циклоиды $(t-\cos t,\sin t)$, $t\in[0,\pi]$?


(d) Непараметризованная кривая, имеющая взаимно однозначную регулярную параметризацию, имеет единственную
(с точностью до сдвига и симметрии) натуральную параметризацию.

(e) Приведите пример непараметризованной кривой с одной `точкой самопересечения',
имеющей не единственную (даже с точностью до сдвига и симметрии) натуральную параметризацию.

\begin{figure}[h]\centering
\includegraphics{p.5}
\caption{Кривизна кривой}\label{f:curture}
\end{figure}


\smallskip
Зафиксируем ориентацию плоскости.
{\bf Кривизной} в точке $s_0$ плоской натуральной параметризованной кривой $\sigma$ называется
число $k(s_0)$, равное по модулю числу $|\sigma''(s_0)|$ и совпадающее
с ним (или противоположное ему) по знаку, если векторы $\sigma'(s_0)$ и
$\sigma''(s_0)$ образуют положительный (отрицательный) базис.
Иными словами, $k(s_0):=\sigma'(s_0)\wedge\sigma''(s_0)$.
{\bf Кривизной} в точке $A=\sigma(s_0)$ непараметризованной кривой $\Gamma$
с натуральной параметризацией $\sigma$ называется число $k(A):=k(s_0)$.

\smallskip
{\bf 3.} (a) Кривизна натуральной параметризованной кривой равна угловой скорости вращения вектора скорости.

(b) {\bf Теорема.} Кривизна непараметризованной кривой, имеющей параметризацию $r$, в точке $r(t)$ равна отношению проекции ускорения $\ddot r(t)$ на прямую, перпендикулярную скорости $\dot r(t)$, к квадрату модуля скорости:
$$k(r(t))=\dfrac{\dot r\wedge\ddot r}{|\dot r|^3}=\dfrac{\ddot x\dot y-\dot x\ddot y}{|\dot x^2+\dot y^2|^{3/2}},$$
где аргумент $t$ функций в правых частях пропущен.

\smallskip
Два подмножества пространства $\R^n$ называются {\it объемлемо изометричными}, если между ними существует
{\it объемлемая изометрия} (движение), т.е. сохраняющее расстояния (в $\R^n$!) отображение $\R^n\to\R^n$,
переводящее первое подмножество во второе.

\smallskip
{\bf 4.} (a) {\bf Теорема.} Для любой функции $\overline k:\R\to\R$ существует
и единственна (с точностью до объемлемой изометрии) плоская натуральная
параметризованная кривая $\sigma:\R\to\R^2$, для которой $k(\sigma(s))=\overline k(s)$.

(b) {\bf Следствие.} Две
плоские несамопересекающиеся непараметризованные кривые объемлемо изометричны тогда и только тогда, когда
функции кривизны их натуральных параметризаций отличаются сдвигом и, возможно,
знаком, т.е. когда существует число $a$, для которого либо
$\sigma_2(s)=\sigma_1(s+a)$ при любом $s$, либо $\sigma_2(s)=\sigma_1(-s+a)$ при любом $s$.

(c) То же, что в (b), для замкнутых кривых.

(d)* Функция $\overline k:S^1\to\R$ кривизны замкнутой плоской параметризованной
кривой $r:S^1\to\R^2$ имеет по крайней мере четыре (нестрогих) экстремума.

(e)* Для любой функции $\overline k:S^1\to\R$, имеющей по крайней мере четыре
(нестрогих) экстремума, существует и единственна (с точностью до объемлемой
изометрии) замкнутая плоская натуральная параметризованная кривая
$\sigma:S^1\to\R^2$, для которой $k(\sigma(s))=\overline k(s)$.

(f) Какие гладкие плоские несамопересекающиеся кривые $N$ имеют следующее
свойство: для любых двух точек $x,y\in N$ существует
движение плоскости, переводящее $x$ в $y$, а $N$ в себя? Ср. [Sk10].

(g) То же для пространственного случая.

\smallskip
Скорость параметризованной кривой и натуральная параметризация определяются для
пространства $\R^3$ (или $\R^n$) аналогично случаю плоскости.

{\bf Кривизной} в точке $s_0$ пространственной натуральной параметризованной
кривой $\sigma$ называется число $k(s_0):=|\sigma''(s_0)|$.
Кривизна пространственной непараметризованной кривой определяется дословно так же, как плоской.

Заметим, что кривизна кривой, рассматриваемой как кривая в плоскости, может
отличаться знаком от кривизны той же кривой, рассматриваемой в пространстве.

\smallskip
{\bf 5.} $k(r(t))=|\dot r\times\ddot r|/|\dot r|^3$ для пространственной параметризованной кривой $r:D\to\R^3$.

\smallskip
{\bf 6.} Найдите кривизну в точке $(1,1,1)$ кривой пересечения поверхностей $z=xy$ и $y^2+z^2=2$.


\subsection{Кручение пространственных кривых}\label{0torsion}

Угловая скорость вращения плоскости равна угловой скорости вращения вектора, нормального к этой плоскости.
Угловая скорость вращения вектора $a(t)$ равна $\left|\left(\dfrac{a(t)}{|a(t)|}\right)'_t\right|$.
(Эти фразы можно считать определениями.)

\smallskip
{\bf 1.} В невесомости мотоциклист едет со скоростью 1 по винтовой линии c параметрами $r$, $v$ и $\omega$
(см. задачу \ref{0curves}.0.b).
Найдите угловую скорость $\Omega(t)$ вращения плоскости колеса (содержащей векторы скорости и ускорения)
в зависимости от времени.

\smallskip
Параметризованная кривая $r:D\to\R^m$ называется {\it бирегулярной}, если векторы $\dot r(t)$ и $\ddot r(t)$ линейно независимы для любого $t$.
Непараметризованная кривая $r$ называется {\it бирегулярной}, если она имеет бирегулярную параметризацию.

\smallskip
{\bf 2.} (a) Приведите пример (гладкой регулярной) непараметризованной кривой, не являющейся бирегулярной.

(b) Две параметризованные кривые с одинаковым образом бирегулярны или нет одновременно.

(c) Регулярная непараметризованная кривая бирегулярна тогда и только тогда, когда ее кривизна отлична от нуля в любой точке.

\smallskip
{\bf 3.} (a) При каких $a$ образ параметризованной кривой $r(t)=(e^t,2e^{-t},e^{at})$ лежит в некоторой плоскости?

(b) Образ параметризованной бирегулярной кривой $r$ лежит в некоторой плоскости тогда и только тогда, когда
$\dot r(t)\wedge\ddot r(t)\wedge\dddot r(t)=0$ для любого $t$.

\smallskip
Пусть $\sigma$ --- бирегулярная натуральная параметризованная кривая.
В оставшейся части этого пункта производные берутся по ее параметру $s$, который всюду опускается.
Обозначим через
$$v:=\sigma',\quad n:=\sigma''/|\sigma''|\quad\text{и}\quad b:=v\times n$$
вектора {\it скорости}, и {\it главной нормали} и {\it бинормали}.
Они образуют {\it репер Френе}.

{\bf Кручением} в точке $s_0$ пространственной бирегулярной натуральной параметризованной кривой $\sigma$ называется
число $-b'\cdot n$.
Это угловая скорость (со знаком) вращения ориентированной соприкасающейся плоскости.
{\it Соприкасающаяся} плоскость --- плоскость, содержащая векторы скорости $\sigma'=v$ и ускорения $\sigma''=kn$;
она ориентирована репером $v,n$.


\smallskip
{\bf 4.} (a) Напишите определение кручения непараметризованной бирегулярной кривой.

(b) Найдите кручение в точке $(-1/2,-1/2,1/4)$ кривой пересечения поверхностей $z=xy$ и $z=(x+1)(y+1)$.

\smallskip
{\bf 5.}
(a) $\sigma''\ne0$ (т.е. определение вектора нормали корректно).

(b) $n'=-kv$ для плоского случая.
Эта формула вместе с формулой $v'=kn$ называются {\it формулами Френе на плоскости.}

(e) $b'=-\varkappa n$.

(c) $\varkappa$ равно проекции вектора $n'$ на направленную ось,
перпендикулярную векторам скорости и ускорения: 
$\varkappa=n'\cdot(v\times n)=v\wedge n\wedge n'$.

(d) $n'=-kv+\varkappa b$.

Формулы из пп. (d) и (e) вместе с формулой $v'=kn$ называются {\it формулами Френе в пространстве}.

\smallskip
{\bf 6.} (a) $\varkappa=\dfrac{\sigma'\wedge\sigma''\wedge\sigma'''}{k^2}$
для бирегулярной натуральной параметризованной кривой $\sigma$ (аргумент $s$ функций в правой части пропущен).

(b) $\varkappa(r(t))=\dfrac{\dot r\wedge\ddot r\wedge\dddot r}{|\dot r\times\ddot r|^2}$
для бирегулярной параметризованной кривой $r:D\to\R^3$
(не обязательно натуральной; аргумент $t$ функций в правой части пропущен).

(c) {\bf Теорема.}  Для любых функций $\overline k:\R\to(0,+\infty)$ и $\overline\varkappa:\R\to\R$ существует и
единственна (с точностью до объемлемой изометрии) пространственная бирегулярная натуральная параметризованная кривая
$\sigma:\R\to\R^3$, для которой $k(\sigma(s))=\overline k(s)$ и $\varkappa(\sigma(s))=\overline \varkappa(s)$.

(d) Сформулируйте и докажите аналог следствия \ref{0curture}.4.b для пространственных кривых.

\newpage
\subsection*{Указания и решения к некоторым задачам}
\addcontentsline{toc}{subsection}{Указания и решения к некоторым задачам}

{\bf \ref{0curves}.0.} (f) $r'(\varphi)=r(\varphi)\ctg\alpha$.

\smallskip
{\bf \ref{0curves}.3.}
(a) Используйте теорему Лагранжа ($y(x_2)-y(x_1)=y'(\theta)(x_2-x_1)$).
Не забудьте аккуратно обосновать переход от супремума к пределу интегральных сумм.

(b) Используйте теорему Коши ($\dfrac{y(t_2)-y(t_1)}{x(t_2)-x(t_1)}=\dfrac{y'(\theta)}{x'(\theta)}$).

\smallskip
{\bf \ref{0curves}.6.} (a) Аналогично (b).

(b) Обозначим $\Gamma:=r_1[a,b]$.
Для $k=1,2$ определим отображение $\overline r_k:[a,b]\to\Gamma$ формулой $\overline r_k(t):=r_k(t)$.
Так как $\overline r_1$ и $\overline r_2$ биекции, то и $\overline r_1\overline r_2^{-1}$ биекция.
Так как $\overline r_1$ и $\overline r_2$ непрерывны (определите, что это такое! при этом используйте расстояния `по плоскости', а не `по кривой $\Gamma$'), то и $\overline r_1\overline r_2^{-1}$ непрерывно (докажите!).
Значит, по теореме о промежуточном значении непрерывной функции $r_1r_2^{-1}=\overline r_1\overline r_2^{-1}$ есть монотонная биекция.

\bigskip
{\bf \ref{0curture}.1.}
(a) Для уравнения винтовой линии из задачи \ref{0curves}.0.b модуль скорости постоянен.

(b) Можно использовать результат задачи \ref{0curture}.2.c, но искать вторую производную в явном виде громоздко.
Лучше сначала доказать формулу \ref{0curture}.3.b.

(c) См. \ref{0curture}.3.a.

(d) Докажите, что этот радиус равен $1/k$, где $k$ --- определенная ниже кривизна.
Для подсчета кривизны докажите формулу \ref{0curture}.3.b

\smallskip
{\bf \ref{0curture}.2.}
(с) Пусть уравнение циклоиды $r(t)=(t-\sin t,\cos t)$.
Вычислите $l(t):=\int_0^t|\dot r(t)|dt$.
Найдите натуральную параметризацию $\sigma$ из уравнения $\sigma(l(t))=r(t)$.

(c') Нет.

(d) Пусть $r:[a,b]\to\R^2$ --- взаимно однозначная параметризация заданной непараметризованной кривой.
Для $t\in [a,b]$ обозначим $l(t):=\int_0^t|\dot r(t)|dt$.
Так как параметризация  $r(t)$  регулярная, то функция $l$ строго монотонна.
Докажите, что параметризация $\sigma$ непараметризованной кривой $r[a,b]$ является натуральной тогда и только тогда, когда $\sigma(l(t))=r(t)$ для любого $t\in [a,b]$.
Тогда существование и единственность натуральной параметризации будут вытекать из строгой монотонности функции $l$.

(e) Пример можно построить в виде `восьмерки'.

\smallskip
{\bf \ref{0curture}.3.} (a)
Пусть $\sigma$ --- наша кривая.
Обозначим $\alpha(s):=\angle(\sigma'(s),Ox)$.
Тогда
\linebreak
$\sigma'(s)=(\cos\alpha(s),\sin\alpha(s))$.
Значит, $\sigma''(s)=(-\sin\alpha(s),\cos\alpha(s))\alpha'(s)$.
Поэтому $\alpha'(s)=k(\sigma(s))$.

(b) В этом решении производные вектор-функции $r$ берутся по
параметру $t$ в точке $t$, а производные вектор-функции $\sigma$ берутся
по параметру $s$ в точке $l(t):=\int_0^t|\dot r(t)|dt$.

Так как $|\sigma'|=1$, то $\dot r=\sigma'|\dot r|$.
(Интегрируя это соотношение, получаем $\sigma(l(t))=r(t)$.)
Дифференцируя это соотношение и используя $\dot l=|\dot r|$, получаем $\ddot r=\sigma''|\dot r|^2+\sigma'|\dot r|'_t$.
Из этого и $\sigma''\perp\dot r$ вытекает $\ddot r\wedge\dot r=\sigma''|\dot r|^2\wedge\dot r=-k|\dot r|^3$.

\smallskip
{\bf \ref{0curture}.4.} (a)
Пусть $\sigma$ --- искомая кривая.
Возьмем систему координат $Oxy$, для которой $\sigma(0)=(0,0)$ и
ось $Ox$ сонаправлена с $\sigma'(0)$.
Определим функции $x,y:\R\to\R$ равенством $\sigma(s)=(x(s),y(s))$.
Обозначим $\alpha(s):=\angle(\sigma'(s),Ox)$.
Тогда $\sigma'(s)=(x'(s),y'(s))=(\cos\alpha(s),\sin\alpha(s))$.
По здаче \ref{0curture}.3.a $\alpha'(s)=k(\sigma(s))$.
Поэтому для искомая кривая однозначно определяется соотношениями
$$\alpha(s)=\int_0^sk(s)ds,\quad x(s)=\int_0^s\cos\alpha(s)ds\quad\mbox{и}
\quad y(s)=\int_0^s\sin\alpha(s)ds.$$

{\bf \ref{0curture}.6.} Тригонометрическая параметризация удобна для вычислений, а радикальная --- нет.

\bigskip
{\bf \ref{0torsion}.1.}
 Если плоскость содержит векторы скорости и ускорения, то нормальный к ней вектор равен векторному произведению
векторов скорости и ускорения.
Уравнение винтовой линии $r(t)=(r\cos\omega t,r\sin\omega t,vt)$.
Обозначим
$\overline\omega:=\dfrac{\omega}{|\dot r(t)|}=\dfrac{\omega}{\sqrt{r^2\omega^2+v^2}}$.
Тогда уравнение движения мотоцикла $\sigma(s)=
(r\cos\overline\omega s,r\sin\overline\omega s,\dfrac{vs\overline\omega}\omega)$.
Находим
$$\sigma'(s)=(r\overline\omega\sin\overline\omega s,
r\overline\omega\cos\overline\omega s,\dfrac{vs\overline\omega}\omega), \quad
\sigma''(s)=(-r\overline\omega^2\cos\overline\omega s,
r\overline\omega^2\sin\overline\omega s,0).$$
Отсюда получаем ответ $\Omega(t)=\dfrac{v\omega}{\sqrt{v^2+r^2\omega^2}}$.

\smallskip
{\bf \ref{0torsion}.2.} (a) Прямая.

(b) Условие бирегулярности локальное.
Поэтому и ввиду регулярности можно считать, что обе кривые взаимно однозначны (инъективны).
Замена параметра меняет только касательную составляющую вектора скорости.

\smallskip
{\bf \ref{0torsion}.4.} (a) {\it Кручением} в точке $A=\sigma(s_0)$ непараметризованной кривой $\Gamma$
с натуральной параметризацией $\sigma$ называется число $\varkappa(A):=\varkappa(s_0)$.
{\it Кручением} в точке $A$ пространственной бирегулярной непараметризованной кривой $\Gamma$ называется кручение  $\varkappa(0)$ такой натуральной параметризации $\sigma$ кривой $\Gamma$, для которой $\sigma(0)=A$.

(b) 0, ибо кривая лежит в плоскости. 

\smallskip
{\bf \ref{0torsion}.5.} (b) Угловая скорость вращения вектора $n$ равна угловой
скорости вращения вектора $v$, поскольку эти векторы перпендикулярны.

(e) Так как $|b|=1$, то $b'\perp b$.
Имеем $b'\cdot v=-(v'\times n)\cdot v-(v\times n')\cdot v=0$.
Значит, $b'\perp v$.
Поэтому $b'||n$.
Отсюда $b'=(b'\cdot n)n=-\varkappa n$.

(c) $\varkappa=-b'\cdot n=-(v'\times n)\cdot n-(v\times n')\cdot n=v\wedge n\wedge n'$.

(d) Так как $|n|=1$, то $n'\perp n$.
Поэтому $n'=(n'\cdot v)v+(n'\cdot b)b=kv+\varkappa b$ по (c).

\smallskip
{\bf \ref{0torsion}.6.} (b) В этом решении производные функции $r$ берутся по параметру $t$ в точке $t$, а
производные вектор-функций $\sigma$ и $b$ берутся по параметру $s$ в точке $l(t):=\int_0^t|\dot r(t)|dt$.

Так как плоскости, образованные парами векторов $(\sigma',\sigma'')$ и $(\dot r,\ddot r)$, совпадают, то
$\dot r\times\ddot r=\nu b$, где $\nu:=|\dot r\times\ddot r|$.
Дифференцируя по $t$, получаем $\dot r\times\dddot r=\nu b'|\dot r|+\nu'_tb$.
Умножая скалярно на $\ddot r$ и используя \ref{0torsion}.5.e, получаем
$$\dot r\wedge\ddot r\wedge\dddot r=\nu|\dot r|b'\cdot\ddot r=\frac{\nu|\dot r|\varkappa\nu}{|\dot r|}=
\nu^2\varkappa.$$

(с) Для заданных функций $\overline k$ и $\overline\varkappa$ тройка векторов $(v,n,b)$ однозначно определена ввиду формул Френе $v'=kn$, $n'=-kv+\varkappa b$ и $b'=-\varkappa n$ (\ref{0torsion}.5.de).

\newpage
\section{Элементы геометрии поверхностей}\label{0sur}

\subsection{Элементы сферической геометрии}\label{0sursph}

{\bf 1.} (a) Человек прошел 10 километров на юг, потом 10 километров на восток и потом 10 километров на север.
В итоге он вернулся в исходное положение.
Найдите все начальные точки, при которых такое могло быть.

(b) Если все время идти на северо-восток, то куда придешь?

(c) {\it Локсодромия} --- траектория путешественника, движущегося по поверхности Земли (которая считается сферой)
все время под углом $\alpha$ к текущему меридиану.
Здесь $\alpha$ --- наперед заданный угол.
Напишите уравнение локсодромии в сферических координатах, если широта путешественника возрастает равномерно.

Интересно, что на работе с такими кривыми основана меркаторовская картографическая проекция.
Ср. с задачей \ref{0curves}.0.f.

\smallskip
Обозначим
$$S^2:=\{(x,y,z)\ |\ x^2+y^2+z^2=1\}.$$
{\it Расстоянием по сфере} между точками сферы называется длина наименьшей дуги большого круга, соединяющей эти точки.
{\it Круг} и {\it окружность} на сфере определяются аналогично случаю плоскости (через расстояние по сфере).

{\it Стереографической проекцией} называется отображение из $S^2-\{(0,0,1)\}$
в плоскость $z=0$, являющееся центральной проекцией из точки $(0,0,1)$.

\smallskip
{\bf 2.} (a) При стереографической проекции окружности переходят в окружности.

(b) В какую точку плоскости переходит при стереографической проекции центр той сферической окружности,
которая переходит в окружность $x^2+(y-2)^2=1$ при стереографической проекции?

\smallskip
{\bf 3.} (a) Грузы масс $m_1, m_2$ на сфере соединены невесомым отрезком сферической прямой.
В какой точке надо закрепить этот отрезок, чтобы полученные `весы' находились в равновесии?
Грузы притягиваются к центру сферы с силами, пропорциональными массам.

(b)* Определите центр масс системы $n$ точек на сфере и докажите теорему о группировке (начните с $n=3$).

\smallskip
Из произвольной точки $M$, лежащей внутри данного трехгранного угла $OABC$, опустим перпендикуляры $MA_1$, $MB_1$ и $MC_1$ на грани $OBC$, $OAC$ и $OAB$, соответственно.
{\it Двойственным} (или полярным) к трехгранному углу $OABC$ называется трехгранный угол $MA_1B_1C_1$ с вершиной $M$.

\smallskip
{\bf 4.} (a) Длина окружности радиуса $R$ на сфере равна $2\pi\sin R$.

(b) Любой трехгранный угол является двойственным к своему двойственному, т.е.
двукратное применение операции построения двойственного угла приводит к исходному трехгранному углу.

(c) Если все двугранные углы трехгранного угла прямые, то и все плоские углы тоже прямые.

(d) Выразите плоские углы двойственного угла через двугранные углы исходного и наоборот.

\smallskip
{\bf 5.} (a) Выполнены ли три признака равенства треугольников для сферических треугольников?

(b) Сферические треугольники равны `по трем углам'.

(c) Если в сферическом треугольнике две стороны равны, то равны и противолежащие им углы (и обратно).

\smallskip
Стороны и противолежащие им углы сферического треугольника обозначаются $\alpha ,\beta ,\gamma$ и $A,B,C$ соответственно.
Там, где это осмысленно, нужно выяснить и доказать, в какое утверждение плоской евклидовой геометрии переходит доказываемая теоpема в пределе при $\alpha^2+\beta^2+\gamma^2\to 0$.

\smallskip
{\bf 6.} (a) $\alpha+\beta\ge\gamma$;  \quad
(b) $\alpha\ge|\beta -\gamma|$; \quad (c) $\alpha+\beta+\gamma \le2\pi$;

(A) $\pi+C\ge A+B$; \quad (B) $C+|A-B|\le\pi$; \quad (C) $A+B+C\ge\pi$.

\smallskip
{\bf 7.} Сумма углов, образованных диагональю прямоугольного параллелепипеда с его ребрами, меньше $\pi$.

\smallskip
{\bf 8.} (a) В тетраэдре $SABC$ все углы при вершине $C$ прямые, а  $\angle BSC=\angle ASC=\pi/4$.
Найдите $\angle ASB$ и $\angle (CBS,ABS)$.

(b) {\it Прямоугольный треугольник.} Если $C=\pi/2$, то
 $\cos\gamma=\cos\alpha\cos\beta$ ({\it теорема Пифагора для сферы}),
$\sin\alpha=\sin\gamma\sin A$, \quad
$\tg\alpha=\tg\gamma\cos B=\tg A\sin\beta$,\quad $\cos B=\sin A\cos\beta$.

(c) {\it Теорема синусов.}
$\dfrac{\sin\alpha}{\sin A}=\dfrac{\sin\beta}{\sin B}=\dfrac{\sin\gamma}{\sin C}$.

(d) {\it Первая теорема косинусов.}
 $\cos\alpha = \cos\beta\cos\gamma+\sin\beta\sin\gamma\cos A$.

(e) {\it Вторая теорема косинусов.}
 $\cos A = -\cos B\cos C + \sin B\sin C\cos\alpha$.

(f)* Если $a,b,c$ углы между ребрами трехгранного угла и противоположными
 гранями, то $\sin\alpha\sin a = \sin\beta\sin b = \sin\gamma\sin c$.

\smallskip
{\bf 9.} В правильной $n$-угольной пирамиде (двугранный) угол при боковом  ребре равен $\varphi$.
Найдите:

($\alpha$) угол $\alpha$ при ребре основания; \quad

($\beta$) плоский угол $\beta$ при вершине;

($\gamma$)* угол $\gamma$ между боковым ребром и основанием.

\smallskip
{\bf 10.} В правильной 6-угольной пирамиде угол при боковом ребре равен $\varphi$.
 Найдите углы между каждой парой граней.

\smallskip
{\bf 11.} Если $AL$ биссектриса сферического треугольника $ABC$, то
 $\dfrac{\sin BL}{\sin BA} = \dfrac{\sin CL}{\sin CA}$.


\smallskip
{\bf 12.} В любом сферическом треугольнике в  одной  точке  пересекаются:

 (a) биссектрисы (т.е. существует вписанная окружность);

 (b) срединные перпендикуляры к сторонам (т.е. существует описанная окружность);

 (c) медианы; \quad

 (d) высоты.

\smallskip
{\bf 13.} Высоты в любом сферическом треугольнике существуют.
Всегда ли каждая высота единственна?
В связи с этим: всегда ли и в каком смысле справедливо утверждение пункта~(d) предыдущей задачи?

\smallskip
{\bf 14.} Сформулируйте и докажите сферические аналоги

(a) критериев вписанности и описанности четырехугольников.

(b) теорем Чевы и Менелая.

\smallskip
{\bf 15.*} Решите систему:
$$\ctg x\ctg y-5= \frac{\cos z}{\sin x\sin y};\quad
\ctg y\ctg z+11 = \frac{\cos x}{\sin y\sin z};\quad
\ctg z\ctg y+7  = \frac{cos y}{\sin x\sin z}.$$

\comment

\subsection{Элементы проективной геометрии}\label{0surpro}

Пусть $\alpha_1$ и $\alpha_2$ --- две плоскости в пространстве, $O$ --- точка,
не лежащая ни на одной из этих плоскостей.
{\it Центральным проектированием} с центром $O$ плоскости $\alpha_1$ на
плоскость $\alpha_2$ называют отображение, которое точке $A_1$ плоскости $\alpha_1$
ставит в соответствие точку пересечения прямой $OA_1$ с плоскостью $\alpha_2$.

\smallskip
{\bf 1.} (a) Центральная проекция плоскости $\alpha_1$
 на параллельную ей плоскость  $\alpha_2$ с центром $O$ равна некоторой гомотетии.

(b) Если плоскости $\alpha_1$ и $\alpha_2$ пересекаются, то центральное
проектирование $\alpha_1$ на $\alpha_2$ с центром $O$ зада\"eт
взаимно-однозначное отображение плоскости $\alpha_1$ без прямой $l_1$ на
плоскость $\alpha_2$ без прямой $l_2$, где $l_1$ и $l_2$ --- прямые пересечения
плоскостей $\alpha_1$ и $\alpha_2$  соответственно с плоскостями, проходящими
через $O$ параллельно $\alpha_1$ и $\alpha_2$.
 При этом на $l_1$ отображение не определено.

Прямую, на которой не определено центральное проектирование, а также прямую точек,
не имеющих прообраза, называют {\it исключительной прямой} данной проекции.

(b) Если прямые параллельны исключительное прямой, то их образы при этом проектировании также параллельны.

\smallskip
{\bf 2.} (a) Через точку на боковом ребре четыр\"ехугольной пирамиды проведите
 плоскость так, чтобы в сечении получился параллелограмм.

(b) Центральным проектированием можно перевести любой выпуклый четырехугольник
в параллелограмм.

\smallskip
{\bf 3.} (a) Пусть $O$ --- точка пересечения диагоналей выпуклого четыр\"ехугольника $ABCD$, а
 $E,F$ --- точки пересечения прямых $AB$ и $CD$, $BC$ и $AD$, соответственно.
 Прямая $EO$ пересекает стороны $AD$ и $BC$ в точках $M$ и $N$,
 прямая $FO$ пересекает стороны $AB$ и $CD$ в точках $P$ и $Q$, соответственно.
 Тогда прямые $PN$, $MQ$ и $EF$ пересекаются в одной точке.

(b) {\it Теорема Паппа.} Точки $C$ и $C'$ лежат на прямых $AB$ и $A'B'$, соответственно.
Тогда три точки пересечения прямых $AB'$ и $A'B$, $BC'$ и $B'C$, $CA'$ и $C'A$ лежат на одной прямой.

(c) {\it Теорема Дезарга.} Прямые $a$, $b$, $c$ пересекаются в одной точке $O$.
 В треугольниках $A_1B_1C_1$ и $A_2B_2C_2$ вершины $A_1$ и $A_2$ лежат на прямой $a$,
 $B_1$ и $B_2$ --- на прямой $b$, $C_1$ и $C_2$ --- на прямой $c$. $A$, $B$, $C$ - точки
 пересечения прямых $B_1C_1$ и $B_2C_2$, $A_1C_1$ и $A_2C_2$, $A_1B_1$ и $A_2B_2$ соответственно.
 Докажите, что точки $A$, $B$ и $C$ лежат на одной прямой.

(d) {\it Обратная теорема Дезарга.}
 В треугольниках $A_1B_1C_1$ и $A_2B_2C_2$ вершины $A_1$ и $A_2$ лежат на прямой $a$,
 $B_1$ и $B_2$ --- на прямой $b$, $C_1$ и $C_2$ --- на прямой $c$. $A$, $B$, $C$ - точки
 пересечения прямых $B_1C_1$ и $B_2C_2$, $A_1C_1$ и $A_2C_2$, $A_1B_1$ и $A_2B_2$ соответственно.
 Точки $A$, $B$ и $C$ лежат на одной прямой.
 Докажите, что прямые $a$, $b$, $c$ пересекаются в одной точке $O$.

\smallskip
{\bf 4.} Если точки $A, B, C, D$ лежат на прямой, параллельной исключительной прямой,
то для их образов $A',B',C',D'$ выполнено $\dfrac{A'B'}{C'D'}=\dfrac{AB}{CD}$.

\smallskip
{\bf 5.} (а) {\it Двойное\/} 
отношение $(ABCD):=\dfrac{\vec{AC}}{\vec{BC}}:\dfrac{\vec{AD}}{\vec{BD}}$
(упорядоченной четверки точек на прямой) сохраняется при центральном проектировании.

(b) Пусть $F$, $O$ --- точки пересечения прямых $AD$ и $BC$, $AC$ и $BD$, соответственно,
а $P$ и $Q$ --- точки пересечения прямой $FO$ c прямыми $AB$ и $CD$, соответственно.
 Тогда $\dfrac{\vec{FP}}{\vec{OP}}=-\dfrac{\vec{FQ}}{\vec{OQ}}$.

\smallskip
{\bf 7.} (a) Точки $a,b,c,d\in\C$ лежат на одной прямой или на одной окружности тогда и только
когда {\it двойное отношение} $(abcd):=\dfrac{c-a}{c-b}:\dfrac{d-a}{d-b}$ вещественно

(b) Напишите комплексное уравнение окружности, проходящей через точки $a,b,c\in\C$.

\smallskip
{\bf 6.} (а) $(ABCX)=(ABCY)\iff X=Y$.

(b) Любую упорядоченную тройку точек прямой можно перевести в любую другую композицией движения в плоскости и центрального или параллельного проектирования в плоскости.
Такое отображение прямой или прямой без точки единственно.

(c) Отображение $P$  прямой или прямой без точки в прямую является композицией движения в плоскости и центрального
или параллельного проектирования в плоскости тогда и только тогда, когда $P$ является дробно-линейным
(т.е. когда существуют  такие $a,b,c,d$, что $P(x)=\dfrac{ax+b}{cx+d}$).

(d) Отображение $P:\R^2-\{(0,y)\ :\ y\in\R\}\to \R^2$ плоскости без прямой в плоскость, заданное формулой  $P(x,y)=(1/x,y/x)$, является центральным проектированием.

\smallskip
Отображение плоскости $\alpha$ без некоторого набора прямых в плоскость $\beta$ называется {\it проективным отображением}, если оно является композицией движений, центральных проектирований и~параллельных проектирований
(т.е. если существуют плоскости $\alpha_1=\alpha, \alpha_2, \ldots,\alpha_n=\beta$ и~отображения~$P_i$ плоскостей
$\alpha_i$ и $\alpha_{i+1}$, каждое из которых является движением, центральных проектированием или параллельным проектированием, прич\"ем данное отображение является их композицией).
Если $\alpha=\beta$, то отображение будем называть {\it проективным преобразованием}.



{\it Проективной плоскостью} $\R P^2$ называется множество всех прямых в $\R^3$, проходящих через точку $(0,0,0)$.

Назовем точки $x,y\in\R^3-\{0\}$ {\it эквивалентными}, если $y=\lambda x$ для некоторого $\lambda\in(0,+\infty)$.
Проективная плоскость находится во взаимно-однозначном соответствии с множеством классов эквивалентности множества $\R^3-\{0\}$.
Класс эквивалентности точки $(x_1,x_2,x_3)$ записывается в виде $(x_1:x_2:x_3)$.
Числа $x_1,x_2,x_3$ называются {\it однородными координатами} точки $(x_1:x_2:x_3)\in\R P^2$.

Проективная плоскость находится во взаимно-однозначном соответствии с некоторой двумерной поверхностью в $\R^4$, см. задачу 2.1.4.c.

Пусть задана упорядоченная тройка $(P,A,B)\in\R^3$ точек, не лежащих на прямой, для которых плоскость $PAB$ не
проходит через точку $(0,0,0)$.
{\it Аффинной картой} проективной плоскости, соответствующей тройке $(P,A,B)$, называется отображение
$f=f_{(P,A,B)}:\R^2\to\R P^2$, которое точке $(x,y)\in\R^2$ ставит в соответствие прямую в $\R^3$,
проходящую через точки $(0,0,0)$ и $P+x\vec{PA}+y\vec{PB}$.
Для определения аффинной карты вместо точек $P,A,B\in\R^3$ можно задавать точку $P$ и пару неколлинеарных
векторов $a=\vec{PA}$, $b=\vec{PB}$.
Пара $(x,y)$ чисел называется {\it координатами точки $f(x,y)$ в аффинной карте}.


\smallskip
{\bf 8.} Выразите координаты в аффинной карте (проективной плоскости), соответствующей точке $Z:=(0,0,1)$ и
векторам $X:=(1,0,0)$, $Y:=(0,1,0)$, через координаты в аффинной карте, соответствующей точке $X$ и векторам $Y$, $Z$.

\smallskip
{\bf 9.} (a) Образ прямой при аффинной карте проективной плоскости (говорят: прямая на аффинной карте проективной плоскости) лежит в множестве, задаваемом в однородных координатах линейным однородным многочленом.
Ровно одна точка этого множества не лежит в образе прямой при аффинной карте.

(b) Сформулируйте и докажите аналог пункта (a) с заменой прямой на кривую второго порядка.

(c) Любое подмножество проективной плоскости, задаваемое в однородных координатах однородным многочленом второй степени,  есть либо образ окружности при аффинной карте, либо объединение двух
прямых, либо прямая,
либо точка, либо пустое множество.



\endcomment

\subsection{Определение и примеры поверхностей}\label{0surexa}

Под  {\it (непараметризованной) поверхностью} далее можно понимать поверхность вращения
(графика бесконечно дифференцируемой положительной функции $f:\R\to\R$) или график (бесконечно
дифференцируемой функции $f:\R^2\to\R$).

Приведем более общее определение.
Пусть $D$ --- декартово произведение интервалов (конечных или бесконечных) на прямой.

{\bf Гладкой регулярной параметризованной поверхностью} называется бесконечно дифференцируемое отображение
$r:D\to\R^3$ (или, что то же самое, упорядоченная тройка отображений $x,y,z:D\to\R$), производная которого невырождена в любой точке (т.е. пара векторов $r_u(u,v),r_v(u,v)$ линейно независима при любых $(u,v)\in D$).
Далее слова {\it гладкая регулярная} опускаются.

{\bf Гладкой непараметризованной поверхностью} называется подмножество $\Pi\subset\R^3$, для любой точки
$P\in\Pi$ которого существует такая ее замкнутая окрестность $OP$ в $\R^3$, что $\Pi\cap OP$
является образом $r([0,1]^2)$ инъективной параметризованной поверхности $r:[0,1]^2\to\R^3$.
Далее слова {\it гладкая непараметризованная} опускаются.

Ср. со сноской в \S\ref{0curves}.

{\it Параметрическим уравнением} или {\it параметризацией} поверхности $\Pi\subset\R^3$ называется такая параметризованная поверхность $r:D\to\R^3$, что $r(D)=\Pi$.

{\it Системой координат} с координатным пространством $D$ на поверхности $\Pi\subset\R^3$
называется такая инъективная параметризованная поверхность $r:D\to\R^3$, что $r(D)\subset\Pi$.

\begin{figure}[h]\centering
\includegraphics{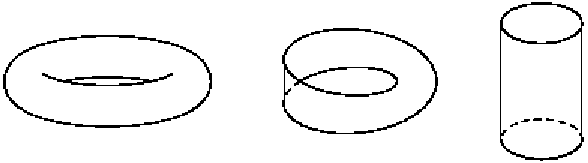}
\caption{Тор, лента Мебиуса и цилиндр}\label{tor}
\end{figure}

\smallskip
{\bf 1.} Нарисуйте следующие подмножества в $\R^3$.
Докажите, что они являются поверхностями.

(a) Квадрат со стороной $a$ на плоскости.  \quad

(b) Боковая поверхность прямого кругового цилиндра радиуса $R$ и высоты $h$ (рис. \ref{tor}).

(c) Боковая поверхность прямого кругового усеченного конуса с радиусами оснований $R,r$ и высотой $h$.

(d) Сфера радиуса $R$.

(e) Неформально, {\it тор} --- поверхность бублика, рис.  \ref{tor}  слева.
Формально, {\it стандартным тором с радиусами $R$ и $r$} называется фигура, образованная вращением окружности $(x-R)^2+y^2=r^2$ вокруг оси $Oy$.
Эта фигура получена из (двумерного) квадрата склейкой  пар его противоположных сторон
`с одинаковыми направлениями', т.е. без поворота.

(f) 
{\it Стандартной лентой Мебиуса} называется поверхность в $\R^3$, заметаемая стержнем длины 2,
равномерно вращающимся относительно своего центра, при равномерном движении этого центра
по окружности радиуса 9, при котором стержень делает пол-оборота, рис. \ref{tor}.
Эта фигура получена из длинной прямоугольной полоски склейкой двух ее противоположных сторон `с противоположным направлением', т.е. с поворотом на $180^\circ$.

(h) Поверхность вращения графика бесконечно дифференцируемой положительной функции $f:\R\to\R$.

(i) Седлообразная поверхность (гиперболический параболоид) $z=xy$.

(j) Двуполостный гиперболоид $z^2=x^2+y^2+1$.

(k) Однополостный гиперболоид $z^2=x^2+y^2-1$.

(l) График бесконечно дифференцируемой функции $f:D\to\R$.

(m) Прообраз нуля при бесконечно дифференцируемой функции  $f:\R^3\to\R$, производная (=градиент) которой ненулевая в каждой точке.


\smallskip
{\bf 2.} (a) Прямая движется поступательно с постоянной скоростью,
пересекая другую прямую под прямым углом, и одновременно
равномерно вращается вокруг этой прямой.
Поверхность, которую описывает движущаяся прямая, называется {\it прямым геликоидом}.
Составьте ее уравнение.
Везде ли это регулярная поверхность?

(b) Рассмотрим две параболы $A(t)=(1,t,t^2)$ и $B(t)=(-1,-t,t^2)$.
Напишите параметрическое уравнение поверхности, заметаемой прямой $A(t)B(t)$.

(с) Напишите неявное уравнение поверхности из (b).

(d) Составьте уравнение поверхности, образованной главными нормалями (т.е., векторами ускорений при равномерном движении) винтовой линии.

(e) Пусть $k>0$ и $\gamma:[0,1]\to\R^3$ --- параметризованная кривая, кривизна которой в любой точке не меньше $k$.
Через каждую точку кривой проведена нормальная плоскость, и в этой плоскости построена окружность с центром на кривой
и заданным радиусом $R$, причем $0<R<1/k$.
`Трубкообразная' поверхность, заметаемая этими окружностями в пространстве называется {\it трубкой} или
{\it каналовой поверхностью}.
Напишите ее уравнение.

(f) Докажите, что любая нормаль к трубке пересекает образ кривой $\gamma$ и перпендикулярна вектору скорости этой кривой.

\smallskip
{\bf 3.*} На круговом торе вращения, кроме параллелей и меридианов, являющихся плоскими окружностями, существует еще два семейства плоских окружностей, называемых {\it окружностями Виларсо}.
Они получаются пересечением тора его касательной плоскостью, касающейся тора в двух точках.
Напишите уравнения этих окружностей.
Проверьте, что все они имеют один и тот же радиус и пересекают все параллели тора под постоянным углом.
({\it Углом} между пересекающимися непараметризованными кривыми в их общей точке называется угол между их касательными в этой точке.)

\smallskip
Поверхности произвольной размерности в евклидовом пространстве $\R^m$ размерности $m$ определяются аналогично.
Они называются также {\it гладкими подмногообразиями} в $\R^m$.

\smallskip
{\bf 4.} 
(a) Любая поверхность в $\R^3$ становится поверхностью в $\R^4$, если рассмотреть $\R^3$ как подмножество в $\R^4$.

(b) Рассмотрим в $\R^4$ окружность $x^2+y^2=1$ и семейство ее нормальных трехмерных плоскостей.
{\it Бутылкой Клейна} (стандартной) $K$ называется поверхность в $\R^4$, заметаемая окружностью $\omega$, центр
которой равномерно описывает окружность $x^2+y^2=1$, а окружность $\omega$ в то же время равномерно поворачивается на угол $\pi$ (поворачивается в движущейся нормальной трехмерной плоскости относительно своего диаметра, движущегося вместе с нормальной трехмерной плоскостью).
Эта фигура получена из прямоугольной полоски такой склейкой ее пар противоположных сторон, при которой одна пара склеивается `с одинаковым направлением', а другая `с противоположным направлением'.

(c) Образ отображения $S^2\to\R^4$, заданного формулой $(x,y,z)\mapsto(x^2,xy,yz,zx)$, является поверхностью.

\smallskip
{\bf 5.} Не являются поверхностями ни объединение двух (или трех) координатных плоскостей в $\R^3$,
ни конус $z^2=x^2+y^2$  в $\R^3$.


\smallskip
{\bf 6.} Пусть $f:\R^n\to\R$ --- гладкая функция и $\grad f(x)\ne0$ для всех точек $x\in f^{-1}(0)$.
Тогда

(a) $f^{-1}(0)$ --- гладкое подмногообразие в $\R^n$. \quad

(b) Это подмногообразие ориентируемо (т.е. коориентируемо, см. \S\ref{0curmain}).

\smallskip
{\bf 7.} Докажите, что следующие множества матриц являются подмногообразиями в множестве $\R^{n^2}$
всех матриц размера $n\times n$:

(а) $GL(n,\R)=\{\mbox{вещественные $n\times n$-матрицы }A\ :\ \det A\ne0\}$.

(b) $SL(n,\R)=\{\mbox{вещественные $n\times n$-матрицы }A\ :\ \det A=1\}$.

(c) $SO(n,\R)=\{\mbox{вещественные $n\times n$-матрицы }A\ :\ AA^T=E,\ \det A=1\}$.

(d) $SU(2)=\{\mbox{комплексные $2\times2$-матрицы }A\ :\ A\overline A^T=E\}$.

(e) $SO(1,1)=\{\mbox{вещественные $2\times2$-матрицы }A\ :
\ AIA^T=\left(\begin{matrix} 1&\hphantom-0\\0&-1 \end{matrix}\right)\}$.

Для каждого из этих многообразий найти размерность, число
компонент связности, а также выяснить является~ли оно ограниченным.



\subsection{Длины кривых на поверхностях и изометрии}\label{0surlen}

{\bf 1.} (a) Найдите формулу для длины параметризованной сферической кривой
$\gamma:[a,b]\to S^2$ в сферических координатах $\varphi,\theta$.

(b) То же в декартовых координатах $x,y$.

(c) То же в {\it стереографических координатах}: паре чисел $(p,q)$ соответствует точка сферы, которая
переходит в точку $(p,q,0)$ при центральной проекции из точки $(0,0,1)$.

(d)* То же в {\it меркаторовских координатах}, которые определяются так.
На карте вводятся прямоугольные координаты $(u,v)$ такие, что любая прямая
на карте соответствует локсодромии (т.е. линии постоянного азимута ---
фиксированного положения стрелки компаса) на поверхности земного шара.

\smallskip
{\bf 2.} Выразите длину кривой
$[a,b]\overset{(u,v)}\to\R^2\overset{r}\to\R^3$ через
функции $u$ и $v$ (в формуле можно использовать не только алгебраические
выражения, но производные и интегралы) для

(a) $r(u,v)=(\cos u,\sin u,v)$ --- цилиндр с цилиндрической системой координат.
\quad

(b) $r(u,v)=(v\cos u,v\sin u,v)$ --- конус с конической системой координат.

(с) $r(u,v)=((2+\cos v)\cos u,(2+\cos v)\sin u,\sin v)$ --- тор с торической
системой координат.

(d) $r(u,v)=(f(v)\cos u,f(v)\sin u,h(v))$ --- поверхность вращения.


\smallskip
{\bf 3.} {\bf Теорема.} Длина кривой
$[a,b]\overset{(u,v)}\to\R^2\overset{r}\to\R^3$ равна
$$\int_a^b\sqrt{r_u^2(u'_t)^2+2r_u\cdot r_vu'_tv'_t+r_v^2(v'_t)^2}\ dt.$$
В этой формуле пропущены аргументы $t$ функций $u,v,u'_t,v'_t$;
\ $r_u=r_u(u(t),v(t))$ и аналогично для $r_v$.

\smallskip
{\it Углом} между пересекающимися параметризованными кривыми в их общей точке называется угол между их касательными в этой точке.

\smallskip
{\bf 4.} (a) Угол между параметризованными кривыми не зависит от их параметризации.

(b) Найдите угол между $r$-образами кривых $v=u+1$ и $v=3-u$ для
$r(u,v)=(v\cos u,v\sin u,v^2)$ ($r$-образом области $D$ является часть
параболоида $z=x^2+y^2$).

\smallskip
Напомним, что два подмножества пространства $\R^n$ называются {\it объемлемо изометричными},
если между ними существует {\it объемлемая изометрия} (движение), т.е. сохраняющее
расстояния (в $\R^n$!) отображение $\R^n\to\R^n$, переводящее первое во второе.

Две (непараметризованные) поверхности называются {\it внутренне изометричными}, если между ними
существует {\it внутренняя изометрия}, т.е. отображение одной в другую, сохраняющее длины всех кривых.

Большая часть дальнейшего материала мотивирована следующими двумя проблемами:
определить, являются ли данные (непараметризованные) поверхности

$\bullet$ внутренне изометричными?

$\bullet$ объемлемо изометричными?

\smallskip
{\bf 5.} (a) Прямоугольник на плоскости внутренне изометричен некоторой части
любого такого цилиндра, у которого диаметр больше стороны прямоугольника.

(a') Прямоугольник на плоскости внутренне изометричен некоторой части прямого кругового конуса.

(b) Единичный квадрат на плоскости внутренне изометричен некоторой части тора,
задаваемого в~$\R^4$ уравнениями $x_1^2+x_2^2=x_3^2+x_4^2=1$.

(c) Прямоугольник на плоскости не является объемлемо изометричным никакой части никакого цилиндра.

\smallskip
{\bf 6.}
(a) Сфера и плоскость не являются внутренне изометричными.

(b) Никакой круг на плоскости не является внутренне изометричным никакой части сферы.

(c) Никакой круг на сфере одного радиуса не является внутренне изометричным
никакой части сферы другого радиуса.

\subsection{Площадь поверхности}\label{0surare}


{\bf Площадью} называется отображение $S$ из семейства всех
двумерных поверхностей в луч $[0,+\infty)$, для которого выполнены следующие
условия.

(аддитивность) Если $\Pi,\Pi'$ и $\Pi\cup\Pi'$ поверхности,
причем $\Pi\cap\Pi'$ является объединением не более чем счетного семейства
кривых, то $S(\Pi\cup\Pi')=S(\Pi)+S(\Pi')$.

(монотонность) Если $f:\Pi\to f(\Pi)$ --- не увеличивающее длины кривых
отображение между поверхностями, то $S(f(\Pi))\le S(\Pi)$.

(нормировка) Площадь единичного квадрата на плоскости равна 1.

\smallskip
{\bf Теорема о площади.} {\it Такое отображение существует и единственно.}

\smallskip
В следующих задачах
\footnote{Решения задач из этой темы {\it особенно полезно} проверить с преподавателем. }
(кроме 4bd) можно пользоваться {\it существованием} из теоремы о площади
(а также аналогичной теоремой и всеми другими результатами для площадей плоских фигур).
В предположении существования можно найти площадь произвольной поверхности и этим доказать {\it единственность}.
Далее, {\it не используя предположение о существовании,} можно проверить
выполнение свойств площади для найденного отображения $S$ и этим доказать существование.

\smallskip
{\bf 1.} Внутренняя изометрия сохраняет площади.


\smallskip
{\bf 2.} (a) Площадь сферического двуугольника с углом $\alpha$ и диаметрально
противоположными вершинами равна $2\alpha$.

(b) {\bf Теорема.} Площадь сферического треугольника с углами $\alpha$, $\beta$ и $\gamma$
равна $\alpha+\beta+\gamma-\pi$.

(c) {\bf Теорема.}
Площадь сферического многоугольника с углами $\alpha_1,\dots\alpha_n$ равна
\linebreak
$\alpha_1+\dots+\alpha_n-(n-2)\pi$.


\smallskip
{\bf 3.} (a) Площадь круга радиуса $R$ на сфере равна $2\pi(1-\cos R)$.

(b) Площадь поверхности, образованной вращением графика функции $f:[a,b]\to [0,+\infty)$ вокруг оси $Ox$, равна
$2\pi\int_a^b f(x)\sqrt{1+f'(x)^2}\ dx$.

(c) {\it Первая теорема Гюльдена.} Если поверхность образована вращением  вокруг оси некоторой линии, лежащей в
одной плоскости с осью и целиком по одну сторону от оси, то площадь поверхности равна произведению длины линии на длину окружности, описанной центром тяжести линии.

(d) {\it Вторая теорема Гюльдена.} Если тело образовано вращением вокруг оси некоторой плоской фигуры, лежащей
в одной плоскости с осью и целиком по одну сторону от оси, то объем этого тела равен произведению площади фигуры на длину окружности, описанной центром тяжести фигуры.

\smallskip
{\bf 4.} (a)* {\bf Теорема.}
Площадь элементарной непараметризованной поверхности $r(D)$ равна
$S(r(D))=\int\int_D |r_u\times r_v| dudv$ [Ra03].

(b) Докажите без использования теоремы о площади, что предыдущее выражение не
зависит от выбора параметризации $r$.

(c) Докажите {\it единственность} в теореме о площади (предполагая существование).

(d)* Докажите {\it существование} в теореме о площади (не предполагая единственности!).
Для проверки монотонности полезно понятие римановой метрики (\S\ref{0polyriem}) [Gr94].

(e) Найдите площадь пересечения цилиндра $x^2+y^2\le1$ и гиперболоида $z=xy$.

\smallskip
Аналогично площади двумерных поверхностей определяется {\bf $n$-мерный объем} $n$-мерных поверхностей.
Аналогично двумерному случаю доказывается теорема о существовании и единственности $n$-мерного объема, а также
следующий результат:

{\it $n$-мерный объем элементарной непараметризованной поверхности $r:D\to\R^m$ равен}
$V(r(D))=\int\dots\int_D |r_{u_1}\wedge\dots\wedge r_{u_n}|du_1\dots du_n.$
Здесь $|r_{u_1}\wedge\dots\wedge r_{u_n}|$ --- $n$-мерный объем $n$-мерного
параллелепипеда, натянутого на векторы $r_{u_1},\dots,r_{u_n}$.

\subsection{Геодезические}\label{0surgeo}

{\it Расстоянием по поверхности} $\Pi$ между точкам $P,X\in\Pi$ называется
инфимум $|P,X|$ длин кривых на этой поверхности, соединяющих $P$ и $X$.

Ясно, что {\it внутренняя изометрия сохраняет расстояния по поверхности.}
Обратное (т.е. то, что отображение, сохраняющее расстояния по поверхности, сохраняет длины всех кривых на
поверхности) доказано ниже с использованием римановой метрики (теорема \ref{0polyriem}.4).

Непараметризованная кривая $\Gamma$ на поверхности $\Pi$ называется {\bf геодезической на поверхности $\Pi$}, если она {\it локально кратчайшая}, т. е. если любая точка $x\in\Gamma$ имеет такую окрестность $U$ в поверхности, что расстояние по поверхности между любыми точками $y_1,y_2\in U\cap\Gamma$ равно длине отрезка кривой $\Gamma$ от $y_1$ до $y_2$.

\smallskip
{\bf 0.} В этой задаче рассматриваются геодезические на поверхностях многогранников.
Их определение аналогично.

(a) Нарисуйте на кубе геодезическую, соединяющую его противоположные вершины.

(b) Нарисуйте геодезическую на прямоугольном параллелепипеде $a\times b\times c$,
соединяющую середины параллельных ребер, не лежащих в одной грани.

(c) Сумма плоских углов выпуклого многогранного угла не превосходит $2\pi$.

(d) Геодезическая на поверхности выпуклого многогранника не проходит через
его вершины (т.е. может в них только начинаться или заканчиваться) и
проходит через ребра по закону `угол падения равен углу отражения'.

\smallskip
{\bf 1.} (a) Кратчайшая кривая, соединяющая заданные две точки, является геодезической.

(b) Геодезическая не обязательно является кратчайшей.

(c)* Любые две точки компактной поверхности можно соединить кратчайшей кривой

(d)* Для любой компактной поверхности существует такое $\varepsilon>0$, что любая
геодезическая длины менее $\varepsilon$ является кратчайшей.

(e) Найдите уравнение (в декартовых координатах в $\R^3$) хотя бы одной геодезической на поверхности $z=xy$.

\smallskip
{\it Производной} отображения $f:\R^k\to\R^n$ в точке $x_0\in\R^k$ называется линейное отображение $A:\R^k\to\R^n$, что
$f(x)=f(x_0)+Ax+o(x)$ при $x\to x_0$.
Будем называть производной также аффинное отображение $x\mapsto f(x_0)+Ax$.

{\it Касательной плоскостью} $T_P=T_P\Pi$ к поверхности $\Pi$ в точке $P\in\Pi$ называется образ плоскости $\R^2$ при производной в точке $r^{-1}(P)$ для некоторого отображения $r:D\to\R^3$ из определения поверхности.

\smallskip
{\bf 2.} (a) Этот образ есть плоскость, проходящая через $P$ и содержащая векторы $r_u(r^{-1}(P))$ и $r_v(r^{-1}(P))$.

(b) Это определение корректно, т.е. не зависит от выбора отображения $r$.

\smallskip
{\bf 3.} (a) Внутренняя изометрия переводит геодезические в геодезические.

(b)* {\bf Теорема.} Образ параметризованной кривой $\gamma:[a,b]\to\Pi$ с постоянной по модулю скоростью является геодезической на поверхности $\Pi$ тогда и только тогда, когда $\gamma''(t)\perp T_{\gamma(t)}\Pi$ для любого $t$
(т.е. когда вектор ускорения $\gamma''(t)$ перпендикулярен плоскости, касательной к поверхности в точке $\gamma(t)$ для любого $t$).

{\it Эвристическое пояснение к прямому доказательству части `только тогда'.}
Пусть, напротив, $\gamma''(t)$ не перпендикулярно поверхности в некоторой точке
$t$.
Тогда проекция кривой $\gamma$ на касательную
плоскость в точке $\gamma(t)$ имеет ненулевую кривизну.
Значит, эту проекцию можно `спрямить'.
Тогда и исходную кривую можно `спрямить'.

{\it Доказательство} получается из уравнения Эйлера-Лагранжа для функционала длины или его квадрата [Ra03].

(c)* При движении по геодезической правое и левое колеса узкого автомобиля
проделают одинаковый путь с точностью до малых порядка квадрата ширины
автомобиля. (Это следует из минимальности геодезических.)

(d)* Нарисуем на яйце (гладкой поверхности, лежащей по одну сторону от любой
касательной плоскости и пересекающей эту плоскость ровно в одной точке;
например, поверхности вращения графика выпуклой функции) произвольную кривую.
Прокатим яйцо по плоскости (без вращения) вдоль этой кривой.
Кривая на яйце
является геодезической тогда и только тогда, когда соответствующая кривая на
плоскости является прямой. (Это следует из предыдущего пункта.)


\smallskip
Параметризованная кривая $\gamma:[a,b]\to\Pi$ называется {\it параметризованной геодезической} на поверхности $\Pi$, если $\gamma''(t)\perp T_{\gamma(t)}\Pi$ для любого $t$.
(Иными словами, если для муравья, живущего в поверхности и не видящего объемлемого пространства и нормальных составляющих векторов, вектор ускорения нулевой, т.е., скорость `постоянна'.)

Теорема 3b означает, что {\it непараметризованная кривая на поверхности является геодезической тогда и только тогда, когда она имеет параметризацию, являющуюся параметризованной геодезической.}
В дальнейших задачах этой теоремой можно пользоваться без доказательства.

Ясно, что параметризованная геодезическая имеет постоянную по модулю скорость.
(От величины постоянного модуля скорости свойство параметризованной кривой быть
параметризованной геодезической не зависит.)
В дальнейшем на
геодезических рассматриваются параметризации с постоянной по модулю скоростью, и под
геодезической понимается либо геодезическая, либо параметризованная геодезическая
(возникающие шероховатости формулировок решатель легко исправит).

\smallskip
{\bf 5.} Найдите все геодезические на поверхностях из задач \ref{0surexa}.1.abcd.
(Для (a,b,d) нужен `явный' ответ, для (c) --- уравнение.)

\smallskip
{\bf 4.} (a) Прямая на поверхности --- геодезическая.

(b) Меридиан поверхности вращения --- геодезическая.

(c) Параллель поверхности вращения является геодезической тогда и только
тогда, когда касательная к меридиану в ее точках параллельна оси вращения.

(d)* Кривая, являющаяся связной компонентой множества неподвижных точек
некоторой изометрии поверхности, является геодезической.


\smallskip
{\bf 6.} 
{\it Теорема Клеро}.
Для параметризованной геодезической $\gamma:\R\to\Pi$ на~поверхности вращения $\Pi$ величина
$r(t)\sin\angle(\gamma'(t),m(t))$ не зависит от $t$.
Здесь $r$ --- расстояние до~оси вращения и $m$ --- меридиан.

\subsection{Уравнение геодезических и экспоненциальное отображение}\label{0surexp}

{\bf 0.} Напишите уравнение параметризованных геодезических на поверхности $z=xy$
(в выбранной Вами системе координат на этой поверхности).

\smallskip
{\bf 1.} Пусть $x_1,x_2:\R\to\R$ и $\gamma=r(x_1,x_2)$ --- параметризованная
кривая на параметризованной поверхности $r:D\to\R^3$.
Будем далее пропускать аргумент $(x_1(t),x_2(t))$ у функции $r$ и ее частных производных, а также обозначать
штрихом дифференцирование по $t$, а индексами у $r$ (но не у $\Gamma$ и $g$) --- частное дифференцирование
по соответствующей переменной (решателю будет полезно первое время кроме индексов писать еще штрихи).

(a) Кривая $\gamma$ является параметризованной геодезической тогда и только тогда, когда
$\gamma''\cdot r_1=\gamma''\cdot r_2=0$ для любого $t$.

(b) $\gamma'=r_1x_1'+r_2x_2'$ --- вектор скорости кривой $\gamma$ в точке $r(x_1(t),x_2(t))$.

(c) $\gamma''=x_1''r_1+x_1'r_1'+x_2''r_2+x_2'r_2'$.

(d) $r_1'=r_{11}x_1'+r_{12}x_2'$.

(e) {\bf Теорема.} Уравнение параметризованной геодезической (на функции $x_1$ и $x_2$ через данную функцию $r$)
$$\begin{cases}
-x_1''=\Gamma^1_{11}(x_1')^2+(\Gamma^1_{21}+\Gamma^1_{12})x_1'x_2'+
\Gamma^1_{22}(x_2')^2\\
-x_2''=\Gamma^2_{11}(x_1')^2+(\Gamma^2_{21}+\Gamma^2_{12})x_1'x_2'+
\Gamma^2_{22}(x_2')^2
\end{cases}\quad\mbox{или}
\quad x_k''+\sum_{i,j}\Gamma^k_{ij}x_i'x_j'=0,\quad\mbox{где}$$
$$g_{ij}=r_i\cdot r_j\quad\mbox{и}\quad
\Gamma^k_{ij}:=\frac{g_{3-k,3-k}r_k\cdot r_{ij}-g_{k,3-k}r_{3-k}\cdot r_{ij}}
{\det g}\quad \text{--- символы Кристоффеля}.$$
\quad

{\bf 2.} (a) {\bf Следствие.}
Через каждую точку в каждом направлении на поверхности проходит ровно одна геодезическая.

(b) Вычислите символы Кристоффеля для сферической системы координат на сфере.


\smallskip
{\it Круг} и {\it окружность} на поверхности определяются аналогично случаю плоскости (через расстояние {\it по поверхности}).

\smallskip
{\bf 3.} (a) Окружность достаточно малого радиуса на поверхности является кривой (значит, ее длина определена).

(b) Круг достаточно малого радиуса на поверхности является поверхностью (значит, его площадь определена).

(с) Обозначим через $S_{\Pi,P}(R)$ площадь круга на $\Pi$ радиуса $R$ с центром в точке $P\in\Pi$, а через $L_{\Pi,P}(R)$ длину окружности на $\Pi$ радиуса $R$ с центром в точке $P\in\Pi$.
Тогда $L_{\Pi,P}(R)=S'_{\Pi,P}(R)$.

\begin{figure}[h]\centering
\includegraphics{p.80}
\caption{Экспоненциальное отображение и квадратичная форма Риччи}\label{f:expric}
\end{figure}

\smallskip
Далее в этом пункте $\Pi\subset\R^m$ --- поверхность размерности $n$.

Для $P\in\Pi$ определим {\it (геодезическое) экспоненциальное отображение}
$$\exp=\exp\phantom{}_P:X\to\Pi\quad\mbox{формулой}
\quad\exp(a):=\gamma_{P,a}(1),$$
где $\gamma_{P,a}:[-1,1]\to\Pi$ --- та геодезическая, для которой
$\gamma_{P,a}(0)=P$ и $\gamma'_{P,a}(0)=a$.
Здесь $X\subset T_P$ --- множество тех векторов $a$, для которых определено $\gamma_{P,a}(1)$.

Для $n=1$ экспоненциальное отображение является натуральной параметризацией.

\smallskip
{\bf 4.} (a) Для единичной сферы~$S^2$ при экспоненциальном отображении
$T_{(1,0,0)}\to S^2$ полярные координаты переходят в
сферические:
$$\exp(1,\rho\cos\varphi,\rho\sin\varphi)=(\cos\rho',\sin\rho'\cos\varphi,\sin\rho'\sin\varphi),
\quad\text{где}\quad \rho':=\frac\pi2-\rho.$$
\quad
(b) Обобщите этот результат для поверхностей вращения.

\smallskip
{\bf 5.}
(a) Прямая, проходящая через $P$, переходит при экспоненциальном отображении $\exp_P$ в геодезическую.

(b) Верно ли, что для любой поверхности и точки на ней образ любой прямой при
экспоненциальном отображении является геодезической?

(c) Образом шара (в $T_P$) радиуса $R$ с центром в $P$ при экспоненциальном
отображении $\exp_P$ является шар (на поверхности) радиуса $R$.

(d) Верно ли, что для любой поверхности и точки на ней образ любого шара при
экспоненциальном отображении является шаром?


\subsection{Параллельный перенос}\label{0surtra}

Все знают, что такое параллельный перенос на плоскости.
Можно определить параллельный перенос и на искривленной поверхности.
Более того, это определение необходимо для решения интересных задач из географии и физики.

Начнем с поверхности многогранника.
(Заметим, что она не является поверхностью в смысле определения,
принятого в данной книге.)
Рассмотрим две соседние грани данного многогранника.
Если повернуть плоскость одной из них вокруг их общего ребра, то она
совместится с плоскостью другой грани [Ta89, рис. 15].
Этот поворот можно рассматривать как перекатывание многогранника с одной
грани на другую через его ребро.
Если на первой грани нарисовать
вектор, то после поворота он отпечатается на плоскости второй грани.
Такой перенос вектора с одной грани на другую называется {\it параллельным
переносом} вектора через ребро [Ta89, рис. 16, 17].

\smallskip
{\bf 1.} (a)
Замкнутая ломаная на выпуклом многограннике не содержит вершин и ограничивает область, сумма плоских углов в вершинах которой равна $S$.
После параллельного переноса вдоль этой ломаной вектор повернется на угол $-S$.

Здесь поверхность многогранника ориентирована, имеется в виду ориентированный угол на ориентированной поверхности и ломаная направлена так, что ее направление вместе с нормальным к ребру ломаной вектором, касающимся грани многогранника, дает данную ориентацию на поверхности многогранника.
Значит, от выбора ориентации поверхности многогранника угол не зависит.

(b)* Этот угол поворота равен (с точностью до $2\pi n$) сумме величин телесных
углов, {\it двойственных} к многогранным углам $A_1,\dots,A_n$.

\smallskip
Приведем два неформальных описания параллельного переноса на поверхностях (уже в смысле этой книги).

Муравей живет в поверхности и не видит объемлемого пространства, в частности, нормальных составляющих векторов.
Касательное к поверхности $\Pi$ векторное поле $v(t)=\vec{A(t)B(t)}$, $A(t)\in\Pi$, {\it параллельно},
если вектор $v'_t(t)$ нулевой, т.е., векторное поле `постоянно', с точки зрения муравья.

Нарисуем на поверхности кривую и проведем касательную к поверхности плоскость в некоторой точке кривой.
Параллельный перенос касательной плоскости вдоль кривой на поверхности --- это качение касательной плоскости
вдоль кривой без проскальзывания (при котором плоскость остается касательной).
Это движение касательной плоскости по неподвижной поверхности задается следующим условием:
{\it мгновенная ось вращения
касательной плоскости касается поверхности и перпендикулярна данной кривой.}

Формально, пусть дана поверхность $\Pi\subset\R^3$ и параметризованная кривая $\gamma:[a,b]\to\Pi$.
Касательное к поверхности $\Pi$ векторное поле $v(t)$ на кривой $\gamma[a,b]$
называется {\bf параллельным вдоль данной кривой}
(в смысле Леви-Чивита), если вектор $v'_t(t)$ перпендикулярен (плоскости,
касательной к) поверхности в точке $\gamma(t)$ при любом $t$.

Вектор $v(\gamma(b))$ называется {\it вектором, полученным из вектора
$v(\gamma(a))$ параллельным переносом вдоль данной кривой}.

\smallskip
{\bf 2.} (a) Какие касательные векторы к плоскости, лежащей в трехмерном
пространстве, получаются друг из друга параллельным переносом?

(b) Поле векторов скорости параметризованной кривой на поверхности
является параллельным вдоль этой кривой тогда и только тогда, когда
эта кривая является параметризованной геодезической.

(c) Параллельность вдоль параметризованной кривой {\it с данным образом}
не зависит от выбора этой кривой.

(d) Результат параллельного переноса вдоль кривой {\it с данными концами}
может зависеть от выбора этой кривой.

\smallskip
Векторное поле на поверхности называется {\it параллельным вдоль
непараметризованной кривой}, если оно параллельно вдоль любой ее параметризации.
Это определение корректно ввиду задачи 2c.

\smallskip
{\bf 3.} (a) Непрерывное семейство векторов одинаковой длины на меридиане
поверхности вращения, касающихся параллелей, параллельно вдоль меридиана.

(b) Дана поверхность вращения гладкой положительной функции $f$.
На какой угол в $\R^3$ повернется вектор, касательный к меридиану, при
параллельном переносе из точки $(a,f(a),0)$ в точку $(b,f(b),0)$ вдоль
меридиана?

(c) При объемлемой изометрии поверхностей семейство векторов, параллельное вдоль
некоторой кривой, переходит в семейство векторов, параллельное вдоль образа этой кривой.

\smallskip
{\bf 4.} Если на данной поверхности семейства векторов $u$ и $v$ параллельны
вдоль данной кривой, то
\quad

(a) $|v(x)|=|v(y)|$. \quad
(b) $u(x)\cdot v(x)=u(y)\cdot v(y)$. \quad
(c) $\angle(u(x),v(x))=\angle(u(y),v(y))$. \quad

(d) Семейства $u+v$ и $3u$ параллельны вдоль той же кривой.

\smallskip
{\bf 5.} (a) {\bf Теорема.} Параллельный перенос на данной поверхности
вдоль данной кривой определяет ортогональное отображение касательных пространств.

(b) Семейство векторов является параллельным вдоль геодезической
тогда и только тогда, когда модуль вектора семейства и угол между вектором
семейства и вектором скорости геодезической постоянны вдоль геодезической.

\smallskip
{\bf 6.} Пусть $v(t)=a_1(t)r_1\bigl(x_1(t),x_2(t)\bigr)+
a_2(t)r_2\bigl(x_1(t),x_2(t)\bigr)$ ---
касательный к поверхности вектор в точке $r(x_1(t),x_2(t))$.
Будем далее обозначать штрихом производную по $t$ и пропускать аргументы
функций.

(a) $v'=a_1'r_1+a_1r_1'+a_2'r_2+a_2r_2'$.

(b) Семейство $v$ параллельно вдоль кривой $\gamma=r(x_1,x_2)$ тогда и только
тогда, когда $v'\cdot r_1=v'\cdot r_2=0$.

(c) {\bf Теорема.} Уравнение параллельного переноса (на функции $a_1$ и $a_2$
через данные функции $r,x_1,x_2$)
$$\begin{cases}
-a_1'=(\Gamma^1_{11}x_1'+\Gamma^1_{21}x_2')a_1+
(\Gamma^1_{12}x_1'+\Gamma^1_{22}x_2')a_2\\
-a_2'=(\Gamma^2_{11}x_1'+\Gamma^2_{21}x_2')a_1+
(\Gamma^2_{12}x_1'+\Gamma^2_{22}x_2')a_2
\end{cases}
\quad\mbox{или}\quad a_k'+\sum_{i,j}\Gamma^k_{ij}x_i'a_j=0.$$
\quad
(d) Любой вектор можно параллельно перенести вдоль любой кривой.

\smallskip
{\it При внутренней изометрии поверхностей семейство векторов, параллельное вдоль
некоторой кривой, переходит в семейство векторов, параллельное вдоль образа этой кривой.}
Это доказано в теореме \ref{0polyriem}.6.c с использованием римановой метрики. 

\begin{figure}[h]\centering
\includegraphics{p.60}
\caption{Параллельный перенос вектора по замкнутому контуру}\label{f:seccur}
\end{figure}

 Объяснение феномена маятника Фуко с использованием параллельного переноса приводится в [Ta89].

\smallskip
{\bf 7.} На какой угол повернется касательный вектор при параллельном переносе
вдоль
\quad

(a) параллели на цилиндре? \quad

(b) контура треугольника с углами $\alpha,\beta,\gamma$ на сфере?

(c) параллели $z=1$ на конусе $z^2=x^2+y^2$? \quad

{\it Указание.} Используйте (доказанную ниже) инвариантность параллельного переноса при внутренней изометрии.

(d) параллели $\theta=\theta_0$ на сфере?

(e) данной параллели данной поверхности вращения?

Здесь поверхность ориентирована, имеется в виду ориентированный угол на ориентированной поверхности, и кривая   направлена так, что ее направление вместе с нормальным к кривой вектором, касающимся поверхности, дает данную ориентацию на поверхности.
Значит, от выбора ориентации поверхности угол не зависит.


\smallskip
{\bf 8.} Пусть $N,N'$ --- поверхности, касающиеся вдоль кривой~$\gamma$.
Докажите, что результат параллельного переноса вдоль кривой~$\gamma$
одинаков для $N$ и $N'$.


\subsection*{Указания и решения к некоторым задачам}
\addcontentsline{toc}{subsection}{Указания и решения к некоторым задачам}

{\bf \ref{0sursph}.1.} (a) Ответ: объединение северного полюса и бесконечного семейства параллелей.

(b) Ответ: на северный полюс. Используйте (с).

(c) Из малого прямоугольного сферического треугольника с катетами `$\cos\theta d\varphi$' и `$d\theta$'
получаем $\theta'=\cos\theta\ctg\alpha$ для $\theta:=\theta(\varphi)$.
(Треугольник можно считать плоским ввиду малости, или использовать задачу 8b и перейти к пределу.)

\smallskip
{\bf \ref{0sursph}.4.} (c) Используйте двойственный трехгранный угол.

\smallskip
{\bf \ref{0sursph}.6.} (a) Опустите перпендикуляр из ребра на грань.

(b) Следует из (a).

(c) Опустите перпендикуляр из вершины на треугольное сечение.

(A,B,C) Используйте (a,b,c) и 4.d.


\smallskip
{\bf \ref{0sursph}.8.}
(c) Возьмем соответствующий трехгранный угол $OA'B'C'$, для которого $OA'=OB'=OC'=1$.
Опустим из $A'$ перпендикуляр $A'A_2$ на плоскость $OB'C'$, а также перпендикуляры $A'B_1$ и $A'C_1$
на прямую $OB'$ и $OC'$, соответственно.
($A_2$ не обязательно лежит внутри угла $B'OC'$, $C_1$ и $B_1$ не обязательно лежат на лучах $OС'$ и $OB'$, соответственно.)
По теореме о трех перпендикулярах $A_2B_1\perp OB'$ и $A_2C_1\perp OC'$.
Значит,  $\angle A'B_1A_2=C$ и $\angle A'C_1A_2=B$.
Тогда
$$A'A_2=A'B_1\sin C = \sin\beta  \sin C\quad\text{и}\quad A'A_2=A'C_1\sin B = \sin\gamma\sin B.$$
\quad
(с) {\it Другое решение.} Выразите $\sin^2A$ из (d).

(d) Возьмем соответствующий трехгранный угол $OA'B'C'$, для которого $OA'=OB'=OC'=1$.
Опустим из $C'$ и $B'$ перпендикуляры $C'C_1$ и $B'B_1$ на прямую $OA'$
($C_1$ и $B_1$ не обязательно лежат на луче $OA'$).
Обозначим $a:=\vec{OA'}$, $b:=\vec{OB'}$, $c:=\vec{OC'}$.
Тогда $a\cdot c=\cos\beta$, $b\cdot c=\cos\alpha$ и $a\cdot b=\cos\gamma$.
Имеем $\vec{C'C_1}=-c+\cos\beta a$ и $\vec{B'B_1}=-b+\cos\gamma a$.
Значит,
$$C'C_1=\sin\beta,\quad B'B_1=\sin\gamma\quad\text{и}\quad\vec{C'C_1}\cdot\vec{B'B_1}=\cos\alpha-\cos\beta\cos\gamma.$$
Тогда $\cos A=\dfrac{\vec{C'C_1}\cdot\vec{B'B_1}}{C'C_1\cdot B'B_1}$.

(e) Используйте (d) и 4.d.

\comment

\bigskip
{\bf \ref{0surpro}.3.} (b) За счет применения центрального проектирования можно считать, что
$BA'||CA'$, $BC'||AC'$ и $AB'||CB'$ (т.е. что прямая $A'B'$ бесконечно удаленная).
Так как точки $A', B', C'$ различны, то три рассматриваемые точки конечны.
Обозначим через $l$ прямую, проведенную через две из них.
Если $l$ пересекает $AC$, то примените теорему о композиции гомотетий.
Если $l||AC$, то примените теорему о композиции переносов.

\smallskip
{\bf \ref{0surpro}.5.} (a) 
Обозначим через $O$ центр проектирования.
Используйте, что
\linebreak
$\dfrac{AC}{BC}=\dfrac{S_{AOC}}{S_{BOC}}=\dfrac{AO\sin\angle AOC}{BO\sin\angle BOC}$.

\smallskip
{\bf \ref{0surpro}.6.} (b,c) Используйте (a) и 5a.

\smallskip
{\bf \ref{0surpro}.8.} См. указание к \ref{0surlen}.1.c.

\smallskip
{\bf \ref{0surpro}.9.} Напишите формулу для $f_{P,A,B}(x,y)$, используя указание к задаче 8.

(c) Докажите, что существуют
{\it другие однородные координаты},
в которых данное множество задается одним из уравнений
$x^2=0$, \ $x^2+y^2=0$, \ $x^2-y^2=0$, \ $x^2+y^2+z^2=0$, \ $x^2+y^2-z^2=0$.

\endcomment

\bigskip
{\bf \ref{0surexa}.1.}
(f) 
$r(u,v)\ =\ \left(\ (9+ u\sin v)\cos 2v,\ (9+u\sin v)\sin2v,\ u\cos v\ \right)$, $u\in[-1,1]$, $v\in[0,\pi]$.
Обозначим через $A(t)$ поворот пространства на $t$ относительно оси $z$.
Тогда $r(u,v)=A(2v)[(9,0,0)+u(\sin v,0,\cos v)]$.
Поэтому $r_u(u,v)=A(2v)(\sin v,0,\cos v)$ сонаправлено со стрежнем.
Раскрывая $r_v(u,v)$ по правилу Лейбница, получаем две ненулевые перпендикулярные составляющие, перпендикулярные стержню.
Поэтому $r_u$ и $r_v$ линейно независимы.

\smallskip
{\bf \ref{0surexa}.2.} (d) $r(u,v)=(\cos t,\sin t,t)-v(\cos t,\sin t,0)$.

(e) $r(u,v)=\gamma(u)+Rn(u)\cos v+Rb(u)\sin v$, где $n(u)$ и $b(u)$ --- вектора нормали и бинормали, см. \S\ref{0torsion}.

\bigskip
{\bf \ref{0surlen}.1.}
Сначала выразите координаты на сфере через декартовы координаты в $\R^3$.

(c) Точки $(x_k,y_k,z_k)\in\R^3$, $k=1,2,3$, лежат на одной прямой тогда и только тогда, когда
$\dfrac{x_1-x_2}{x_1-x_3}=\dfrac{y_1-y_2}{y_1-y_3}=\dfrac{z_1-z_2}{z_1-z_3}$.
Здесь в знаменателе допускается 0:
первое равенство нужно понимать как наглядную запись равенства $(x_1-x_2)(y_1-y_3)=(y_1-y_2)(x_1-x_3)$.

(d) Используйте результат задачи \ref{0surexa}.1.c.

\smallskip
{\bf \ref{0surlen}.5.} (a,a',b) Задайте изометрию формулой. Используйте теорему 3.

(c) При движении пространства плоскости переходят в плоскости.

\smallskip
{\bf \ref{0surlen}.6.} Используйте задачу \ref{0surgeo}.5.d.


(c) Аналогично (a,b).

\bigskip
{\bf \ref{0surgeo}.0.} (c) Проведите через вершину угла луч, лежащий внутри угла,
и через точку на этом луче (отличную от вершины угла) плоскость, перпендикулярную лучу.

{\bf \ref{0surgeo}.1.} (e) $x=z=0$. 

\smallskip
{\bf \ref{0surgeo}.5.} Для сферы используйте задачи 4d или \ref{0surexp}.2.a.
Для цилиндра и конуса используйте инвариантность при изометрии.

\smallskip
{\bf \ref{0surgeo}.6.}  Утверждение равносильно тому, что $e\wedge\gamma\wedge\gamma'=const$, где $\gamma$ --- геодезическая и $e$ --- единичный вектор, параллельный оси вращения.
Последнее условие доказывается дифференцированием по правилу Лейбница, т.к. $e\wedge\gamma\wedge\gamma''=0$.

\bigskip
{\bf \ref{0surexp}.0.} См. 1.

\smallskip
{\bf \ref{0surexp}.3.} (b) Следует из 5c.

\bigskip
{\bf \ref{0surtra}.7.} (d) Используйте (c) и задачу 8.

\newpage
\section{Числовые кривизны поверхностей}\label{0cur}

\subsection{Скалярная кривизна}\label{0cursca}

{\bf Скалярной кривизной} поверхности $\Pi$ во внутренней точке $P$ называется число
$$\tau=\tau_{\Pi,P}:=6\lim\limits_{R\to0}\frac{2\pi R-L_{\Pi,P}(R)}{\pi R^3}.$$
Далее для скалярной и других кривизн $P$ и $\Pi$ пропускаются из обозначений, поскольку ясны из контекста.

\begin{figure}[h]\centering
\includegraphics{p.10}
\caption{Скалярная кривизна поверхности}\label{f:scalcur}
\end{figure}


\smallskip
{\bf 1.} (abcd) Найдите скалярную кривизну в точках поверхностей из задач \ref{0surexa}.1.abcd.

(e)*  Скалярная кривизна в точке $(0,0,0)$ отрицательна для седлообразной поверхности $z=xy$.

\smallskip
{\bf 2.} (a) Как изменяется скалярная кривизна при гомотетии пространства?

(b) Внутренняя изометрия сохраняет скалярную кривизну.

\smallskip
Для решения задач 3, 4, 6 нужны задача \ref{0surexp}.3.c и \ref{0polyriem}.8.a, см. задачу \ref{0polyriem}.9.a.

\smallskip
{\bf 3.} (a) $\lim\limits_{R\to0}\dfrac{2\pi R-L(R)}{R^2}=0.$

(c)* $L_{\Pi,P}(R)=2\pi R-\dfrac{\pi\tau R^3}6+O(R^5)$ [Gr90, BBB06].

(d)* Существует ли отображение поверхностей, сохраняющее скалярную кривизну, но не являющееся внутренней изометрией?

\smallskip
{\bf Теорема.} Элементарная непараметризованная двумерная поверхность внутренне изометрична некоторой части
плоскости тогда и только тогда, когда ее скалярная (или секционная, или гауссова, см. далее) кривизна
равна нулю в каждой точке [MF04, Ra04].

\smallskip
Для доказательства этой просто формулируемой теоремы, как и для получения
формул для вычисления скалярной кривизны (задача \ref{0polyriem}.9), нужно изучить \S\ref{0poly}.

\smallskip
{\bf 4.} Обозначим через $S_{\Pi,P}(R)$ площадь круга на $\Pi$ радиуса $R$ с центром в точке $P\in\Pi$.

(a) $\tau=24\lim\limits_{R\to0}\dfrac{\pi R^2-S_{\Pi,P}(R)}{\pi R^4}$.

(b)* $S_{\Pi,P}(R)=\pi R^2-\dfrac{\pi\tau}{24}R^4+O(R^6)$ [Gr90].


\smallskip
Для трехмерной поверхности $\Pi\subset\R^m$ и точки $P\in\Pi$ обозначим через
$S_{\Pi,P}(R)$ площадь сферы на $\Pi$ радиуса $R$ с центром в точке $P\in\Pi$.
{\bf Скалярной кривизной} поверхности $\Pi$ в точке $P$ называется число
$$\tau=\tau_{\Pi,P}:=6\lim\limits_{R\to0}\frac{4\pi R^2-S_{\Pi,P}(R)}{(4/3)\pi R^4}.$$

{\bf 5.} Вычислите скалярную кривизну точек следующих трехмерных поверхностей в $\R^4$:

(a) гиперплоскости $\R^3$;  \quad (b) цилиндра $S^1\times\R^2$; \quad (c) сферы $S^3$;  \quad

(d)* цилиндра $S^2\times\R$; \quad (e) конуса $t^2=x^2+y^2+z^2$;  \quad (f)* тора $S^2\times S^1$.

\smallskip
Пусть $n\ge2$ и $\Pi$ --- $n$-мерная поверхность в пространстве $\R^m$.
Обозначим через

$\bullet$ $V_n$ $n$-мерный объем (\S\ref{0surare}) шара радиуса 1 в $\R^n$,

$\bullet$ $S_n$ $(n-1)$-мерный объем поверхности шара радиуса 1 в $\R^n$,

$\bullet$ $V_{\Pi,P}(R)$ $n$-мерный объем шара на $\Pi$ радиуса $R$ с центром
в точке $P\in\Pi$,

$\bullet$ $S_{\Pi,P}(R)$ $(n-1)$-мерный объем поверхности шара на $\Pi$
радиуса $R$ с центром в точке $P\in\Pi$.

{\bf Скалярной кривизной} поверхности $\Pi$ в точке $P$ называется число
$$\tau=\tau_{\Pi,P}:=6\lim\limits_{R\to0}
\frac{S_nR^{n-1}-S_{\Pi,P}(R)}{V_nR^{n+1}}.$$

{\bf 6.} (b)  $S_{\Pi,P}(R)=V'_{\Pi,P}(R)$.

(c) $S_n=nV_n$.

(d) $\tau=6(n+2)\lim\limits_{R\to0}\dfrac{V_nR^n-V_{\Pi,P}(R)}{V_nR^{n+2}}$.

(e)* $V_P(R)=V_nR^n-\dfrac\tau{6(n+2)}V_nR^{n+2}+O(R^{n+4})$ [Gr90].

\smallskip
{\bf 7.} {\bf Теорема.} Внутренняя изометрия сохраняет скалярную кривизну.

(Докажите с использованием $n$-мерного аналога теоремы о площади.)

\subsection{Главные кривизны}\label{0curmain}


{\bf 1.} Пусть задана система точек $A_1,\dots,A_s$ с массами $m_1,\dots,m_s$.
{\it Моментом инерции}  этой системы относительно прямой $l$  называется число
$I(l)=m_1|A_1l|^2+\dots+m_s|A_sl|^2$, где $|A_il|$ --- расстояние от точки $A_i$ до прямой $l$.

(a) Пусть $I_+$ и $I_-$ --- наибольшее и наименьшее значения моментов инерции
относительно прямых на плоскости, проходящих через фиксированную точку $O$ (возможно, $I_+=I_-$).
Возьмем одну из прямых $l_+$, для которой $I(l_+)=I_+$.
Тогда $I(l)=I_+\cos^2\varphi+I_-\sin^2\varphi$, где $\varphi=\angle(ll_+)$.

(b) В пространстве существуют три такие попарно перпендикулярные прямые $l_1,l_2,l_3$, проходящие через
$O$, что для любой прямой $l$, проходящей через $O$, выполнено
$I(l)=I(l_1)\cos^2(ll_1)+I(l_2)\cos^2(ll_2)+I(l_3)\cos^2(ll_3)$.

(c)* Могут ли прямые с таким свойством не быть перпендикулярными для некоторых $A_1,\dots,A_s$, $m_1,\dots,m_s$?

\smallskip
{\bf Коориентацией} поверхности $\Pi$ называется семейство единичных векторов $n(P)$, нормальных к поверхности
(т.е. перпендикулярных к касательной плоскости в точке $P$, см. определение в пункте \ref{0surgeo}) и непрерывно зависящих от точки $P\in\Pi$.

{\it Кривизной (непараметризованной) кривой на коориентированной поверхности}
называется проекция ускорения на нормаль при движении по этой кривой с единичной скоростью.

{\it Кривизна в указанном смысле} совпадает по модулю, но не обязательно по
знаку, с {\it кривизной} соответствующей непараметризованной кривой.

\begin{figure}[h]\centering
\includegraphics{p.20}
\caption{Нормальное и `косое' сечения поверхности}\label{f:maincur}
\end{figure}

\smallskip
Простейшие инварианты объемлемой изометрии появились еще в XVIII веке при
решении следующей проблемы (рис. \ref{f:maincur}).
Выберем коориентированную поверхность $\Pi$ и точку $P$ на ней.
Как зависит от плоскости $\alpha$, проходящей через точку $P$, кривизна в точке
$P$ (непараметризованной) кривой $\alpha\cap\Pi$?

{\bf Главными кривизнами} $\lambda_+$ и $\lambda_-$ коориентированной
поверхности $\Pi$ в точке $P\in\Pi$ называются наибольшее и наименьшее значения
кривизн в точке $P$ (непараметризованных) кривых пересечения
поверхности с плоскостями, проведенными через нормаль в точке $P$.

{\it Главным направлением} коориентированной поверхности $\Pi$ в точке
$P\in\Pi$, отвечающим данной главной кривизне $\lambda_\pm$, называется
направление той прямой (в касательной плоскости к $\Pi$ в точке $P$),
для которой кривизна пересечения поверхности с
плоскостью, проходящую через эту прямую и нормаль, равна $\lambda_\pm$.


\smallskip
{\bf 2.} Задав коориентацию, найдите главные кривизны и главные направления
в точках поверхностей из задачи \ref{0surexa}.1.abcde*f.

\smallskip
{\bf 3.} Как изменяются главные кривизны при

(a) изменении коориентации на противоположную? \qquad
(b) гомотетии пространства?

\smallskip
{\bf 4.} (a) Объемлемая изометрия сохраняет главные кривизны.

(b) Внутренняя изометрия может не сохранять главные кривизны.

\smallskip
{\bf 5.} {\bf Теорема.}
(a) Если в точке $P$ поверхности $\Pi$ главные кривизны одного знака, то для некоторого $\varepsilon>0$
пересечение $\Pi$ и шара $B(P,\varepsilon)$ в $\R^3$ с центром в $P$ и радиуса $\varepsilon$ лежит
по одну сторону от касательной плоскости к $\Pi$ в точке $P$.

(b) Если в точке $P$ поверхности $\Pi$ главные кривизны разного знака, то ни для какого
$\varepsilon>0$ пересечение $\Pi\cap B(P,\varepsilon)$ не лежит по одну сторону
от касательной плоскости к $\Pi$ в точке $P$.


\smallskip
{\bf 6.} (0) На бесконечном круговом конусе с углом раствора (т.е. максимальным
углом между образующими) $\pi/2$ взята точка $P$ (отличная от вершины).
Через нормаль к конусу в точке $P$ проведена плоскость под углом
$\pi/3$ к образующей конуса (направленной от вершины).
Найдите кривизну в точке $P$ кривой пересечения плоскости и конуса.

(a) Найдите кривизну $k(\varphi)$ в начале координат кривой пересечения поверхности
$z=ax^2+2bxy+cy^2$ с плоскостью, проведенной через ось $Oz$ под углом $\varphi$ к оси $Ox$.

(b) Предположим, что поверхность $z=f(x,y)$ касается плоскости $Oxy$ в начале координат $O$,
т.е. что $f(0,0)=f_x(0,0)=f_y(0,0)=0$.
Обозначим через
$h=\left(\begin{array}{cc} f_{xx} & f_{xy} \\ f_{xy} & f_{yy}
\end{array}\right)$
гессиан.
Тогда главные кривизны являются корнями уравнения
$$(\lambda-f_{xx})(\lambda-f_{yy})=f_{xy}^2,\quad\mbox{т.е.}\quad
\det(h-\lambda E)=0,$$
а главные направления соответствуют cобственными векторам гессиана.
(Инвариантное определение и геометрический смысл соответствующего оператора приведены в \S\ref{0polywein}.)

(c) {\bf Формула Эйлера.}
Пусть $\Pi\subset\R^3$ --- коориентированная поверхность и $P\in\Pi$.
Пусть $k(\varphi)$ --- кривизна в точке $P$ кривой пересечения поверхности с
плоскостью, проведенной через нормаль в точке $P$ под углом $\varphi$ к
тому лучу, для которого эта кривизна максимальна (т.е. к главному направлению,
отвечающему $\lambda_+$).
Тогда
\linebreak
$k(\varphi)=\lambda_+\cos^2\varphi+\lambda_-\sin^2\varphi$.

(d) Если главные кривизны различны, то главные направления ортогональны.

(e)* Как вычислять главные кривизны для поверхностей, заданных в
параметрическом виде?

\smallskip
{\bf 7.} (0) На бесконечном круговом конусе с углом раствора (т.е. максимальным
углом между образующими) $\pi/2$ взята точка $P$ (отличная от вершины).
Через точку $P$ проведена плоскость под углом $\pi/3$ к образующей конуса (направленной от вершины) и под углом 
$\pi/3$ к нормали к конусу в точке $P$.
Найдите кривизну в точке $P$ кривой пересечения плоскости и конуса.

(a) Как отличаются (для поверхностей, рассмотренных в задаче 2) кривизны
кривых пересечения поверхности с двумя плоскостями ($\alpha$ и $\alpha_n$ на
рис. \ref{f:maincur}),
содержащими точку $P$ и пересекающими касательную плоскость  к поверхности по
одной и
той же прямой, одна из которых проходит через нормаль, а другая
под углом $\theta$ к нормали?

Указание к (0,a): если не получается, то см. далее.

(b) Проекция на нормаль в точке $P=\gamma(0)$ ускорения параметризованной
кривой $\gamma$ на поверхности $\Pi$ зависит только от скорости
$\gamma'(0)$ этой кривой в точке $P$.

(c) {\bf Теорема Менье.}
Обозначим через $k$ и $k_n$ кривизны кривых пересечения поверхности с двумя
плоскостями ($\alpha$ и $\alpha_n$ на рисунке),
содержащими точку $P$ и пересекающими касательную плоскость  к поверхности по
одной и той же прямой, одна из которых проходит через нормаль, а другая
под углом $\theta$ к нормали.
Тогда $k\cos\theta=k_n$.

(d) Определим отображение $d(\varepsilon):\Pi\to\R^3$ формулой
$d(\varepsilon)(P)=P+\varepsilon n$.
Тогда проекция из (b) равна $\lim\limits_{\varepsilon\to0}
\dfrac{[d(\varepsilon)'_P(a)]^2-a^2}{2\varepsilon}$, где $a=\gamma'(0)$.

(e) Если главные кривизны в каждой точке поверхности равны нулю, то эта поверхность является частью плоскости.

\smallskip
Можно доказать, что $\tau=2\lambda_+\lambda_-$.

Главные кривизны многомерных коориентированных поверхностей определяются более сложно, см. \S\ref{0polywein}.


\subsection{Полная средняя кривизна}\label{0curmean}

{\it $\varepsilon$-окрестностью} фигуры $M$ (на плоскости или в пространстве)
называется множество $M_\varepsilon$ точек, удаленных от некоторой точки
фигуры $M$ не более, чем на $\varepsilon$:
$$M_\varepsilon:=\{x\ :\ |x-y|<\varepsilon\mbox{ для некоторой }y\in M\}.$$

{\bf 0.} Нарисуйте $\varepsilon$-окрестность в плоскости и найдите ее
периметр и площадь для

(a) квадрата со стороной 1; \quad

(b) выпуклого многоугольника площади $S$ и периметра $P$.

(c) выпуклого множества площади $S$ и периметра $P$.

\smallskip
{\bf 1.} Нарисуйте $\varepsilon$-окрестность в пространстве и найдите ее объем
и площадь ее поверхности  для

(a) куба с ребром 1;

(b) правильной треугольной призмы с длинами ребер 1;

(с) правильного тетраэдра с ребром 1;

(d) произвольного выпуклого многогранника (решите сами, какие нужно задавать
данные; аккуратно докажите Ваше утверждение о старшем коэффициенте).


\smallskip
{\bf 2.} (a) Коробки имеют форму прямоугольных параллелепипедов.
Можно ли в одной коробке пронести другую коробку с большей суммой измерений по длине, ширине и высоте?

(b) Если выпуклый многогранник $M$ с длинами ребер $l_i$ и двугранными углами
$\alpha_i$ содержится в шаре радиуса $R$, то $\sum l_i(\pi-\alpha_i)\le8\pi R$.

\smallskip
Поверхность (граница) фигуры $F$ обозначается $\partial F$.

{\it Полной средней кривизной} выпуклого многогранника $M$ называется число
$$H(\partial M):=\lim\limits_{\varepsilon\to0}
\frac{S(\partial M_\varepsilon)-S(\partial M)}{\varepsilon}.$$
В задаче 2d Вы доказали, что $H(\partial M)=\sum l_i(\pi-\alpha_i)$ для выпуклого многогранника $M$ с длинами ребер $l_i$ и двугранными углами $\alpha_i$.

Теперь рассмотрим поверхность $\Pi\subset\R^3$ с коориентацией $n:\Pi\to\R^3$.
(Определение коориентации см. в пункте \ref{0curmain}.)
Обозначим через $\Pi_{n,\varepsilon}:=\{P+\varepsilon n(P)\}_{P\in\Pi}$ поверхность, образованную концами векторов $\varepsilon n(P)$, отложенных от точек $P$ поверхности $\Pi$ (рис. \ref{f:meancur}).

\begin{figure}[h]\centering
\includegraphics{p.30}
\caption{Сдвиг поверхности вдоль семейства нормалей}\label{f:meancur}
\end{figure}

\smallskip
{\bf 3.} (a) Если $r:D\to\R^3$ --- параметризованная поверхность, $n:r(D)\to S^2$ --- коориентация и
$r_{n,\varepsilon}(u,v):=r(u,v)+\varepsilon n(r(u,v))$, то $r(D)_{n,\varepsilon}=r_{n,\varepsilon}(D)$.

(b) $\Pi_{n,\varepsilon}$ действительно поверхность.

\smallskip
{\bf Полной средней кривизной} коориентированной поверхности $\Pi$ называется число
$$\displaystyle H(\Pi,n):=\lim\limits_{\varepsilon\to0}\frac{S(\Pi_{\varepsilon,n})-S(\Pi)}\varepsilon.$$
Можно эвристически `доказать', что полная средняя кривизна мыльной
пленки (т. е. поверхности минимальной площади с данной границей) равна 0.

\smallskip
{\bf 4.} Задав коориентацию, найдите полную среднюю кривизну поверхностей из задачи \ref{0surexa}.1abcde.

\smallskip
{\bf 5.} Как изменяется полная средняя кривизна при

(a) изменении коориентации на противоположную?
\qquad
(b) гомотетии пространства?

\smallskip
{\bf 6.} (a) Полная средняя кривизна аддитивна, т.е. $H(\Pi_1\cup\Pi_2)=H(\Pi_1)+H(\Pi_2)$, если $\partial\Pi_1$ и
$\partial\Pi_2$ --- замкнутые кривые, пересекающиеся по кривой (из обозначений пропущены коориентации на $\Pi_1$ и $\Pi_2$, полученные сужением некоторой коориентации на $\Pi_1\cup \Pi_2$).

(b) Объемлемая изометрия сохраняет полную среднюю кривизну.

(c) Внутренняя изометрия может не сохранять полную среднюю кривизну.

\smallskip
{\bf 7.} Используем обозначения из задачи 3,5.a.
Будем пропускать в формулах аргумент $(u,v)$ функций $r_u$ и $r_v$, а также вместо $n(r(u,v))$ писать $n$.

(a) $S(r(D))=\int\int_D r_u\wedge r_v\wedge n\ dudv$ для коориентации $n(r(u,v))=r_u\times r_v/|r_u\times r_v|$.

(a') $S(r(D))=-\int\int_D r_u\wedge r_v\wedge n\ dudv$ для коориентации $n(r(u,v))=-r_u\times r_v/|r_u\times r_v|$.

(b) Вектор $n(r(u,v))$ перпендикулярен поверхности $r(D)_{n,\varepsilon}$ в ее точке $r(u,v)+\varepsilon n(r(u,v))$.
Т.е. поле нормалей $m(n(r(u,v))):=n(r(u,v))$ задает коориентацию поверхности $r(D)_{n,\varepsilon}$.

(с) $H(r(D),n)=\int\int_D(r_u\wedge n_v\wedge n+n_u\wedge r_v\wedge n)dudv$ для коориентации
\linebreak
$n(r(u,v)):=r_u\times r_v/|r_u\times r_v|$.
В частности, предел, определяющий полную среднюю кривизну, действительно существует.

\smallskip
{\bf 8.*} Для выпуклого ограниченного множества $M\subset\R^3$ и его
$\varepsilon$-окрестности $M_\varepsilon$
$$(a)\qquad V(M_\varepsilon)=V(M)+S(\partial M)\varepsilon+
\frac12H(\partial M)\varepsilon^2+\frac{4\pi}3\varepsilon^3.$$
$$(b)\qquad S(\partial M_\varepsilon)=S(\partial M)+H(\partial M)\varepsilon+
4\pi\varepsilon^2.$$
$$(c)\qquad H(\partial M)=
2\lim\limits_{\varepsilon\to0}
\dfrac{V(M_\varepsilon)-V(M)-S(\partial M)\varepsilon}{\varepsilon^2}=
\lim\limits_{R\to\infty}\dfrac{V(M_R)-4\pi R^3/3}{R^2}.$$

\subsection{Средняя кривизна в точке}\label{0curmeanpt}

Приведем `физическое' определение.
Возьмем распределение масс на коориентированной поверхности, при котором масса каждого ее куска равна полной
средней кривизне этого куска (таким образом, масса куска может быть отрицательной).
Тогда {\it средней кривизной} поверхности в точке называется плотность в этой точке.

Формально, {\bf средней кривизной} коориентированной поверхности $\Pi$ в точке $P$ называется число
$$H=H_{\Pi,P,n}:=\lim\limits_{\diam(\Pi_P)\to0}\frac{H(\Pi_P,n)}{S(\Pi_P)},$$
где $\Pi_P$ --- образы всевозможных прямоугольников, содержащие точку $P$, при всевозможных
параметризациях поверхности.
Полученное число является {\it плотностью} полной средней кривизны относительно площади.

\smallskip
{\bf 1.} Задав коориентацию, найдите среднюю кривизну в точках поверхностей из задачи \ref{0surexa}.1.abcd.

\smallskip
{\bf 2.} Задав коориентацию, найдите знак средней кривизны точек

(a) тора; \quad (b) поверхности вращения функции $f$.

\smallskip
{\bf 3.} Как изменяется средняя кривизна в точке при гомотетии пространства?

\smallskip
{\bf 4.} (a) Напишите определение предела $\lim\limits_{\diam(\Pi_P)\to0}\dfrac{H(\Pi_P,n)}{S(\Pi_P)}$
`на языке $\varepsilon$-$\delta$'.

(b) $H(r(D),n)=\int\int_DH_{r(D),r(u,v),n}|r_u\times r_v|dudv$.
Эта формула дает эквивалентное определение полной средней кривизны, наиболее точно формализующее вышеприведенное  `физическое' определение.

(c) Напишите определения полной средней кривизны плоской кривой и средней кривизны плоской кривой в точке.

(d)* Последняя равна обычной кривизне.

\smallskip
{\bf 5.} (abc*) {\bf Теорема.}
$$H_P=\frac{r_u\wedge n_v\wedge n+n_u\wedge r_v\wedge n}{|r_u\times r_v|}=
-f_{xx}-f_{yy}=
\frac{(r_v^2r_{uu}+r_u^2r_{vv}-2(r_u\cdot r_v)r_{uv})\wedge r_u\wedge r_v}
{|r_u\times r_v|^3},$$
где первая и третья формулы выполнены для коориентации
$n=r_u\times r_v/|r_u\times r_v|$,
а вторая формула выполнена в точке $O$ для поверхности $z=f(x,y)$, касающейся
плоскости $Oxy$ в начале координат и той коориентации, для которой
$n_O=(0,0,1)$.

При доказательстве первой формулы не забудьте доказать существование предела,
определяющего среднюю кривизну в точке.

(d) Вычислите среднюю кривизну в каждой точке поверхности вращения положительной функции $f$,
если нормали направлены к оси вращения.

\smallskip
{\bf 6.} Найдите среднюю кривизну поверхности $z=xy$ в точке $(1,1,1)$.

\smallskip
{\bf 7.} {\bf Теорема.} Минус половина средней кривизны равна полусумме
главных кривизн и равна среднему значению кривизны нормального сечения:
$$-\frac H2=\frac{\lambda_++\lambda_-}2=\frac1\pi\int_0^\pi k(\varphi)d\varphi.$$

\subsection{Полная гауссова кривизна}\label{0curgau}

В этом пункте $\Pi\subset\R^3$ --- поверхность с коориентацией $n:\Pi\to\R^3$.
Поверхность $\Pi_{n,\varepsilon}$ определена в \S\ref{0curmean}.

Коориентированная поверхность $\Pi\subset\R^3$ называется {\it выпуклой}, если
лучи, определенные (закрепленными) нормалями в разных точках, не пересекаются.
(Это определение не согласуется с определением выпуклости множества, но это
не должно привести к путанице.)
Следующий материал интересен даже для выпуклых поверхностей.

{\it Полной гауссовой кривизной} выпуклой коориентированной поверхности $\Pi$ называется число
$$K(\Pi,n)=\lim\limits_{R\to\infty}\frac{S(\Pi_{n,R})}{R^2}.$$

{\bf 1.} Задав коориентацию, найдите полную гауссову кривизну поверхностей из задачи \ref{0surexa}.1.abcd.

(e) то же для произвольной замкнутой выпуклой (т.е. ограничивающей выпуклое тело) поверхности.

\smallskip
{\bf Полной гауссовой кривизной} поверхности $\Pi$ с коориентацией $n$ называется число $K(\Pi,n)$, для которого
$$S(\Pi_{n,\varepsilon})=S(\Pi)+H(\Pi,n)\varepsilon+K(\Pi,n)\varepsilon^2\qquad(*)$$
при всех тех $\varepsilon$, при которых для любых различных $P,Q\in\Pi$ выполнено
$P+\varepsilon n(P)\ne Q+\varepsilon n(Q)$ (т.е. концы нормальных векторов длины $\varepsilon$,
начала которых --- различные точки поверхности, различны).

\smallskip
{\bf 2.} (a) Это определение совпадает с предыдущим для выпуклых поверхностей.

(b) Найдите полную гауссову кривизну тора.

\smallskip
{\bf 3.} Как изменяется полная гауссова кривизна при

(a) изменении коориентации на противоположную?
\quad
(b) гомотетии пространства?

\smallskip
{\bf 4.} (a) Полная гауссова кривизна аддитивна, т.е. $K(\Pi_1\cup\Pi_2)=K(\Pi_1)+K(\Pi_2)$,
если $\partial\Pi_1$ и $\partial\Pi_2$ --- замкнутые кривые, пересекающиеся по кривой (из обозначений пропущены коориентации на $\Pi_1$ и $\Pi_2$, полученные сужением некоторой коориентации на $\Pi_1\cup \Pi_2$).

(b) Объемлемая изометрия сохраняет полную гауссову кривизну.

\smallskip
Теорема Egregium Гаусса утверждает, что {\it внутренняя изометрия сохраняет полную гауссову кривизну.}
Она доказана в задаче \ref{0polyriem}.6.e с использованием римановой метрики и теоремы Гаусса-Бонне. 

\begin{figure}[h]\centering
\includegraphics{p.40}
\caption{Сферический образ поверхности}\label{f:sphimag}
\end{figure}

Если отложить вектор $n(P)$ от начала координат, то его конец будет лежать на единичной сфере.
Построенное отображение $n:\Pi\to S^2\subset\R^3$ называется {\bf сферическим} или {\it гауссовым}.
Cферическим отображением называется также композиция $n\circ r:D\to S^2\subset\R^3$.
Поверхность $n(\Pi)\subset S^2\subset\R^3$ называется {\bf сферическим} или {\it гауссовым}
{\bf образом} коориентированной поверхности $\Pi$ (рис. \ref{f:sphimag}).

\smallskip
{\bf 5.} Площадь сферического образа выпуклой поверхности равна ее полной гауссовой кривизне.

\smallskip
Пусть нормали к различным точкам коориентированной поверхности не сонаправлены.
Определим {\it площадь сферического отображения (со знаком)} как площадь сферического образа со знаком плюс (или минус),
если при обходе границы поверхности по часовой стрелке (относительно нормалей) граница сферического образа
обходится по часовой стрелке (или против часовой стрелки).

\smallskip
{\bf 6.} Площадь сферического отображения отрицательна для седлообразной поверхности $z=x^2-y^2$, $x^2+y^2\le1$.

\smallskip
Если коориентированную поверхность можно разбить на конечное число частей, на каждой из которых
нормали к различным точкам не сонаправлены, то площадью ее сферического отображения называется
сумма площадей его сужений на части.
Корректность этого определения, т.е. независимость от разбиения поверхности,
фактически доказывается через его эквивалентность следующему определению.

Далее в этом пункте $r:D\to\R^3$ ---  параметризованная поверхность (не обязательно инъективная!) с коориентацией
$n:r(D)\to\R^3$.
Ее {\bf площадью} называется число
$$S(r):=\int\int_D r_u\wedge r_v\wedge n\ dudv.$$
Здесь и далее в этом пункте пропускаем в формулах аргумент $(u,v)$ функций $r_u$ и $r_v$, а также вместо $n(r(u,v))$ пишем $n$.

\smallskip
{\bf 7.} (a) Это число $S(r)$ может быть разным для разных $r$ при одинаковом $r(D)$.

(b) Вектор $n(r(u,v))$ перпендикулярен поверхности $n(r(D))$ в ее точке $n(r(u,v))$.
Т.е. поле нормалей $m(n(r(u,v))):=n(r(u,v))$ задает коориентацию поверхности $n(r(D))$.



\smallskip
{\bf 8.} (a) При данных $\Pi$ и коориентации $n:\Pi\to\R^3$ площадь сферического отображения $n\circ r:D\to\R^3$
 не зависит от инъективной параметризованной поверхности $r:D\to\R^3$, для которой $\Pi=r(D)$.

(b) {\bf Теорема.} Площадь сферического отображения коориентированной поверхности равна ее полной гауссовой кривизне:
$K(r(D),n)=\int\int_Dn_u\wedge n_v\wedge n\ dudv$ для коориентации $n=r_u\times r_v/|r_u\times r_v|$.

В частности, определение гауссовой кривизны осмысленно, т.е. $S(\Pi_\varepsilon)$ действительно выражается
формулой (*) с некоторыми (не зависящими от $\varepsilon$) $H(\Pi,n)$ и $K(\Pi,n)$.

\subsection{Гауссова кривизна в точке}\label{0curgaupt}


Вновь начнем с `физического' определения.
Возьмем распределение масс на коориентированной поверхности, при котором масса
каждого ее куска равна полной гауссовой кривизне этого куска (таким образом,
масса куска может быть отрицательной).
Тогда {\it гауссовой кривизной} поверхности в точке называется плотность в этой точке.

Формально, {\bf гауссовой кривизной} коориентированной непараметризованной
поверхности $\Pi$ в точке $P$ называется число
$$\displaystyle K=K_{\Pi,P}:=
\lim\limits_{\diam(\Pi_P)\to0}\frac{K(\Pi_P)}{S(\Pi_P)}$$
где $\Pi_P$ --- образы всевозможных прямоугольников, содержащие точку $P$, при всевозможных
параметризациях поверхности.
Это {\it плотность} полной гауссовой кривизны относительно площади.

\smallskip
{\bf 1.} Найдите гауссову кривизну в точках поверхностей из задач \ref{0surexa}.1.abcd.

\smallskip
{\bf 2.} Найдите знак гауссовой кривизны точек
\quad

(a) тора; \quad (b) поверхности вращения.

\smallskip
{\bf 3.} (a) Напишите определения `окружностного образа' плоской кривой,
полной гауссовой кривизны плоской кривой и гауссовой кривизны плоской кривой в точке.

(b) Последняя равна обычной кривизне.

\smallskip
{\bf 4.} Как изменяется гауссова кривизна в точке при гомотетии пространства?

\smallskip
{\bf 5.} (abc) {\bf Теорема.}
$$K_P=\frac{n_u\wedge n_v\wedge n}{|r_u\times r_v|}=f_{xx}f_{yy}-f_{xy}^2=
\frac{(r_{uu}\wedge r_u\wedge r_v)(r_{vv}\wedge r_u\wedge r_v)-
(r_{uv}\wedge r_u\wedge r_v)^2}{|r_u\times r_v|^4},$$
где первая формула выполнены для коориентации $n=r_u\times r_v/|r_u\times r_v|$,
а вторая формула выполнена в точке $O$ для поверхности $z=f(x,y)$, касающейся плоскости $Oxy$ в начале координат.

При доказательстве первой формулы не забудьте доказать существование предела,
определяющего гауссову кривизну в точке.

(d) Вычислите гауссову кривизну в каждой точке поверхности вращения положительной функции $f$.

\smallskip
{\bf 6.} Найдите гауссову кривизну поверхности $z=xy$ в точке $(1,1,1)$.

\smallskip
{\bf 7.} (a) {\bf Следствие.} Главные кривизны являются корнями уравнения $\lambda^2+H\lambda+K=0$.

(b) {\bf Следствие.} Гауссова кривизна равна произведению главных кривизн: $K=\lambda_+\lambda_-$.

(c) В теореме \ref{0curmain}.5 можно заменить условие одинаковости (различности) знака главных кривизн на $K>0$ ($K<0$).

(d) Если гауссова и средняя кривизны в каждой точке поверхности равны нулю, то эта поверхность является частью плоскости.

\smallskip
{\bf Теорема.} Для двумерной поверхности в $\R^3$ имеем $\tau=2K$
(где $\tau$ --- скалярная кривизна, см. \S\ref{0cursca}).

\smallskip
Это следует из Теоремы Гаусса-Бонне (задача \ref{0cursec}.2.e) и указания к задаче \ref{0polytr}.8.


\subsection{Секционная кривизна}\label{0cursec}

{\bf Полной секционной кривизной} $\sigma(\Pi)$ двумерной поверхности $\Pi$ с гладкой границей $\partial\Pi$ называется угол между касательным вектором в точке границы $\partial\Pi$ и вектором, полученным из него параллельным переносом вдоль кривой $\partial\Pi$.
Здесь поверхность $\Pi$ ориентирована, имеется в виду ориентированный угол на ориентированной поверхности, и кривая  $\partial\Pi$ направлена так, что ее направление вместе с нормальным к $\partial\Pi$ вектором, касающимся поверхности $\Pi$, дает данную ориентацию на $\Pi$.


\smallskip
{\bf 1.} (a) Объемлемая изометрия сохраняет полную секционную кривизну.

(b) {\bf Теорема.} Внутренняя изометрия сохраняет полную секционную кривизну.

{\it Указание.} Используйте (доказанную ниже) инвариантность параллельного переноса при внутренней изометрии.

(c) Полная секционная кривизна аддитивна, т.е. $\sigma(\Pi_1\cup\Pi_2)=\sigma(\Pi_1)+\sigma(\Pi_2)$, если $\partial\Pi_1$ и $\partial\Pi_2$ --- замкнутые кривые, пересекающиеся по кривой.
(На $\Pi_1$ и $\Pi_2$ рассматриваются ориентации, полученные сужением некоторой ориентации на $\Pi_1\cup \Pi_2$).

(d) Полная секционная кривизна не зависит от ориентации.

\smallskip
{\bf Производной отображения $n:\Pi\to S^2$ из поверхности в сферу} называется семейство отображений
$dn_P:T_P\to T_{n(P)}$ касательной плоскости $T_P$ (в точке $P$ к $\Pi$) в
касательную плоскость $T_{n(P)}$ (в точке $n(P)$ к $S^2$), для которого
$$n(Q)=n(P)+dn_P(\pr(Q-P))+o(|Q-P|)\quad\text{при}\quad\Pi\ni Q\to P.$$
Здесь $\pr:\R^m\to T_P$ --- ортогональная проекция.

\smallskip
{\bf 2.} Определение сферического отображения $n:\Pi\to S^2$ дано в \S\ref{0curgau}.

(0) Найдите производную в точке $(0,0,1)$ вращения $R_{x=0}^{\pi/2}:S^2\to S^2$.

(a) Вектор, касающийся поверхности $\Pi$ в ее точке $P$, касается сферы $S^2$ в точке $n(P)$.

(b) Для любой точки $X\in\Pi$ вектор $dn_P(PX)$ равен вектору $PX$ как свободный вектор.


(c) При семействе отображений $dn_P$ (оно называется {\it сферическим отображением касательных пространств})
семейство векторов, параллельное вдоль некоторой кривой, переходит в семейство векторов,
параллельное вдоль ее образа при сферическом отображении.

(d) Угол поворота вектора при параллельном переносе вдоль замкнутой кривой
равен углу поворота вектора при параллельном переносе вдоль сферического
образа этой кривой: $\sigma(\Pi)=\sigma(n(\Pi))$.

(e) {\bf Теорема Гаусса-Бонне.} Угол поворота касательного к двумерной
поверхности в $\R^3$ вектора при параллельном переносе вдоль границы поверхности равен
полной гауссовой кривизне этой поверхности: $\sigma(\Pi)=K(\Pi)$.
(См. также теорему \ref{0curgau}.8.b.)

\smallskip
Возьмем распределение масс на коориентированной поверхности, при котором масса каждого ее куска равна полной
секционной кривизне этого куска (таким образом, масса куска может быть отрицательной).
Тогда {\it секционной кривизной} поверхности в точке называется плотность в этой точке.

Формально, {\bf секционной кривизной} поверхности $\Pi$ в точке $P$ называется число
$$\displaystyle\sigma(P):=
\lim\limits_{\diam(\Pi_P)\to0}\frac{\sigma(\Pi_P)}{S(\Pi_P)},$$
где $\Pi_P$ --- образы всевозможных прямоугольников, содержащие точку $P$,
при всевозможных параметризациях поверхности.
Это {\it плотность} полной секционной кривизны относительно площади.

\smallskip
{\bf 3.} Как изменяется секционная кривизна в точке при гомотетии пространства?

\smallskip
Для двумерной поверхности в $\R^m$ выполнено $\tau=2\sigma$
(где $\tau$ --- скалярная кривизна, см. \S\ref{0cursca} и указание к задаче \ref{0polytr}.8).
Итак, для двумерной поверхности в $\R^3$ скалярная, секционная и гауссова
кривизны совпадают (с точностью до множителя 2).
Заметим, что здесь скалярная и секционная кривизны определены для двумерной поверхности в $\R^m$ при $m>3$,
а гауссова --- нет.

\subsection*{Указания и решения к некоторым задачам}
\addcontentsline{toc}{subsection}{Указания и решения к некоторым задачам}

{\bf \ref{0cursca}.1.} (b,c) Используйте задачи 2.b и \ref{0surlen}.5.aa'.

\smallskip
{\bf \ref{0cursca}.2.} (a) Уменьшается в $k^2$ раз, где $k$ --- коэффициент гомотетии.

\smallskip
{\bf \ref{0cursca}.5.} (b) Используйте задачу 2.b.

\bigskip
{\bf \ref{0curmain}.1.} (a) Утверждение вытекает из того, что момент инерции есть функция вида
$$A\cos^2\varphi+2B\cos\varphi\sin\varphi+C\sin^2\varphi= P\cos2\varphi+R\sin2\varphi+S= T\cos2(\varphi+\varphi_0)+S.$$
\ \quad
(b) Докажите, что $I(l)$ есть квадратичная форма от направляющего вектора прямой $l$.

\smallskip
{\bf \ref{0curmain}.2.} (f) Глобально главные кривизны не определены.


\smallskip
{\bf \ref{0curmain}.6.} (a) Уравнение сечения есть $$\gamma(t)=(t\cos\varphi,t\sin\varphi,t^2(a\cos^2\varphi+2b\cos\varphi\sin\varphi+c\sin^2\varphi) ).$$
Примените теорему \ref{0curture}.3.b.

(b) Проверяется вычислениями, аналогичными (a).
Невычислительное доказательство этого результата (а также формул для $H$ и $K$
далее) получается, если интерпретировать гессиан как матрицу {\it второго
дифференциала} функции $f$ (или {\it второй квадратичной формы} задаваемой ей поверхности, см. \S\ref{0polywein}).

(c) Достаточно доказать для поверхности $z=f(x,y)$, касающейся плоскости $Oxy$ в начале координат $O$.
См. указание к задаче \ref{0curmain}.1.a.

(d) Следует из (c).

(e) Напишите уравнение касательной плоскости и нормали для поверхности $r(u,v)$ и используйте (b).

\smallskip
{\bf \ref{0curmain}.7.} (a) Ответ приведен в (c).

(b) Обозначим $n=n(\gamma(t))$.
Тогда
$$n\cdot\gamma'=0\quad\Rightarrow\quad n'\cdot\gamma'+n\cdot\gamma''=0
\quad\Rightarrow\quad n\cdot\gamma''=-\gamma'\cdot\partial n(P)/\partial\gamma'.$$
Здесь $\partial n(P)/\partial\gamma'$
есть производная векторного поля $n(P)$ в точке $P=\gamma(0)$  по направлению вектора $\gamma'(0)$.

{\it Другое указание.}
Можно считать, что поверхность задана уравнением $z=f(x,y)$,
$f(0,0)=f_x(0,0)=f_y(0,0)=0$, а уравнение плоскости $z=x\ctg\theta$.
На кривой пересечения рассмотрим параметр $y$.
Дифференцируя, получаем $z_y=f_xx_y+f_y$ и $z_y=x_y\ctg\theta$.
Поэтому в точке $(0,0,0)$ имеем $y_y=1$, $z_y=0$ и $x_y=0$, т.е.
скорость кривой пересечения единичная.
Проекция ускорения кривой пересечения на ось $Oz$ равна
$z_{yy}=(f_{xy}+f_{xx}x_y)x_y+f_xx_{yy}+f_{yx}x_y+f_{yy}$.
В точке $(0,0,0)$ имеем $z_{yy}=f_{yy}$, что не зависит от $\theta$.
Поэтому проекция на ось $Oz$ ускорения кривой пересечения в начале координат
не зависит от $\theta$. А поскольку это ускорение лежит в проведенной
плоскости, оно равно $k(\varphi,\theta)=k(\varphi,0)/\cos\theta$.

(c) Следует из (b).

(e) По формуле Эйлера и теореме Менье любое плоское сечение есть прямая.

 \bigskip
{\bf \ref{0curmean}.0.} (b), (c) Ответ: $P(M_\varepsilon)=P+2\pi\varepsilon$,
\ $S(M_\varepsilon)=S+P\varepsilon+\pi\varepsilon^2$.

\smallskip
{\bf \ref{0curmean}.1.} (d) Ответ приведен после задачи 3.

Перенесем параллельно каждый сферический сектор (являющийся частью
$\varepsilon$-окрестности) так, чтобы вершина перешла в начало координат.
Докажите, что почти все лучи, выходящие из начала координат, пересекают
ровно один из перенесенных сферических секторов.

Другое решение получается, если построить на сфере точки, соответствующие
(нормалям к) граням многогранника и дуги, соответствующие (нормалям к) ребрам
многогранника.

\smallskip
{\bf \ref{0curmean}.2.} Используйте результат задачи 1d.

Заметим, что существует тетраэдр, содержащий тетраэдр большего периметра.

\smallskip
{\bf \ref{0curmean}.4.} (e) Найдите и используйте формулу для площади поверхности тора радиусов $R$ и $r$.
(Например, используя задачу \ref{0surare}.3.b или \ref{0surare}.3.c.)

\smallskip
{\bf \ref{0curmean}.5.} (a) Изменяется знак.

Эвристическое рассуждение: изменение нормали `равносильно' перемене мест поверхностей $\Pi_\varepsilon$ и $\Pi$.

Доказательство получается из формулы задачи 7c и аналогичной формулы со знаком `минус' для противоположной нормали,
ср. с задачей 7a'.

\smallskip
{\bf \ref{0curmean}.6.} Используйте аддитивность площади и $(\Pi_1\cup\Pi_2)_\varepsilon=\Pi_{1,\varepsilon}\cup\Pi_{2,\varepsilon}$,
$(\Pi_1\cap\Pi_2)_\varepsilon=\Pi_{1,\varepsilon}\cap\Pi_{2,\varepsilon}$.

\smallskip
{\bf \ref{0curmean}.7.} (a) $r_u\wedge r_v\wedge n=(r_u\times r_v)\cdot n$.

(c) Используйте задачу 3,5.a.

\bigskip
{\bf \ref{0curmeanpt}.2.} Ввиду инвариантности и аддитивности полной средней кривизны, средняя кривизна
поверхности вращения в точке $P$ равна $\lim\limits_{x\to0}\frac{H(\Pi_x)}{S(\Pi_x)}$,
где $\Pi_x$ --- поверхность, образованная вращением $x$-окрестности точки $P$.
Величины $S(\Pi_x)$ и $H(\Pi_x)$ можно вычислить, используя задачу \ref{0surare}.3.b или \ref{0surare}.3.c.

\smallskip
{\bf \ref{0curmeanpt}.5} (a) Примените теорему о среднем.

(b) Возьмем $r(u,v):=(u,v,f(u,v))$.
Будем опускать аргумент $(0,0)$ функций.
Тогда $r_u=(1,0,f_u)$ и $r_v=(0,1,f_v)$.
Так как $|n(u,v)|=1$, то $n_v\perp n=(0,0,1)$.
Значит, $n_v=(n_{v,1},n_{v,2},0)$.
Поэтому $r_u\wedge n_v\wedge n= n_{v,2}= n_v\cdot r_v\overset{*}= -n\cdot r_{vv}=-f_{vv}$.
Здесь равенство (*) получено дифференцированием равенства $n(u,v)\cdot r_v(u,v)=0$ по $v$ в точке $(0,0)$.

(c) Обозначим $m:=r_u\times r_v$.
Тогда $n=m/|m|$.
Значит, $\displaystyle n_v=\frac{m_v}{|m|}-\frac{m|m|_v}{|m|^2}$.
Поэтому $r_u\wedge n_v\wedge n=r_u\wedge m_v\wedge m/|m|^2$.

Раскрывая $m_v$ по формуле Лейбница и учитывая, что
$a\wedge(b\times x)\wedge(a\times y)=(a\cdot b)x\wedge a\wedge y$,
получаем искомую формулу.

Или сначала напишите уравнение касательной плоскости и нормали для поверхности $r(u,v)$.

\smallskip
{\bf \ref{0curmeanpt}.6} Используйте первый абзац указания к 5с. 
Подставляйте конкретные значения в те формулы, которые уже не нужно будет дифференцировать.

\bigskip
{\bf \ref{0curgau}.1.} (e) Разбейте поверхность на две части.
Используйте 2a и 8b.

\smallskip
{\bf \ref{0curgau}.2.} (a) Перейдите к пределу при $\varepsilon\to\infty$.

(b) См. указание к задаче \ref{0curmean}.4.e.

\smallskip
{\bf \ref{0curgau}.3.} (a) Не меняется. Доказательство аналогично задаче \ref{0curmean}.5.a.

\smallskip
{\bf \ref{0curgau}.4.} (a) Используйте задачу \ref{0curmean}.6 и указание к ней.

\smallskip
{\bf \ref{0curgau}.5.} Используйте 8b.

\smallskip
{\bf \ref{0curgau}.8.} (b) Аналогично задаче \ref{0curmean}.7.c используйте задачу \ref{0curmean}.3.a.

\begin{figure}[h]\centering
\includegraphics{p.50}
\caption{Гауссова кривизна точек тора}\label{f:gaucur}
\end{figure}

\bigskip
{\bf \ref{0curgaupt}.2.} Аналогично задаче \ref{0curmeanpt}.2.

(a) Ответ: см. рис. \ref{f:gaucur}.

\smallskip
{\bf \ref{0curgaupt}.5.} См. указание к \ref{0curmeanpt}.

(c) Имеем $K=n_u\wedge n_v\wedge n/|m|=m_u\wedge m_v\wedge m/|m|^4$.
Нужны тождества
$$(a\times x)\wedge(b\times x)\wedge(c\times x)=0\quad \text{и}\quad
(a\times x)\wedge(b\times y)\wedge(x\times y)=(a\wedge x\wedge y)(b\wedge x\wedge y).$$

{\bf \ref{0curgaupt}.6} См. указание к \ref{0curmeanpt}.6.

\smallskip
{\bf \ref{0curgaupt}.7.} (d) Следует из (a) и задачи \ref{0curmain}.7.e.

\bigskip
{\bf \ref{0cursec}.2.} (d) Ввиду (c) достаточно доказать эту теорему для сферы.
Это делается при помощи аппроксимации сферической области сферическими многоугольниками.

(e) Следует из (d).

(e) {\it Другое указание.} Пусть
$v_P$ --- семейство единичных векторов на $\partial\Pi$, параллельное вдоль
$\partial\Pi$, и $a_P$ и $b_P$ --- произвольная ортонормированная пара векторных
полей на $\Pi$.
Имеем (опуская аргумент $P$)
$$-\sin\angle(v,a)d\angle(v,a)=d(\cos\angle(v,a))=d(v\cdot a)=
v\cdot da+dv\cdot a=v\cdot da=\sin\angle(v,a)b\cdot da.$$
Тогда
$$\sigma(\Pi)=
\int\limits_{\partial\Pi}d\angle(v,a)=
-\int\limits_{\partial\Pi}b\cdot da=
-\int\limits_{\partial D}b\cdot a_udu+b\cdot a_vdv=
\int\int_D\left(\frac{\partial(b\cdot a_v)}{\partial u}-\frac{\partial(b\cdot a_u)}{\partial v}\right)dudv=$$
$$=\int\int_D(b_u\cdot a_v-b_v\cdot a_u)dudv=
\int\int_Dn_u\wedge n_v\wedge n\ dudv=K(\Pi).
$$
Здесь предпоследнее равенство справедливо, поскольку
$$\left(\begin{array}{c} a\\b\\n\end{array}\right)_u=
\left(\begin{array}{ccc} 0&\omega_1&\omega_2\\
-\omega_1&0&\omega_3\\
-\omega_2&-\omega_3&0\end{array}\right)
\left(\begin{array}{c}a\\b\\n\end{array}\right)
\quad\mbox{и}\quad
\left(\begin{array}{c} a\\b\\n\end{array}\right)_v=
\left(\begin{array}{ccc} 0&\omega_1'&\omega_2'\\ -\omega_1'&0&\omega_3'\\
-\omega_2'&-\omega_3'&0\end{array}\right)\left(\begin{array}{c}a\\b\\n
\end{array}\right),
$$
откуда $b_u\cdot a_v-b_v\cdot a_u=\omega_2\omega_3'-\omega_2'\omega_3=
n_u\wedge n_v\wedge n$.

\newpage
\section{Полилинейные кривизны поверхностей}\label{0poly}

\subsection{Риманова метрика. Применение к внутренним изометриям.}\label{0polyriem}

Для поверхности $\Pi$ и точки $P\in\Pi$ обозначим через $T_P=T_{P,\Pi}$
касательную плоскость к $\Pi$ в точке $P$.
Напомним, что $D$ --- декартово произведение интервалов (конечных или бесконечных) на прямой.

{\bf Римановой метрикой}
\footnote{Более точно, римановой метрикой, индуцированоой из $\R^3$.}
(или первой квадратичной формой) на $\Pi$ называется семейство билинейных форм
$$g_P:T_P\times T_P\to\R\quad(P\in\Pi),\quad\mbox{определенных формулой}
\quad g_P(a,b)=a\cdot b.$$

{\bf 1.} (a) Длина образа кривой $\gamma:[a,b]\to\Pi$ на поверхности $\Pi$ равна
$\int_a^b\sqrt{g_{\gamma(t)}(\gamma'(t),\gamma'(t))}\ dt$.

(b) Косинус угла между параметризованными кривыми $\gamma,\beta:[-1,1]\to\Pi$ на поверхности $\Pi$ в точке
$P=\gamma(0)=\beta(0)$ равен $\dfrac{g_P(\gamma',\beta')}{\sqrt{g_P(\gamma',\gamma')g_P(\beta',\beta')}}.$
Здесь $\gamma':=\gamma'(0)$ и $\beta':=\beta'(0)$.

\smallskip
{\bf 2.} Риманова метрика симметрична и положительно определена.

\smallskip
{\bf 3.} (a) {\bf Теорема.} Матрица римановой метрики поверхности $r(D)$ в точке $P=r(u_1,u_2)\in r(D)$ в стандартном базисе $(r_1,r_2)$ есть матрица скалярных произведений (т.е. матрица Грама) этого базиса: $g_{ij}=r_i\cdot r_j$.
(Здесь аргумент $(u_1,u_2)$ функций $r_1,r_2$ пропущен.)
 
(b) Вычислите матрицу римановой метрики поверхности $r(D)$ в точке $r(u,v)$ в стандартном базисе для 
параметризованной поверхности $r(u,v)=(u,v,f(u,v))$ (через функцию $f$ и ее частные производные).

(c) То же для параметризованной поверхности $r(u,v)=(x(u,v),y(u,v),z(u,v))$
(через функции $x,y,z$ и их частные производные).

(d) {\bf Теорема.} Для разных параметризаций $r,\widetilde r:D\to\R^3$ одной непараметризованной поверхности и
соответствующих матриц $G,\widetilde G$ римановой метрики поверхности $r(D)=\widetilde r(D)$ в точке 
$r(u_0,v_0)=\widetilde r(\widetilde{u_0},\widetilde{v_0})$ в базисах $(r_u,r_v)$ и $(\widetilde r_u,\widetilde r_v)$ выполнено $\widetilde G=J^TGJ$, где $J=(r^{-1}\circ\widetilde r)'$.

\smallskip
Напомним, что {\it внутренней изометрией} называется отображение поверхностей, сохраняющее длины всех кривых. 

\smallskip
{\bf 4.} {\bf Теорема.} Следующие три условия на отображение поверхностей равносильны: 

(I) отображение является внутренней изометрией.  

(R) отображение переводит (точнее, его производная увлекает) риманову метрику на второй в риманову метрику на первой.

(D) отображение сохраняет расстояния.

\smallskip
{\bf 5.} (a) Внутренняя изометрия (точнее, ее производная) сохраняет длины касательных векторов.

(b) Внутренняя изометрия сохраняет углы между кривыми.

(c) Для поверхности $\Pi$ и точки $P\in\Pi$ определим функцию
$$f=f_{\Pi,P}:\Pi\to\R\quad\mbox{формулой}\quad f(X)=|P,X|^2,$$
где $|P,X|$ --- расстояние по поверхности.
Тогда $f'(P)=0$ и второй дифференциал функции $f$ совпадает с римановой
метрикой.
 
\smallskip
{\bf 6.} (a) Для параметризованной поверхности $r:D\to\R^3$ выразите $r_k\cdot r_{ij}$ через $g_{ij}:=r_i\cdot r_j$ и их производные.

(b) Внутренняя изометрия сохраняет символы Кристоффеля. 

(c) {\bf Теорема.} При внутренней изометрии поверхностей параметризованные
геодезические переходят в параметризованные геодезические.

(Докажите без использования задачи \ref{0surgeo}.3.b.)

(d) {\bf Теорема.} При внутренней изометрии поверхностей семейство векторов, параллельное вдоль некоторой
кривой, переходит в семейство векторов, параллельное вдоль ее образа.

(e) {\bf Теорема Egregium Гаусса.} Внутренняя изометрия сохраняет гауссову кривизну.

\smallskip
Риманова метрика на поверхности задается сопоставлением каждой параметризации $r:D\to\R^3$ некоторого куска этой
поверхности семейства матриц $G=g_{ij}$, в стандартном базисе $(0,1),(1,0)$, билинейных форм
$$\overline g:\R^2\times\R^2\to\R\quad (X\in D),\quad\mbox{определенных формулой}
\quad \overline g(a,b)=r'(X)a\cdot r'(X)b.$$
(Эти билинейные формы обозначаются также $\overline g=r'(X)^*g$ и называются {\it обратными $r$-образами} римановой метрики на $r(D)$.)
Эти матрицы должны быть связаны на пересечениях кусков, как в задаче 3d ниже.
Заметим, что на {\it всей} поверхности (например, на сфере или листе Мебиуса) риманову метрику нельзя задать
семейством матриц.

\smallskip
{\bf 7.} (a) Длина образа кривой $\gamma=(u_1,u_2):[a,b]\to D$
на параметризованной поверхности $r:D\to\R^3$ равна
$$\int_a^b\sqrt{\sum_{i,j}g_{ij}u_i'\cdot u_j'}\ dt=
\int_a^b\sqrt{\overline g_{r(\gamma)}(\gamma',\gamma')}\ dt=
\int_a^b\sqrt{g_{r(\gamma)}(r'\gamma',r'\gamma')}\ dt.$$
Здесь пропущены аргументы $t$ функций $\gamma,u_1,u_2$ и аргумент $(u_1(t),u_2(t))$ функций $r',g_{ij}$.

(b) {\bf Теорема.} Матрица $g_{ij}$ в точке $(u_1,u_2)$ в стандартном базисе есть матрица скалярных
произведений (т.е. матрица Грама) базиса $(r_1,r_2)$: $g_{ij}=r_i\cdot r_j$.
 
(c) $\det G=g_{11}g_{22}-g_{12}^2>0$ в любой точке.
(Здесь и далее через $\det G$ обозначается определитель {\it матрицы} $g_{ij}$ (но не билинейной формы $g$ или $\widetilde g=r'(X)^*g$.)

(d) {\bf Теорема.} Для разных параметризаций $r,\widetilde r$ одной непараметризованной поверхности и
соответствующих матриц $G,\widetilde G$ выполнено $\widetilde G=J^TGJ$, где $J=(r^{-1}\circ\widetilde r)'$.

(e) Пусть $\gamma,\beta:[-1,1]\to D$ --- параметризованные кривые,
причем $\gamma(0)=\beta(0)=X$.
Обозначим $(a_1,a_2):=\gamma'(0)$ и $(b_1,b_2):=\beta'(0)$.
Тогда косинус угла между кривыми $r\circ\gamma$ и $r\circ\beta$ (на
параметризованной поверхности $r:D\to\R^3$) в точке $r(X)$ равен
$$\frac{\sum_{i,j}g_{ij}a_ib_j}{\sqrt{(\sum_{i,j}g_{ij}a_ia_j)(\sum_{i,j}g_{ij}b_ib_j)}}=
\frac{\overline g(\gamma',\beta')}
{\sqrt{\overline g(\gamma',\gamma')\overline g(\beta',\beta')}}=
\dfrac{g_r(r'\gamma',r'\beta')}
{\sqrt{g_r(r'\gamma',r'\gamma')g_r(r'\beta',r'\beta')}}.$$
Здесь $\gamma'$ и $\beta'$ берутся в точке $0$, а $r,r',g_{ij}$ в точке $X$.


(f) $|r_u\times r_v|^2=\det G$.

(g) Площадь поверхности $r(D)$ равна
$\int\int\limits_D\sqrt{\det G_{(u,v)}}\ dudv$.

\smallskip
{\bf 8.} (a) Система координат экспоненциального отображения (\S\ref{0surexp}) является
{\it евклидовой} в точке $P$, т.е. $g_P(u,v)=u\cdot v$ и $g_X(u,v)'_X|_{X=P}=0$
(или, в ортонормированном базисе в касательной плоскости, $g_{ij}(P)=\delta_{ij}$ и $g_{ij}'(P)=0$).

(b) Внутренняя изометрия сохраняет экспоненциальное отображение: если
$f:\Pi\to\Pi_1$ --- внутренняя изометрия, то $\exp_{f(P)}(f'(P)u)=f(\exp_Pu)$.

\smallskip
{\bf 9.} (a) Существует предел, определяющий скалярную кривизну (см. формулу перед задачей \S\ref{0cursca}.6.a).

(b) В системе координат экспоненциального отображения $\tau=\sum_{i,k}g_{ii,kk}$.

(c) Найдите $\tau$ для поверхности $r:D\to\R^m$.


\subsection{Оператор Вейнгартена (вторая квадратичная форма)}\label{0polywein}

{\bf Оператором Вейнгартена} (или оператором формы) коориентированной
поверхности $\Pi\subset\R^3$ называется семейство операторов

\begin{figure}[h]\centering
\includegraphics{p.70}
\caption{Оператор Вейнгартена}\label{f:wein}
\end{figure}

$$\widetilde q_P:T_P\to T_P\quad (P\in\Pi),\quad\mbox{определенных формулой}
\quad\widetilde q_P(a):=-\partial n/\partial a=-(n_{\gamma_a(t)})'_t|_{t=0}.$$
Здесь $\gamma_a:[-1,1]\to\Pi$ --- такая кривая, что $\gamma_a(0)=P$ и
$\gamma_a'(0)=a$.
(Эта {\it деривационная формула Вейнгартена} записывается также в виде
$n_i=-q_i^1r_1-q_i^2r_2$.)

\smallskip
{\bf 1.} (a) Приведенное определение корректно, т.е. $\partial n/\partial a$
лежит в $T_P$ и не зависит от выбора кривой $\gamma_a$.

(b) Найдите оператор Вейнгартена сферы и цилиндра (т.е. для каждой точки найдите его матрицу в выбранном Вами базисе).

(c) Для заданной параметризации $r:D\to\R^3$ поверхности $r(D)$ найдите матрицу оператора Вейнгартена в точке $r(u,v)$ в базисе $r_u(u,v),r_v(u,v)$.


\smallskip
{\it Второй квадратичной формой} коориентированной непараметризованной
поверхности $\Pi\subset\R^3$ называется семейство билинейных форм
$$q_P:T_P\times T_P\to\R\quad (P\in\Pi),\quad\mbox{определенных формулой}
\quad q_P(a,b):=\widetilde q_P(a)\cdot b=-b\cdot\partial n/\partial a.$$
(Формально, это семейство лучше было бы называть {\it второй билинейной формой}).

\smallskip
{\bf 2.} (a) Вторая квадратичная форма единичной сферы равна ее римановой метрике.

(b) Найдите вторую квадратичную форму цилиндра.

(c) {\bf Теорема.} Матрица второй квадратичной формы поверхности
$z=f(x,y)$, касающейся плоскости $Oxy$ в начале координат $O=P$, в стандартном
базисе является гессианом функции $f$.

(d) Для заданной параметризации $r:D\to\R^3$ поверхности $r(D)$ найдите матрицу второй квадратичной формы в точке $r(u,v)$ в базисе $r_u(u,v),r_v(u,v)$.

\smallskip
{\bf 3.} Рассмотрим коориентированную поверхность $\Pi\subset\R^3$.

(a) Проекция на нормаль $n(P)$ в точке $P=\gamma(0)$ ускорения $\gamma''(0)$
параметризованной кривой $\gamma$ на поверхности равна второй квадратичной
форме от вектора скорости этой кривой в точке $P$:
$n(P)\cdot\gamma''(0)=q_P(\gamma'(0),\gamma'(0))$.

(b) Определим отображение $d(\varepsilon):\Pi\to\R^3$ формулой $d(\varepsilon)(P)=P+\varepsilon n(P)$.
Тогда
$$2q_P(a,b)=\lim\limits_{\varepsilon\to0}
\dfrac{d(\varepsilon)'_P(a)\cdot d(\varepsilon)'_P(b)-a\cdot b}\varepsilon.$$

Вторая квадратичная форма на коориентированной поверхности задается сопоставлением каждой параметризации $r:D\to\R^3$ некоторого куска этой поверхности семейства матриц $q_{ij}$, в стандартном базисе $(0,1),(1,0)$, билинейных форм 
билинейных форм
$$\overline q:\R^2\times\R^2\to\R\quad (X\in D),\quad\mbox{определенных формулой}
\quad \overline q(a,b)=q_{r(X)}(r'(X)a,r'(X)b).$$
(Эти билинейные формы обозначаются также $\overline q=r'(X)^*q$ и называются {\it обратными $r$-образами} билинейной формы на $r(D)$.)
Эти матрицы должны быть связаны на пересечениях кусков, как в задаче \ref{0polyriem}.3.d.

\smallskip
{\bf 4.} Вычислите матрицу $q_{ij}$ в стандартном базисе для  \quad

(a) $r(u,v)=(\cos u\cos v,\cos u\sin v,\sin u)$; \quad

(b) $r(u,v)=(a\cos u\cos v,a\cos u\sin v,c\sin u)$; \quad

(c) $r(u,v)=((2+\cos u)\cos v,(2+\cos u)\sin v,\sin u)$;   \quad

(d) $r(u,v)=(f(u)\cos v,f(u)\sin v,\sin u)$.  \quad

\smallskip
{\bf 5.} {\bf Теорема.} Главные кривизны и направления в точке на поверхности
являются собственными числами и направлениями оператора Вейнгартена (или пары
первой и второй квадратичных форм $g$ и $q$, т.е. корнями уравнения
$\det(g-\lambda q)=0$) в этой точке.

\smallskip
{\it Главными кривизнами} и {\it главными направлениями} многомерной
поверхности в $\R^m$ называются собственные значения и собственные векторы ее
оператора Вейнгартена (который определяется аналогично).

\smallskip
{\bf 6.} Сформулируйте и докажите аналоги известных Вам теорем о главных
кривизнах и направлениях для трехмерных поверхностей в $\R^4$.

\smallskip
{\bf 7.*} Пусть $q_{ij}$ --- матрица второй квадратичной формы в стандартном
базисе.

(a) {\it Деривационные формулы Гаусса.}
$r_{ij}=\Gamma^1_{ij}r_1+\Gamma^2_{ij}r_2+q_{ij}n$ [Ra03 (528)].

(b) {\bf Теорема Бонне.}
Две элементарные непараметризованные поверхности в $\R^3$ объемлемо изометричны
тогда и только тогда, когда они имеют параметризации, индуцирующие одинаковые
первые и  вторые квадратичные формы (или одинаковые римановы метрики и
операторы Вейнгартена) [Ra03, \S81].

Указание. Примените теорему единственности для системы уравнений, составленных
из деривационных формул Гаусса и Вейнгартена.

Заметим, что реализуются не все пары форм, а только удовлетворяющие
{\it уравнениям Гаусса и Петерсона-Кодацци} [Ra03, \S82,\S83].



\subsection{Билинейная форма Риччи}\label{0polric}

 Билинейная форма Риччи описывает искажение объема при экспоненциальном отображении.
{\bf Билинейной формой (тензором) Риччи} $n$-мерной поверхности
$\Pi\subset\R^m$ в точке $P\in\Pi$ называется такая симметричная билинейная
форма $\rho=\rho_P:T_P\times T_P\to\R$, что
$$V(\exp(A))=V(A)-\frac16\int_A\rho(u,u)du+O(h^{n+3})\quad\mbox{при}
\quad h=\diam(A\cup P)\to0$$
по измеримым множествам $A\subset T_P$.

Или, эквивалентно, что для любого единичного $n$-мерного
куба $A\subset T_P$ с вершиной в $P$ выполнено
$$V(\exp(hA))=h^n-\frac{h^{n+2}}6\int_A\rho(u,u)du+O(h^{n+3})\quad\mbox{при}
\quad h\to0.$$
Аналогичная формула справедлива с заменой куба $A$ на любое измеримое множество
$A\subset T_P$ и $h^n$ на $h^nV(A)$ ($h^{n+2}$ и $h^{n+3}$ не меняются).

\smallskip
{\bf 1.} (a) Такая симметричная билинейная форма существует и единственна.


(b) В системе координат экспоненциального отображения $\rho_{kl}=\sum_ig_{ii,kl}$.

(c) Внутренняя изометрия сохраняет билинейную форму Риччи.

(d) Найдите компоненты $\rho_{ij}$ для поверхности $r:D\to\R^m$ в базисе $(r_u,r_v)$.

\smallskip
{\bf 2.} Для симметричной билинейной формы $\omega:\R^n\times\R^n\to\R$

(a) существует и единственен такой оператор $\widetilde\omega:\R^n\to\R^n$, что
$\widetilde\omega(u)\cdot v=\omega(u,v)$ для обычного скалярного произведения в $\R^n$.

(b) $(n+2)\int_{B^n}\omega(u,u)du = V_n\tr\widetilde\omega = V_n\sum_i\omega(e_i,e_i)$, где $B^n$ --- единичный шар в $\R^n$ и $V_n$ --- его $n$-мерный объем.

(c) $\int_A\omega(u,u)du = \frac13\tr\widetilde\omega+\frac14\sum_{i<j}\omega(e_i,e_j)$,
если куб $A$ натянут на ортонормированный базис $e_1,\dots,e_n$.

\smallskip
{\bf 3.}  Обозначим через $\Pi\subset\R^m$ поверхность размерности $n$.

(a) {\bf Теорема.} $\tau=\tr\widetilde\rho$.

Здесь оператор $\widetilde\rho:T_P\to T_P$ определен соотношением $\widetilde\rho(u)\cdot v=\rho(u,v)$
для любых $u,v\in T_P$.

(b) $2\rho_P(u,u)=\sum_i\tau_{\exp_P\left<u,e_i\right>,P}$, где $|u|=1$ и
$e_1,\dots,e_{n-1}$ --- ортонормированный базис в ортогональном дополнении к
$u$ в $T_P$ и $\left<u,e_i\right>$ --- двумерная плоскость в $T_P$,
натянутая на векторы $u$ и $e_i$ [BBB06].

(c) Если $e_1,\dots,e_n$ --- ортонормированный базис в $T_P$, то
$\tau_{\Pi,P}=\sum_{i<j}\tau_{\exp_P\left<e_i,e_j\right>,P}$.



\subsection{Тензор кривизны Римана}\label{0polytr}

 В этом пункте $\Pi\subset\R^m$ --- поверхность размерности $n$.

Пусть $A\subset T_P$ --- область с кусочно-гладкой границей $\partial A$,
содержащей точку $P$.
Обозначим через $\sigma(A):T_P\to T_P$ линейный оператор, сопоставляющий
вектору $x\in T_P$ вектор, полученный из $x$ параллельным переносом вдоль
ориентированной кривой $\exp_P(\partial A)$.
Область $A$ считается настолько малой, что $\exp_P(\partial A)$ определено.

\begin{figure}[h]\centering
\includegraphics{p.90}
\caption{Оператор $\sigma(A)$ и оператор секционной кривизны}\label{f:secriem}
\end{figure}

{\bf Оператором секционной кривизны} поверхности $\Pi$ в точке $P\in\Pi$,
отвечающим паре $u,v\in T_P$ линейно независимых векторов,
называется такой линейный оператор $R(u,v)=R(u,v)_P:T_P\to T_P$, что
для параллелограмма $A_{u,v}$, натянутого на $u,v$, имеем
$$\sigma(hA_{u,v})=E+h^2R(u,v)+o(h^2)\quad\mbox{при}\quad h\to0.$$
Здесь ориентированная кривая $\partial(hA_{u,v})$ выходит из $P$ в направлении вектора $u$.

Если $u$ и $v$ линейно зависимы, то положим $R(u,v)=0$.

Для ориентированной двумерной поверхности $\Pi$ и точки $P\in\Pi$ обозначим
через $R_P^\alpha:T_P\to T_P$ поворот на $\alpha$ относительно точки $P$.

В задачах 1.b и 2.b надо, в частности, доказать, что $R(u,v)$ существует.

\smallskip
{\bf 1.} (a) Для $n=2$ имеем $\sigma(hA_{u,v})=R_P^{(1+o(1))S(\exp(hA_{u,v}))}$.

(b) Для $n=2$ имеем $R(u,v)_P=(u\wedge v)\sigma_PR_P^{\pi/2}$.


(c) При параллельном переносе вдоль кривой, лежащей на двумерной подповерхности
в нашей поверхности, параллельность к подповерхности касательного к
поверхности вектора не обязательно сохраняется.

(d)* $u\cdot R(u,v)v=\sigma_{\exp\left<u,v\right>}u\wedge v$.
(Число $\sigma_{\exp\left<u,v\right>}$ называется {\it секционной кривизной}
поверхности $\Pi$ в точке $P$, отвечающей паре линейно независимых векторов $u,v$.)

\smallskip
{\bf 2.} (a) Верно ли, что для $n=3$ выполнено
$R(u,v)=\sigma_{\exp\left<u,v\right>}u\wedge vR_P^{\pi/2}\circ \pr$, где $\left<u,v\right>$ --- плоскость
в $T_P$, натянутая на векторы $u,v$, a $\pr:T_P\to\left<u,v\right>$ --- ортогональная проекция?

(b) Для стандартной сферы~$S^n=\{x_0^2+\dots+x_n^2=1\}$ в точке $(0,\dots,0,1)$ выполнено
$R(u,v)w=u(v\cdot w)-v(u\cdot w)$, где $u,v,w$ --- касательные векторы к~$S^n$.

\smallskip
{\bf 3.} (a)* Оператор $R(u,v)$ секционной кривизны существует и единственен.

(b) Он линейно зависит от $u,v$.


(d) Внутренняя изометрия сохраняет оператор секционной кривизны.

\smallskip
{\bf 4.} {\bf Теорема о симметриях тензора Римана.}

(a) $R(u,v)$ кососимметричен, т.е.
$[R(u,v)x]\cdot y=-[R(u,v)y]\cdot x$.

(b) $R(u,v)$ кососимметричен по $u,v$: $R(u,v)=-R(v,u)$.

(c) {\it Тождество Бьянки (алгебраическое).} $R(u,v)x+R(v,x)u+R(x,u)v=0$.

(d) $[R(u,v)x]\cdot y=[R(x,y)u]\cdot v$.
\qquad

\smallskip
{\bf 5.} {\bf Теорема.}
(a) Существует и единственно такое 4-линейное отображение $R=R_P:(T_P)^4\to\R$, что
$$\exp\phantom{}'(u)x\cdot\exp\phantom{}'(u)y=x\cdot y+\frac13R(u,x,u,y)+o(|u|^2)$$
при любых постоянных $x,y\in T_P$.
Здесь $\exp'(u):T_P\to T_{\exp(u)}$.


(c)* $R_P(u,x)v\cdot y=R_P(u,x,v,y)$.

(d) В системе координат экспоненциального отображения $R_{ijkl}=g_{ij,kl}$.

(Это дает эквивалентное определение оператора секционной кривизны и, значит, тензора кривизны Римана, см. ниже.
Именно так его определял Риман [Ca28, 9.1, стр. 204].)

(e) Найдите компоненты $R_{ijkl}$ для поверхности $r:D\to\R^m$ в базисе $(r_u,r_v)$.

\smallskip
{\bf 6.*} (a) [MF04, Ra04]
{\bf Теорема.} Следующие условия на элементарную $n$-мерную поверхность равносильны:

(1) она изометрична некоторой части пространства $\R^n$;

(2) параллельный перенос по замкнутому контуру переводит каждый вектор в себя;

(3) все ее операторы секционной кривизны в любой точке нулевые (или один тензор
кривизны Римана нулевой, см. ниже).

(b) Если на односвязной поверхности выполнено (3), то поверхность параллелизуема.

\smallskip
{\bf Тензором кривизны Римана} в точке $P$ (многомерной)
поверхности $\Pi\subset\R^m$ называется трилинейное отображение
$$R:(T_P)^3\to T_P,\quad\mbox{определенное формулой}\quad R(u,v,x):=R(u,v)x.$$

{\bf 7.} (a) Вычислите компоненты $R^i_{jkl}$ на~двумерной сфере в~сферических координатах.

(b) {\bf Теорема.} Для двумерной поверхности в $\R^m$ тензор Римана
выражается через секционную кривизну: $R(u,v,x)=\sigma u\wedge vR_P^{\pi/2}(x)$.

(Или, в системе координат экспоненциального отображения, $R_{1212}=\sigma\det(g_{ij})$;
остальные компоненты нулевые или равны $\pm R_{1212}$ ввиду симметрий тензора Римана.)

(c) {\bf Теорема.} Для трехмерной поверхности в $\R^m$ тензор Римана
выражается через билинейную форму Риччи (и скалярную кривизну):
$$R(u,v)x=\rho(u,x)v-\rho(v,x)u+(u\cdot x)\widetilde\rho(v)
-(v\cdot x)\widetilde\rho(u)+\frac\tau2[(v\cdot x)u-(u\cdot x)v].$$
\quad
(d)* {\bf Теорема.} Билинейная форма Риччи равна свертке тензора Римана:
\linebreak
$\rho(u,v)=\sum_i R(e_i,u)v\cdot e_i$.

\smallskip
{\bf 8.}* {\bf Теорема.} Для двумерной поверхности в $\R^m$ имеем
$$S(\exp(A))=S(A)-\frac\tau{12}\int_Au^2du+o(h^4)\quad\mbox{при}\quad h=\diam(A\cup P)\to0$$
по измеримым множествам $A\subset T_P$.
Иными словами, билинейная форма Риччи пропорциональна римановой метрике
с~коэффициентом~$\tau/2$: $2\rho(u,v)=\tau u\cdot v$ [Gr90, MF04, Ra04].

\smallskip
{\it Замечание.} В определении оператора секционной кривизны (и некоторых задачах этого, но не предыдущего, пункта) экспоненциальное отображение можно заменить на произвольное отображение $f:A\to\Pi$
($A\subset T_P$), для которого $f(tv)'_t|_{t=0}=v$ при любом $v\in A$
(ввиду гладкости достаточно выполнения этого условия для базисных векторов).
Например, на отображение, локально обратное проекции на касательную плоскость.


\subsection*{Указания и решения к некоторым задачам}
\addcontentsline{toc}{subsection}{Указания и решения к некоторым задачам}

{\bf \ref{0polyriem}.1.} (a) Следует из задачи \ref{0curves}.3.b, ср. с задачей \ref{0surlen}.3.  

\smallskip
{\bf \ref{0polyriem}.2.} Ибо $a\cdot b=b\cdot a$ и $a\cdot a=|a|^2>0$ для $a\ne0$.

\smallskip
{\bf \ref{0polyriem}.3.} (d) Ибо риманова метрика --- билинейная форма.

\smallskip
{\bf \ref{0polyriem}.4.} Часть `$(R)\Rightarrow (I)$' вытекает из задачи 1a.
Часть `$(I)\Rightarrow(R)$' вытекает из задачи 5a и тождества поляризации $2g(x,y)=g(x+y,x+y)-g(x,x)-g(y,y)$.
Часть `$(I)\Rightarrow(D)$' вытекает из определения расстояния.
Часть `$(D)\Rightarrow(R)$' вытекает из задачи 5c.

\smallskip
{\bf \ref{0polyriem}.5.} (a) 
Пусть $r:\Pi_1\to\Pi_2$ --- внутренняя изометрия, $P\in\Pi$ и $X\in T_P\Pi$. 
Возьмем любую кривую $\gamma:[0,1]\to\Pi$, для которой $\gamma(0)=P$ и $\gamma'(0)=\vec{PX}$. 
Тогда $r'(P)\vec{PX}=(r\circ\gamma)'(0)$. 
Так как $r$ внутренняя изометрия, то $L(\gamma[0,t])=L((r\circ\gamma)[0,t])$. 
Дифференцируя по $t$ это равенство (по правилу дифференцирования интеграла из задачи \ref{0curves}.3.b по верхнему пределу) получаем $|\gamma'(0)|=|(r\circ\gamma)'(0)|$. 
Т.е. $PX=|r'(P)\vec{PX}|$. 

(b) Следует из (a) и задачи 1b. 

(c) Используйте задачу \ref{0surexp}.2.a.

\smallskip
{\bf \ref{0polyriem}.6.} (a) Ответ. Знак скалярного произведения опускается.
$$2r_ir_{ii}=(g_{ii})'_i,\quad 2r_ir_{ij}=(g_{ii})'_j,\quad
2r_ir_{jj}=2(g_{ij})'_j-(g_{jj})'_i, \quad\mbox{где}\quad i\ne j.$$
Другой способ:
$g_{ij}=r_ir_j$, значит $2r_kr_{ij}=(g_{ki})'_j+(g_{kj})'_i-(g_{ij})'_k$.

(b) Следует из (a) и теоремы 4. 

(с) Следует из (b) и теоремы \ref{0surexp}.1.e.

(d) Следует из (b) и теоремы \ref{0surtra}.6.c.

(e) Следует из (d) и теоремы \ref{0cursec}.2.e Гаусса-Бонне.

\smallskip
{\bf \ref{0polyriem}.8.} (a)
Обозначим $e:=\exp_P:T_P\to\Pi$.
Тогда $g_A(u,v)=e'_Au\cdot e'_Av$ для $A,u,v\in T_P$.
(Или, в координатах, $g_{ij,A}=e'_Aa_i\cdot e'_Aa_j$, где $T_P$
отождествлено с $\R^n$ при помощи некоторого ортонормированного базиса
$a_1,\dots,a_n$.)

Так как $e'(P)=\id T_P$, то  $g_P(u,v)=u\cdot v$.
(Или, в координатах, $g_{ij}(0)=a_i\cdot a_j=\delta_{ij}$.)

Так как $\gamma_u(t):=e(P+ut)$ --- геодезическая, то
$e''_P(u,u)=\gamma_u''(0)\perp T_P$ для любого $u\in T_P$.
Тогда $e''_P(u,v)\perp T_P$ для любых $u,v\in T_P$.
Значит, для $u,v,w\in T_P$ имеем
$$g_{P+w}(u,v)'_w=g_{P+wt}(u,v)'_t|_{t=0}=(e'_{P+wt}u\cdot e'_{P+wt}v)'_t|_{t=0}
=e''_P(w,u)\cdot e'_Pv+e'_Pu\cdot e''_P(w,v)=0,$$
поскольку $(e'_{P+wt})'_t|_{t=0}u=e''_P(w,u)$.

{\it Замечание.}
То же доказательство равенства $g_{ij}'=0$ можно изложить и координатно.
Кривая $\gamma(t):=e(at,bt,0,\dots,0)$ является геодезической для любых $a$ и
$b$.
Тогда $\gamma''(t)\cdot e_k=(a^2e_{11}+2abe_{12}+b^2e_{22})\cdot e_k=0$
для любого $k$. (Догадайтесь сами, что такое $e_k$.)
Отсюда $e_{11}\cdot e_k=e_{12}\cdot e_k=e_{22}\cdot e_k=0$.
Аналогично $e_{ij}\cdot e_k=0$ для любых $i,j,k$.
Тогда в этом базисе $(g_{ij})'_k=(e_i\cdot e_j)'_k=
e_{ik}\cdot e_j+e_{jk}\cdot e_i=0$.

{\it Замечание.}
Равенство $g_{ij}'(P)=0$ равносильно тому, что в $\varepsilon$-окрестности
точки $P$ с точностью до $o(\varepsilon)$ параллельный перенос коммутирует с
экспоненциальным отображением, или ортогональности оператора $e'_u$.

\smallskip
{\bf \ref{0polyriem}.9.}
(a) Возьмем ортонормированный базис $v_1,\dots,v_n$ в $T_P$.
Поскольку геодезическая система координат евклидова (8a), то $g_{ij}(u)=\delta_{ij}+o(|u|)$.
Поэтому
$$\sqrt{\det g_{ij}(u)}=\sqrt{1+o_1(|u|)}=1+o_2(|u|).$$
Значит,
$$V_{\Pi,P}(R)=\int_{B(R)}\sqrt{\det g_{ij}(u)}du=V_nR^n+o(R^{n+1}).$$
Ввиду бесконечной дифференцируемости функции $V_{\Pi,P}(R)$ получаем требуемое.

(b) Аналогично, начиная с равенства $g_{ij}(u)=\delta_{ij}+\sum_{k,l}g_{ij,kl}u_ku_l+o(|u|^2)$.
См. детали в решении задачи \ref{0polric}.1.b; используйте \ref{0polric}.2.b.

(c) Используйте (b).

\bigskip
{\bf \ref{0polywein}.3.} (a) Доказано в задаче \ref{0curmain}.7.b.

\bigskip
{\bf \ref{0polric}.1.}
Возьмем ортонормированный базис $v_1,\dots,v_n$ в $T_P$.

(a) Аналогично задаче \ref{0polyriem}.9.a $\sqrt{\det g_{ij}(u)}=1+o(|u|)$.
Так как $V(\exp(A))=\int_A\sqrt{\det g_{ij}(u)}du$, то
ввиду бесконечной дифференцируемости функции $V(\exp(hA))$ получаем требуемое.

(b) Поскольку геодезическая система координат евклидова (\ref{0polyriem}.8.a), то
$g_{ij}(u)=\delta_{ij}+R_{ij}(u,u)+o(|u|^2)$, где $R_{ij}$ --- квадратичная
форма по $u$ (точнее, $R_{ij}(u,u)=\sum_{k,l} g_{ij,kl}u_ku_l$).
Поэтому
$$\sqrt{\det g_{ij}(u)}=
\sqrt{1+2\sum_iR_{ii}(u,u)+o_1(|u|^2)}=1+\sum_iR_{ii}(u,u)+o_2(|u|^2).$$
Так как $V(\exp(A))=\int_A\sqrt{\det g_{ij}(u)}du$,
то формула из определения билинейной формы Риччи будет верна, если положить
$\rho(u,v):=\sum_iR_{ii}(u,v)=\sum_{i,j,k}g_{ii,kl}u_kv_l.$

(c) Очевидно по определению. А также следует из (b).

(d) Используйте (b).

{\it Замечание к (a).}
Используем указание к задаче \ref{0polyriem}.8.a.
Поскольку геодезическая система координат евклидова, то по формуле Тейлора
$e'_uv_i=v_i+\dfrac{e'''_u(u,u,v_i)}2+o(|u|^2)$,\quad
$$2R_{ij}(u,u)=v_j\cdot e'''_u(u,u,v_i)+v_i\cdot e'''_u(u,u,v_j)\quad\mbox{и}
\quad 2\rho(u,u)=\sum_iv_i\cdot e'''_u(u,u,v_i).$$
Объясним смысл выражения $e'''_u(u,u,v_i)$.
Обозначим через $[L,L']$ векторное пространство линейных операторов из
векторного пространства $L$ в векторное пространство $L'$ (все векторные
пространства рассматриваются над $\R$).
Тогда $e':T_P\to [T_P,T_P]$ ($e'$ не обязательно линейно).
Значит, $e''_u:T_P\to [T_P,T_P]$ --- линейный оператор, или
отображение $e''_u:T_P\times T_P\to T_P$, линейное по второму аргументу.
Тогда $e'':T_P\to [T_P,[T_P,T_P]]$ ($e''$ не обязательно линейно).
Поэтому $e'''_u:T_P\to [T_P,[T_P,T_P]]$ --- линейный оператор, или
отображение $e'''_u:(T_P)^3\to T_P$, линейное по второму и третьему аргументам.

\smallskip
{\bf \ref{0polric}.3.} (a) Следует из 2b.

(c) Следует из (a,b).

\bigskip
{\bf \ref{0polytr}.1.}
(a) По определению секционной кривизны (\S\ref{0cursec}) и оператора $\sigma(A)$.

(b) Следует из (a) и $S(\exp(hA_{u,v}))=(1+o(1))S(hA_{u,v})=h^2u\wedge v+o(h^2)$.

(c) Рассмотрите сферу в $\R^3$ (или окружность в $\R^2$).

\smallskip
{\bf \ref{0polytr}.3.} (a) Возьмем `естественную' кривую
$$\gamma:[0,4h]\to D,\quad\mbox{для которой}\quad\gamma(0)=\gamma(4h)=P\quad\mbox{и}\quad
\gamma[0,4h]=\exp\phantom{}_P(h\partial A_{u,v}).$$
Обозначим через $w(t)=\sum_ia_i(t)r_i(\gamma(t))$ результат переноса вектора
$w(0)$ вдоль отрезка $\gamma[0,t]$.
Далее пропускаем аргумент $t$ функций $w(t)$, $a(t)$, $\gamma(t)$ и их
производных, а также аргумент $\gamma(t)$ отображений $r_i$ и $dr_i=r_i'$.
Далее штрих обозначает производную по $t$.

Напомним, что $r_i'(\gamma(t)):\R^n\to T_{\gamma(t)}$ --- линейный оператор.
Для любого $j$ имеем
$$0=w'\cdot r_j=\sum_i[a_i'r_i\cdot r_j+a_i(r_i'\gamma')\cdot r_j].$$
Для фиксированной параметризации $r$ это уравнение параллельного переноса можно
представить в виде
$a'=F_\gamma(\gamma')a$, где $F_P(x):\R^n\to\R^n$ --- линейный оператор.
Тогда $w(1)=\exp[-\int
_0^{4h}F_\gamma(\gamma')dt]w(0)$.
Имеем
$$\int\limits_0^{4h}F_\gamma(\gamma')dt=
\int\limits_0^hF_{vh-vt}(v)dt+\int\limits_0^hF_{uh+vh-ut}(u)dt-
\int\limits_0^hF_{uh+vt}(v)dt-\int\limits_0^hF_{ut}(u)dt.$$
Теперь нужное утверждение получатся путем взятия начальных членов ряда Тейлора
экспоненты оператора и оператора $F_P(x)$.

(b) $\sigma(hA_{u_1+u_2,v})=\sigma(hA_{u_1,v})\sigma(hA_{u_2,v})$.

(d) Используйте инвариантность параллельного переноса при внутренней изометрии.

\smallskip
{\bf \ref{0polytr}.4.} (a) Следует из ортогональности оператора $\sigma(A)$.

(b) Имеем $\sigma(hA_{v,u})\sigma(hA_{u,v})=E$.
Разлагая с точностью до $o(h^2)$, получаем требуемое.



(d) Следует из (abc).

\smallskip
{\bf \ref{0polytr}.5.} (a) Следует из евклидовости метрики.


(d) Следует из (a,c).

(e) Используйте (d).

\smallskip
{\bf \ref{0polytr}.7.} (a) Используйте 7b.

(b) Вытекает из 1a.

(c) Вытекает из 7d.

(d) Следует из 5cd и задачи \ref{0polric}.1.b.


\smallskip
{\bf \ref{0polytr}.8.} Из 1a и 7d вытекает $\rho=\sigma g$, т.е. $\widetilde\rho=\sigma\id$.
Отсюда и из $\tau=\tr\widetilde\rho$ (\ref{0polric}.3.a) получаем $\tau=2\sigma$ и $2\rho=\tau g$.

Было бы интересно найти прямое доказательство (ср. с \ref{0polric}.1.b, \ref{0polric}.3.a).


\newpage
\section{Ковариантное дифференцирование}\label{0di}

\subsection{Примеры тензорных полей}\label{0diex}

 \hfill{\it Собираются, стягиваются с разных мест вызываемые предметы,}

\hfill{\it причем иным приходится преодолевать не только даль, но и давность...}

\hfill{\it В. Набоков, Королек.}

\smallskip
{\it Векторным полем} на поверхности $\Pi$ называется семейство
касательных векторов $v_P\in T_P$ ($P\in\Pi$), непрерывное по $P$.

\smallskip
{\bf 1.} На сфере без северного и южного полюсов задано векторное поле.
При стереографической проекции из {\it северного} полюса (точнее, при
отображении, индуцированном этой проекцией) это поле переходит в постоянное
векторное поле на плоскости.
Переходит ли это поле  в постоянное векторное поле на плоскости
при стереографической проекции из {\it южного} полюса?

\smallskip
{\it Операторным полем}  на поверхности $\Pi$ называется
семейство линейных операторов $A_P:T_P\to T_P$ ($P\in\Pi$), непрерывное по $P$.

{\it Ковекторным полем}  на поверхности $\Pi$ называется
семейство ковекторов (=линейных функционалов) $\varphi_P:T_P\to\R$ ($P\in\Pi$),
непрерывное по $P$.

{\it Полем билинейных отображений (=форм)} на поверхности $\Pi$ называется
семейство билинейных отображений $\omega_P:T_P\times T_P\to\R$ ($P\in\Pi$),
непрерывное по $P$.

{\it Полем $k$-линейных отображений (=форм)} на поверхности $\Pi$
называется семейство $k$-линейных отображений $\omega_P:(T_P)^k\to\R$
($P\in\Pi$), непрерывное по $P$.

Более точно, определенные выше векторные поля называются {\it касательными}
векторными полями.
Аналогичное замечание справедливо для операторных и других рассматриваемых
полей.

\smallskip
{\bf 2.} Существует ли на торе в $\R^3$ поле
\quad

(a) ненулевых векторов? \quad

(b) ненулевых ковекторов? \quad

(с) невырожденных операторов? \quad

(d) положительно определенных симметричных билинейных форм (билинейная форма
$B:V\times V\to\R$ называется {\it положительно определенной},
если $B(a,a)>0$ для любого $a\in V-\{0\}$)?

(e) невырожденных кососимметричных билинейных форм (билинейная форма
$B:V\times V\to\R$ называется {\it невырожденной}, если
для любого $a\in V-\{0\}$ найдется такой $x\in V$, что $B(a,x)\ne0$)?

(f)* операторов $I_P:T_P\to T_P$, для которых $I_P^2=-E$?

\smallskip
{\bf 3.} (abcdef) То же, что и в задаче 2, но для листа Мебиуса.

\smallskip
{\bf 4.} (a)* {\bf Теорема о еже.} {\it На сфере $S^2$ не существует
(касательного) векторного поля из ненулевых векторов.}

(b) Выведите из (a), что на сфере $S^2$ не существует ковекторного поля из
ненулевых  ковекторов.

(def) То же, что в задаче 2def, но для сферы $S^2$.

\smallskip
{\bf 5.} (abef) То же, что в 2abef, но для сферы
$S^3=\{(x,y,z,w)\in\R^4\ |\ x^2+y^2+z^2+w^2=1\}$, плюс

(g) невырожденных трилинейных форм;

(h) невырожденных векторных произведений.


\subsection{Ковариантное дифференцирование функций}\label{0difun}

{\bf 0.} Для функции $f:\R^2\to\R$ напишите определения частной производной, производной в направлении вектора
$(a,b)$, градиента и производной --- все в точке $(x_0,y_0)$.
Имеют ли эти определения смысл для функции $f:S^2\to\R$?

\smallskip
Здесь и далее $\Pi\subset\R^m$ --- поверхность и $f:\Pi\to\R$ --- функция.

Для $P\in\Pi$ и $u\in T_P$ через $\gamma_u=\gamma_{u,P}:[-1,1]\to\Pi$
обозначается произвольная кривая, для которой $\gamma_u(0)=P$ и
$\gamma_u'(0)=u$.

{\bf Производной функции $f:\Pi\to\R$ в точке $P\in\Pi$ по направлению
касательного вектора $u\in T_P$} (точнее, по касательному вектору $u\in T_P$)
называется число $(\nabla_uf)_P:=[f(\gamma_u(t))]'_t|_{t=0}$.

\smallskip
{\bf 1.} (a) Приведенное определение производной корректно, т.е. не зависит от
выбора кривой $\gamma_u$.
(Рекомендую сначала решить следующие пункты в предположении корректности.)

(b) Зададим функцию $f:S^2\to\R$ формулой $f(P)=\sin\angle POZ$, где
$Z=(0,0,1)$ и $O=(0,0,0)$.
Найдите производную этой функции в точке $(x,y,z)=(1/2,0,\sqrt3/2)$ по
направлению касательного вектора $(0,1,0)$.

(c) $\nabla_u(f_1+f_2)=\nabla_uf_1+\nabla_uf_2$.
\quad

(d) $\nabla_u(fg)=(\nabla_uf)g+f\nabla_ug$.

(e) Пусть $r:D\to\Pi$ --- параметризация поверхности.
Выразите $(\nabla_uf)_P$ через координаты $(a,b)$ вектора $u$ в базисе
$(r_x,r_y)$.


\smallskip
{\bf Производной (=дифференциалом) функции $f:\Pi\to\R$} называется семейство
(=поле) $\nabla f$ ковекторов (=линейных функционалов)
$$\{(\nabla f)_P\}_{P\in\Pi},\quad\mbox{заданных формулой}\quad
(\nabla f)_P(u):=(\nabla_uf)_P.$$

{\bf 2.} (a) Отображение $(\nabla f)_P:T_P\to\R$ действительно является
линейным функционалом, т.е.
$(\nabla f)_P(\lambda a+\mu b)=\lambda(\nabla f)_P(a)+\mu(\nabla f)_P(b)$.

(b) Найдите координаты производной функции из задачи 1b в произвольной
точке $P$ сферы и некотором базисе в $T_P^*$ (выберите и укажите базис сами).

(c) $\nabla(f_1+f_2)=\nabla f_1+\nabla f_2$.
\quad

(d) $\nabla(fg)=(\nabla f)g+f\nabla g$.

(e) Пусть $r:D\to\Pi$ --- параметризация поверхности.
Найдите координаты линейного функционала $(\nabla f)_P$ в базисе $(r_x,r_y)_P$.

(f) Пусть $\varphi:D\to D$ --- замена координат.
Выразите базис $(\varphi\circ r)_x,(\varphi\circ r)_y$ через производную
отображения $\varphi$ и базис $r_x,r_y$.

(g) Как преобразуются координаты производной при замене координат
$\varphi:D\to D$?

\smallskip
{\bf 3.} (a) Существует единственный касательный {\it вектор градиента}
$(\grad f)_P$ в точке $P$, для которого
$(\nabla_uf)_P=u\cdot(\grad f)_P$ при любом $u$.

(b) Запишите градиент функции в~полярных координатах в~$\R^2$,
в~сферических координатах на сфере~$S^2$ и в~сферических координатах
в~$\R^3$ (соответственно для функций $\R^2\to\R$, $S^2\to\R$ и $\R^3\to\R$).

(c) Направление наибольшего роста функции в~некоторой точке
задается вектором ее градиента в~этой точке.

(d) {\it Линией уровня} функции $f:\Pi\to\R$ называется множество $f^{-1}(c)$,
где $c\in\R$. Градиент перпендикулярен линии уровня.

(e) Как преобразуются координаты градиента при замене координат?

(f) Запишите в~произвольной системе координат формулу для производной
функции~$f$ в~направлении вектора градиента функции~$g$.

\subsection{Коммутатор векторных полей}\label{0dicom}

 Пусть $u$ и $v$ --- векторные поля на плоскости (или в~$\R^n$, или на
$n$-мерной поверхности в $\R^m$).
При каких условиях существует система координат $r:\R^2\to\R^2$ (или
$r:\R^2\to\Pi$), для которой эти поля являются координатными (т.е.
$u(r(X))=r'_X(1,0)$ и $v(r(X))=r'_X(0,1)$)?
Решение этой просто формулируемой, но важной задачи приводит к следующему
понятию.

{\bf Коммутатором} векторных полей $u$ и $v$ на поверхности $\Pi$ называется
такое векторное поле $[u,v]$, что
$$\nabla_u\nabla_vf-\nabla_v\nabla_uf=\nabla_{[u,v]}f\quad\mbox{для любой
функции}\quad f:\Pi\to\R.$$

{\bf 1.} (a) Такое поле $[u,v]$ существует и единственно.

(b) Коммутатор обладает свойствами

$\bullet$ $[u,v]=-[v,u]$,

$\bullet$ $[\lambda u,v]=\lambda [u,v]$ и

$\bullet$ $[u_1+u_2,v]=[u_1,v]+[u_2,v]$.

(c) Векторное поле~$v$ в~$\R^n$ с~декартовыми координатами
$(x^1,\dots,x^n)$ называется {\it линейным}, если
$v^i(x^1,\dots,x^n)=A^i_kx^k$, где $A$ --- некоторая постоянная матрица.
Докажите, что коммутатор линейных векторных полей есть снова линейное
векторное поле, и выразите его матрицу через матрицы исходных полей.

(d) Найдите выражение для коммутатора в произвольных координатах.

\smallskip
{\bf 2.} Пусть $u_1$ и $u_2$ --- векторные поля на~$\R^n$.

(a) Обозначим через

$\bullet$ $a_1(t)$ интегральную кривую поля $u_1$, для которой $a_1(0)=P$,

$\bullet$ $a_2(t)$ интегральную кривую поля $u_2$, для которой $a_2(0)=P$,

$\bullet$ $b_{1,s}(t)$ интегральную кривую поля $u_1$, для которой
$b_{1,s}(0)=a_2(s)$,

$\bullet$ $b_{2,t}(s)$ интегральную кривую поля $u_2$, для которой
$b_{2,t}(0)=a_1(t)$.

Докажите, что $[u_1,u_2]=0$ тогда и только тогда, когда $b_{1,s}(t)=b_{2,t}(s)$
для любых $P\in\R^n$ и достаточно малых $t,s\in\R$.

(b) Система координат $r:\R^n\to\R^n$, для которой
$u_1(r(X))=r'_X(1,0,0,0,\dots,0)$, $u_2(r(X))=r'_X(0,1,0,0\dots,0)$ и т.д.
существует тогда и только тогда, когда все эти
векторные поля $u_i$ линейно независимы в каждой точке и все их попарные
коммутаторы нулевые.

\smallskip
{\bf 3.} (a) Пусть $u,v$ --- векторные поля на $\R^m$, касающиеся
поверхности~$\Pi\subset\R^m$.
Тогда векторное поле $[u,v]$ тоже касается поверхности~$\Pi$ и его ограничение
$[u,v]|_\Pi$ на~$\Pi$ совпадает с~полем~$[u|_\Pi,v|_\Pi]$, где $u|_\Pi$ и
$v|_\Pi$ --- ограничения на $\Pi$ полей $u$ и $v$.

(b) На~единичной сфере
$$S^3=\{(x,y,z,w)\in\R^4\ |\ x^2+y^2+z^2+w^2=1\}$$
рассмотрим векторные поля
$$u=(-y,x,-w,z), \quad v_1=(-w,-z,y,x)\quad\mbox{и}\quad v_2=(-w,z,-y,x).$$
Вычислите коммутаторы $[u,v_1]$ и $[u,v_2]$.

\smallskip
{\bf 4.} (a) Любое ли ненулевое векторное поле на $\R^m$ можно 'выпрямить',
т.~е. найти систему координат, в~которой компоненты этого поля будут постоянны?

(b) А ковекторное?

(с)* Пусть $A$ --- операторное поле на~$\R^n$.
Каждой паре векторных полей~$u,v$ на~$\R^n$ сопоставим векторное поле
$$N(u,v)=A^2[u,v]-A[Au,v]-A[u,Av]+[Au,Av].$$
Ясно, что отображение $N:T_P\times T_P\to T_P$ является билинейным.
Докажите, что если существует система координат, в~которой матрицы операторов
семейства~$A$ одинаковы, то~$N=0$.

\smallskip
{\bf 5.} (a) Приведите пример двух коммутирующих векторных полей на~$S^3$,
линейно независимых во~всех точках.

(b)* Любые три попарно коммутирующие векторные поля на~трехмерной
сфере~$S^3$ линейно зависимы в~некоторой точке сферы~$S^3$.

\subsection{Ковариантное дифференцирование векторных полей}\label{0divec}


{\bf 0.} Найдите координаты (в выбранной Вами системе координат) векторного поля на
единичной сфере без северного полюса, которое при стереографической проекции
из северного полюса (точнее, при отображении, индуцированном этой проекцией)
переходит в постоянное векторное поле $(0,1)$ на плоскости $z=1$.

\smallskip
Обозначим через $\pr_{T_P}$ ортогональную проекцию на касательную плоскость $T_P$.

{\bf Ковариантной производной векторного поля $v$ на поверхности $\Pi$ в
точке $P\in\Pi$ по направлению касательного вектора $u\in T_P$} называется
вектор
$$(\nabla_uv)_P:=\pr\phantom{}_{T_P}(v_{\gamma_u(t)})'_t|_{t=0}.$$

{\bf 1.} (b) Найдите ковариантную производную векторного поля
$v(r,\varphi)=(\cos\varphi,-(\sin\varphi)/r)$ на плоскости в точке $(0,2)$
в направлении вектора $(1,1)$.

(a,c,d1,d2,e1,e2,g) Сформулируйте и докажите аналоги задач \ref{0difun}.1.a,c,d,e,g для векторных полей.

У задачи 1d два аналога:
$$\nabla_{ u}(f v)=(\nabla_{ u}f) v+f\nabla_{ u} v\quad
\mbox{и}\quad\nabla_{ u}( v_1\cdot v_2)=
(\nabla_{ u} v_1)\cdot v_2+ v_1\cdot\nabla_{ u} v_2.$$
Формулу, аналогичную 1e, найдите

(e1) при условии $g_{ij}=\delta_{ij}$ в данной точке.

(e2) для общего случая.

(f) Даны линейно независимые вектора $u,v$ и $x$ в $\R^3$.
Выразите через их попарные скалярные произведения коэффициенты
разложения по базису $u,v$ ортогональной проекции вектора $x$ на плоскость,
содержащую вектора $u$ и $v$.


\smallskip
{\bf 2.} (a) Кривая на поверхности является геодезической тогда и только тогда,
когда равна нулю ковариантная производная вектора ее скорости вдоль нее.

(b) Семейство векторов является параллельным вдоль кривой на поверхности
тогда и только тогда, когда ковариантная производная этого семейства вдоль этой
кривой равна нулю.

(c) $[u,v]=\nabla_uv-\nabla_vu$.

\smallskip
{\bf 3.} (a) Если риманова метрика локально евклидова (см. определение в задаче \ref{0polyriem}.8.a), то
$\nabla_u\nabla_vw-\nabla_v\nabla_uw=\nabla_{[u,v]}w$.

(b)* $R(u,v)w=\nabla_v\nabla_uw-\nabla_u\nabla_vw+\nabla_{[u,v]}w$ (см. определение тензора Римана в \S\ref{0polytr}).

(с) Для поверхности $r:D\to\Pi$ найдите компоненты $R^i_{jkl}$ тензора Римана в базисе $(r_x,r_y)$.


\smallskip
{\bf Ковариантной производной векторного поля $v$} на поверхности $\Pi$
называется семейство линейных операторов
$$(\nabla v)_P:T_P\to T_P,\quad\mbox{заданных формулой}
\quad(\nabla v)_P( u):=\nabla_{ u} v.$$

{\bf 4.}
(b) Найдите ковариантную производную векторного поля
$v(r,\varphi)=(\cos\varphi,-(\sin\varphi)/r)$ на плоскости.

(a,c,d1,d2,e1,e2,g) Сформулируйте и докажите аналоги задач \ref{0difun}.2.a,c,d,e,g для векторных полей.


 \subsection{Ковариантное дифференцирование тензорных полей}\label{0ditens}

 Пусть $\gamma:[a,b]\to\Pi$ --- параметризованная кривая на поверхности $\Pi$.
Операторное поле $A(t):T_{\gamma(t)}\to T_{\gamma(t)}$ называется
{\bf параллельным вдоль данной кривой} (в смысле Леви-Чивита), если
для любого векторного поля $v(t)\in T_{\gamma(t)}$, параллельного вдоль кривой
$\gamma$, векторное поле $A(v(t))$ параллельно вдоль кривой $\gamma$.

\smallskip
{\bf 0.} (a) Придумайте пример касательного операторного поля, параллельного вдоль любой кривой.

(b) Для любых точки $P$ ориентированной единичной сферы в $\R^3$ и касательного
вектора $v$ в этой точке обозначим через $R^{\pi/2}_P(v)$ вектор, полученный из
$v$ поворотом в касательной плоскости на $\pi/2$ в положительном направлении.
Тогда операторное поле $R^{\pi/2}_P$ параллельно вдоль любой кривой на сфере.

(с) Сформулируйте и докажите аналоги теорем из \S\ref{0surtra} для операторных полей.

\smallskip
Пусть $u,v\in T_P$ --- касательные векторы.
Здесь и далее $v_u(t)\in T_{\gamma_u(t)}$ --- касательное векторное
поле, параллельное вдоль кривой $\gamma_u$, для которого $ v_u(0)= v$.

{\bf Ковариантной производной операторного поля $A$ на поверхности $\Pi$ в
точке $P\in\Pi$ по направлению $ u\in T_P$} называется оператор,
сопоставляющий вектору $v$ вектор
$$(\nabla_{ u}A)_P( v):=
\pr\phantom{}_{T_P}[A_{\gamma_u(t)}( v_u(t))'_t|_{t=0}].$$
{\bf Ковариантной производной} операторного поля $A$ на
поверхности $\Pi$ называется семейство билинейных отображений
$$(\nabla A)_P:T_P\times T_P\to T_P,\quad\mbox{заданных формулой}
\quad (\nabla A)_P( u, v):=(\nabla_{ u}A)( v).$$

{\bf 1.} (b) Для любых точки $P=(x,y,z)$ ориентированной единичной сферы в
$\R^3$ и касательного вектора $v$ в этой точке обозначим через $A_P(v)$ вектор,
полученный из $xv$ поворотом в касательной плоскости
на $\pi/2$ в положительном направлении.
Найдите ковариантную производную операторного поля $A$ на сфере в точке
$(1,0,0)$ по направлению вектора $(0,1,0)$.

(a,c,d1,d2,e1,e2,g) Сформулируйте и докажите аналоги задач \ref{0divec}.1.a,c,d1,d2,e1,e2,g для операторных полей.

(f) {\it Формула Лейбница.}
Если $v$ --- векторное поле (не обязательно параллельное вдоль кривой
$\gamma_u$), то
$(\nabla_uA)_P(v_P)=[\nabla_u(A(v))]_P-A_P([\nabla_uv]_P).$
Напомним, что $[\nabla_u(A(v))]_P=
\pr\phantom{}_{T_P}[A_{\gamma_u(t)}(v_{\gamma_u(t)})'_t|_{t=0}]$.

\smallskip
{\bf 2.} (a,c,d1,d2,e1,e2,g) Сформулируйте и докажите аналоги задач \ref{0divec}.4.a,c,d1,d2,e1,e2,g
для операторных полей.
В e1,e2 обязательно укажите базис!

\smallskip
{\bf 3.} (a) Дайте определение параллельности ковекторного поля вдоль кривой на поверхности.

(b) Сформулируйте и докажите аналоги теорем из \S\ref{0surtra} для ковекторных полей.

\smallskip
{\bf Ковариантной производной} ковекторного поля $\varphi$ на
поверхности $\Pi$ называется поле $(\nabla\varphi)_P:T_P\times T_P\to\R$
($P\in\Pi$) билинейных отображений, заданных формулой
$$(\nabla\varphi)_P( u, v)=(\nabla_{ u}\varphi)_P( v):=
\varphi_{\gamma_u(t)}( v_u(t))'_t|_{t=0}.$$

{\bf 4.} (b) Найдите матрицу (с указанием базиса) ковариантной производной
ковекторного поля, являющегося производной функции $f(x,y,z)=x$ на сфере
$x^2+y^2+z^2=1$ в произвольной точке сферы.

(a,c,d1,d2,e1,e2,g) Сформулируйте и докажите аналоги задач 2a,c,d1,d2,e1,e2,g
для ковекторных полей.

\smallskip
{\bf 5.} {\it Формула Лейбница.} Если $v(t)\in T_{\gamma(t)}$ --- произвольное
векторное поле (не обязательно параллельное вдоль кривой $\gamma_u$), для
которого $v(0)=v$, то
$$(\nabla_u\varphi)(v)=\varphi_{\gamma_u(t)}(v(t))'_t|_{t=0}-\varphi_P(\pr_{T_P}v'(0)).$$

{\bf 6.} (a) Дайте определение параллельности поля билинейных форм вдоль кривой на поверхности.

(b) Сформулируйте и докажите аналоги теорем из \S\ref{0surtra} для полей билинейных форм.

\smallskip
{\bf Ковариантной производной поля $\omega$ билинейных отображений} на
поверхности $\Pi$ называется семейство трилинейных отображений
$(\nabla\omega)_P:T_P\times T_P\times T_P\to\R$ ($P\in\Pi$), заданных формулой
$$(\nabla\omega)_P(u,v,w)=(\nabla_u\omega)_P(v,w):=
\omega_{\gamma_u(t)}(v_u(t),w_u(t))'_t|_{t=0}.$$

{\bf 7.} (b) Найдите какой-нибудь ряд (из трех элементов) трехмерной матрицы
(с указанием базиса) ковариантной
производной билинейной формы $x_1x_2dx_3\wedge dx_4$ на сфере
$x_1^2+x_2^2+x_3^2+x_4^2=1$ в точке $(0,0,0,1)$ сферы.
По определению, эта билинейная форма сопоставляет паре касательных векторов
$(a,b)$ в точке $(x_1,x_2,x_3,x_4)$ сферы число $x_1x_2(a_3b_4-a_4b_3)$.

(a,c,d1,d2,e1,e2,g) Сформулируйте и докажите аналоги задачи 2a,c,d1,d2,e1,e2,g
для полей билинейных отображений.

(h) Найдите значение формы $\varphi d\varphi\wedge d\theta$ на паре касательных
векторов $(a,b,c)$ и $(a',b',c')$ в точке $(0,\frac1{\sqrt2},\frac1{\sqrt2})$
единичной сферы в $\R^3$.

\smallskip
{\it Симметризацией} поля $\psi_P:T_P\times T_P\to\R$ билинейных
отображений называется поле билинейных отображений
$(\Alt\psi)_P:T_P\times T_P\to\R$, определенных формулой
$2(\Alt\psi)_P( u, v):=\psi_P( u, v)-\psi_P( v, u)$.

\smallskip
{\bf 8.} (a) Система дифференциальных уравнений
$\partial f/\partial x^i=\varphi_i(x_1,\dots,x_n)$, $i=1,2,\dots,n$, разрешима
тогда и только тогда, когда $\Alt(\nabla\varphi)=0$
(где $\varphi$ рассматривается как ковекторное поле).

(b) Выведите формулу для $\Alt(\nabla\varphi)$ в произвольных координатах.

(c) $\Alt(\nabla\varphi)=\nabla_u(\varphi( v))-\nabla_v(\varphi( u))- \varphi([ u, v])$.

\smallskip
{\it Дифференциалом} поля $\omega$ дифференциальных $k$-форм (=кососимметричных
$k$-линейных отображений) называется поле
$d\omega$ дифференциальных $(k+1)$-форм, определенных формулой
$$(d\omega)_P(u_0,\dots,u_k):=\Sigma_{\sigma\in S_{k+1}}
(-1)^{sgn\sigma}(\nabla\omega)_P(u_{\sigma(0)},\dots,u_{\sigma(k)}).$$

{\bf 9.} (a) $d\varphi=\Alt(\nabla\varphi)$ для ковекторного поля $\varphi$.

(b) Выведите формулу для $d\omega$ в произвольных координатах.

\subsection{Общее понятие тензора}

{\bf Тензором типа $(p,q)$} на линейном пространстве $V$ (над $\R$)
называется полилинейное отображение $V^p\times (V^*)^q\to\R$.

\smallskip
{\bf 1.}
Найдите размерность линейного пространства всех кососимметрических тензоров на $n$-мерном линейном пространстве.
Тензор называется {\it кососимметрическим}, если как при перестановке любых двух аргументов из $V^p$,
так и при перестановке любых двух аргументов из $(V^*)^q$, его значение меняет знак.

\smallskip
{\bf 2.}
(a) Линейное пространство всех линейных отображений $V^p\to V^q$ канонически (т.е. независимо от выбора базиса) изоморфно линейному пространству тензоров типа $(p,q)$ на $V$, где $V^0:=\R$.

Каноничность изоморфизма означает независимость от выбора базиса.
Из каноничности следует наличие изоморфизма не только для одного $V$, а даже для семейства --- например, для семейства касательных пространств к поверхности в $\R^m$.

(b) Смешанное произведение в $\R^3$ является тензором типа $(3,0)$.
Векторное произведение в $\R^3$ переходит в тензор типа $(2,1)$ при изоморфизме из (a).
Последний тензор также называют векторным произведением.

(c) {\it Операция опускания и поднимания индексов.}
Для любого линейного подпространства $V\subset\R^m$ линейное пространство всех тензоров типа $(p,q)$ на
$V$ канонически изоморфно линейному пространству всех тензоров типа $(p-1,q+1)$ на $V$.

(d) Смешанное произведение переходит в векторное при изоморфизме из (c)
(т.е. эти тензоры связаны операцией опускания и поднимания индексов).

(e) Обозначим через $T_m^n$ линейное пространство всех тензоров типа $(m,n)$.
Линейное пространство всех линейных отображений $T_m^n\to T_p^q$ канонически  изоморфно
линейному пространству $T_x^y$ для некоторых $x,y$.

\smallskip
Базисом в $V^*$, {\bf двойственным} к базису $e_1,\dots,e_n$ в $V$, называется набор ковекторов
$f^1,\dots,f^n$, определенных формулами $f^i(e_j)=\delta_{ij}$.

\smallskip
{\bf 3.} Это действительно базис в $V^*$.

\smallskip
{\bf Компонентами} тензора $T:V^p\times (V^*)^q\to\R$ в базисе $e_1,\dots,e_n$ в $V$ называется массив чисел $$T_{i_1,\dots,i_p}^{j_1,\dots,j_q}(e_1,\dots,e_n):=T(e_{i_1},\dots,e_{i_p},f^{j_1},\dots,f^{j_q}),$$
где $f^1,\dots,f^n$ --- базис в $V^*$, двойственный к базису $e_1,\dots,e_n$ в $V$.

\smallskip
{\bf 4.}
(a) Выпишите компоненты смешанного и векторного произведений в
стандартном базисе в $\R^3$.

(b) Для любого набора компонент существует и единственен соответствующий тензор.

\smallskip
{\bf 5.} 
Обозначим через $e_1,\dots,e_n$ базис в линейном пространстве $V$ и через $f^1,\dots,f^n$ двойственный базис в $V^*$.
Для базиса $e_1,\dots,e_n$ в $V$ тензор $T:V^p\times (V^*)^q\to\R$ обозначается $$\sum\limits_{i_1,\dots,i_p,j_1,\dots,j_q} T(e_{i_1},\dots,e_{i_p},f^{j_1},\dots,f^{j_q})
f^{i_1}\otimes\dots \otimes f^{i_p}\otimes e_{j_1}\otimes\dots \otimes e_{j_q}.$$
\ \quad
(a) Найдите значение тензора $f^1\otimes e_2 + f^2\otimes (e_1+3e_3)$ на паре
$e_1+5e_2+4e_3$, $f^1+f^2+f^3$.

(b) Найдите компоненты тензора $f^1\otimes f^2\otimes e_1 + f^1\otimes f^2\otimes e_2$ в базисе
$(\widetilde e_1,\widetilde e_2):= (e_1,e_2)\left(\begin{matrix} 1&1\\ 2&3 \end{matrix}\right)$.

(c) Как меняются компоненты тензора типа $(p,q)$ при замене базиса с матрицей $A$?

\smallskip
{\bf 6.} Пусть $f:\Pi\to\R$ --- гладкая функция на поверхности в $\R^3$, $P\in\Pi$ и $f'(P)=0$.
Докажите, что вторые частные производные $(\partial^2f/\partial x^i\partial x^j)|_P$ являются компонентами
некоторого тензора ранга~$2$ и определить тип этого тензора.

\smallskip
{\bf 7.} Дайте определение параллельности вдоль кривой на поверхности и
ковариантной производной поля

(a) $k$-линейных отображений.
\qquad
(b) тензоров типа $(p,q)$.

\smallskip
{\bf 8.*} $\nabla_mR^i_{jkl}+\nabla_kR^i_{jlm}+\nabla_lR^i_{jmk}=0$ (см. определение тензора Римана в \S\ref{0polytr}).

\smallskip
{\bf 9.} (a,c,d1,d2,e1,e2,g) Сформулируйте и докажите аналоги задач
\ref{0divec}.2.a,c,d1,d2,e1,e2,g для тензоров типа $(p,q)$.

\smallskip
См. подробнее http://ru.wikipedia.org/wiki/Тензор (24.11.2013 там опечатка) и
\linebreak
http://en.wikipedia.org/wiki/Tensor (24.11.2013 там опечатка исправлена).





\subsection*{Указания и решения к некоторым задачам}
\addcontentsline{toc}{subsection}{Указания и решения к некоторым задачам}

{\bf Ответ к \ref{0diex}.2a, 3a и 5a.} Да.

\smallskip
{\bf Указание к пунктам (b) всех задач из \S\ref{0diex}.} Если $\{v_P\}$ векторное поле, то
ковекторное поле можно определить формулой $\varphi_P(u):=u\cdot v_P$.

\smallskip
{\bf Указание к пунктам (c) всех задач из \S\ref{0diex}.} Да, $A_P(v):=v$.

\smallskip
{\bf Указание к пунктам (d) всех задач из \S\ref{0diex}.} Да, $B_P(u,v):=u\cdot v$.

\smallskip
{\bf \ref{0diex}.2e,4e.} Да. Возьмите ориентированную площадь.

\smallskip
{\bf \ref{0diex}.2f,4f.} Да.
Возьмите поворот касательной плоскости на $\pi/2$ для выбранной ориентации тора или сферы.

\smallskip
{\bf Ответ к \ref{0diex}.3ef,5e.} Нет.

\smallskip
{\bf \ref{0diex}.2,5.} Постройте два или три векторных поля, линейно независимых в каждой точке.
Тогда нужные объекты достаточно построить в одной точке.

\bigskip
{\bf \ref{0difun}.1.} Обозначим
$\gamma_u(t)=r(a(t),b(t))$, $a=a'(0)$, $b=b'(0)$.
Тогда $(\nabla_uf)_P=[f(r(a(t),b(t)))]'_t|_{t=0}$.

(e) Ответ: если
$$u=ar_x+br_y,\quad\mbox{то}\quad
(\nabla_uf)_P=a(f\circ r)_x|_{r^{-1}(P)}+b(f\circ r)_y|_{r^{-1}(P)}.$$

{\bf \ref{0difun}.2.} (e) Ответ: $(\nabla f)_P=((f\circ r)_x|_{r^{-1}(P)},(f\circ r)_y|_{r^{-1}(P)})$.

(g) От производной требуется только то, что она является ковектором.
Используйте предыдущий пункт.

\bigskip
{\bf \ref{0dicom}.5.b.} Иначе $S^3\cong S^1\times S^1\times S^1$.

\bigskip
{\bf \ref{0divec}.1.}
(e1) Пусть $v_{r(x,y)}=p(x,y)r_x(x,y)+q(x,y)r_y(x,y)$ --- векторное поле
на поверхности $r(x,y)$ и $\gamma_u(t)=r(a(t),b(t))$.
Тогда
$$ u=r_xa'(0)+r_yb'(0)\quad\mbox{и}\quad \nabla_{ u} v=
r_x\left([v_{r(a(t),b(t))}]'_t|_{t=0}\cdot r_x\right)+
r_y\left([v_{r(a(t),b(t))}]'_t|_{t=0}\cdot r_y\right).$$
\quad
(e2) Аналогично (e1) используя (f).
Или используйте (e1) и (g).

(g) От производной требуется только то, что она является оператором.
Используйте закон изменения базиса в касательном пространстве при замене
переменных на поверхности (т.е. задачу \ref{0difun}.2.f).

\smallskip
{\bf \ref{0divec}.3.c.} Используйте (b).

\bigskip
{\bf Общее указание к \ref{0ditens}.}
Аналогично результатам о ковариантном дифференцировании функций и векторных полей.

\smallskip
{\bf \ref{0ditens}.4e1,4e2,7e1,7e2.}
Используйте уравнение параллельного переноса векторных полей.

\newpage
\section{Обобщение на римановы многообразия}

\subsection{Элементы гиперболической геометрии Лобачевского}

 Назовем {\it плоскостью Лобачевского} половинку $z\ge0$ двуполостного
гиперболоида $z^2=x^2+y^2+1$.
Назовем {\it прямыми Лобачевского} сечения этой половинки плоскостями,
проходящими через начало координат.

\smallskip
{\bf 1.} (a) Через точку плоскости Лобачевского, не лежащую на данной прямой
Лобачевского, проходит более одной прямой Лобачевского, не пересекающей
данной прямой Лобачевского.

(b) Через любые две точки плоскости Лобачевского проходит ровно одна прямая
Лобачевского.

(c) Для любой кривой $(x(t),y(t),z(t))$ на плоскости
Лобачевского $x_t(t)^2+y_t(t)^2-z_t(t)^2>0$.

\smallskip
{\it Длиной Лобачевского} кривой $(x(t),y(t),z(t))$, $t\in[a,b]$, на плоскости
Лобачевского называется число
$\int_a^b\sqrt{x_t(t)^2+y_t(t)^2-z_t(t)^2}\ dt$.
(Или, выражаясь научно, назовем {\it римановой метрикой Лобачевского}
сужение псевдоримановой метрики $ds^2=dx^2+dy^2-dz^2$ в $\R^3$ на плоскость
Лобачевского.)

Далее риманова метрика $g_{ij}$ {\it записывается} в виде
$ds^2=g_{ij}dx^idx^j$.

\smallskip
{\bf 2.} (a) Плоскость Лобачевского изометрична верхней полуплоскости с
римановой метрикой Лобачевского $ds^2=-4\dfrac{dw d\bar w}{(w-\bar w)^2}$
(модель Пуанкаре в верхней полуплоскости).
Далее {\it плоскостью Лобачевского} называется верхняя полуплоскость с
римановой метрикой Лобачевского.

(b) Плоскость Лобачевского инвариантна относительно преобразований $p(w)=w+a$
и $q(w)=-1/w$.

(c) $\ch|z_1,z_2|=1+\dfrac{|z_1-z_2|^2}{2\mbox{Im}z_1\mbox{Im}z_2}$.

Указание: точки $z_1$ и $z_2$ переводятся изометрией в точки с одинаковыми
абсциссами.

(c) Выведите {\it теорему Пифагора} для плоскости Лобачевского:
$\ch c=\ch a\ch b$

Указание: можно считать $C=i$, $A=ki$ и $B=\cos\varphi+i\sin\varphi$.

(d) Окружность Лобачевского является евклидовой окружностью.

(e)* Найдите длину окружности Лобачевского радиуса $R$.

Указание: движением плоскости Лобачевского центр окружности
можно перевести в центр модели Пуанкаре в круге, затем найти связь евклидова
радиуса и радиуса Лобачевского...

(f)* Сфера, плоскость и плоскость Лобачевского попарно локально не изометричны.

(g)* Любая внутренняя изометрия метрики Лобачевского, сохраняющая ориентацию,
является дробно-линейным
преобразованием $f(z)=\dfrac{az+b}{cz+d}$ с определителем $ad-bc=1$ (и
вещественными $a,b,c,d$).


\subsection{Геометрия на римановых многообразиях}

{\bf Римановым многообразием} (локальным) называется пара $(M,g)$ из открытого
множества в $\R^n$ и поля невырожденных симметричных билинейных форм $g$ на нем.

Это поле называется {\it римановой метрикой.}

При этом изометрического вложения $M\subset\R^m$ не задано!

{\it Длины кривых и площади} определяются через риманову метрику формулами,
полученными ранее.

{\it Касательным пространством $T_P$ в точке $P\in M$} называется пространство
$\R^n$.

{\it Скалярная кривизна, геодезические, экспоненциальное отображение,
тензор Риччи, ковариантное
дифференцирование функций, касательные векторные и тензорные поля} на $(M,g)$
определяются дословно так же, как для поверхностей в $\R^m$.

\smallskip
{\bf 1.} Вычислите скалярную кривизну в точках

(a) плоскости Лобачевского, т.е. верхней полуплоскости с римановой метрикой
$ds^2=-4\dfrac{dw d\bar w}{(w-\bar w)^2}$.
\quad

(b) плоскости с римановой метрикой $ds^2=\lambda(x,y)(dx^2+dy^2)$.

(c)* пространства $\R^n$ с римановой метрикой
$ds^2=\dfrac{dx_1^2+\dots+dx_n^2}{(1+x_1^2+\dots+x_n^2)^2}$.

(d)* произвольного риманова многообразия.

\smallskip
{\bf 2.} (a)* Уравнение геодезической
$$x_k''+\sum_{i,j}\Gamma^k_{ij}x_i'x_j'=0,\quad\mbox{где}\quad
2\Gamma^k_{ij}=\sum_l g^{kl}([g_{lj}]_{x_i}+[g_{li}]_{x_j}-[g_{ij}]_{x_l}).$$
\quad (b) Через каждую точку в каждом направлении проходит ровно одна
геодезическая.

\smallskip
{\bf 3.} Найдите геодезические на верхней полуплоскости с римановой метрикой
\quad

(a) $ds^2=y(dx^2+dy^2)$; \quad (b) $ds^2=-4\dfrac{dw d\bar w}{(w-\bar w)^2}$;
\quad (c) $ds^2=\dfrac{dx^2+dy^2}{x^2+y^2}$.

\smallskip
{\bf 4.*} Найдите все функции $\lambda(x,y,z)$, для которых все кривые
$\{y=c_1,z=c_2\}$ являются геодезическими римановой метрики
$e^{\lambda(x,y,z)}(dx^2+dy^2+dz^2)$ на $\R^3$.


\smallskip
{\bf 5.} Сформулируйте и докажите аналоги всех определений и теорем из \S\ref{0polric}.


\smallskip
Касательное к $M$ векторное поле $v$ на кривой
$\exp_P\circ\gamma:[0,1]\to T_P\to M$ называется
{\bf параллельным вдоль кривой $\exp_P\circ\gamma$}, если его прообраз
$[\exp_P'(\gamma(t))]^{-1}v_{\exp_P(\gamma(t))}$ при производной
экспоненциального отображения параллелен.

\smallskip
{\bf 6.} (a) Это определение совпадает с прежним для поверхностей в $\R^m$.

(b) Выпишите явно и решите уравнение параллельного переноса вдоль
горизонтальной евклидовой прямой для плоскости Лобачевского.

(c) То же вдоль данной евклидовой окружности.

\smallskip
Определение {\it параллельного переноса, секционной кривизны, тензора Римана и
ковариантной производной поля $k$-линейных форм} на $(M,g)$ повторяет приведенное выше.

\smallskip
{\bf 7.*} (a) Для двумерного риманова многообразия $\tau=2\sigma$, т.е.
угол поворота касательного вектора при параллельном переносе вдоль замкнутой
кривой (ориентированной согласованно с ориентацией многообразия) и
ограничивающей область $A$ равен $\frac12\int_A\tau dS$.

(b) Сформулируйте и докажите аналоги всех определений и теорем из \S\ref{0polytr}.


\smallskip
{\bf Ковариантной производной векторного поля $v$ на $M$ в точке $P\in M$
по направлению вектора $u\in T_P$} называется вектор
$\nabla_u v:=([\exp_P'(ut)]^{-1}v_{\exp_P(ut)})'_t|_{t=0}$,
т.е. производная (в точке $0\in\R^n$ по направлению вектора $u$)
прообраза векторного поля $v$ при производной экспоненциального отображения.

\smallskip
{\bf 8.*} (a) Это определение совпадает с прежним для поверхностей в $\R^m$.

(b) Напишите определение ковариантной производной векторного поля на римановом многообразии.

(c) Матрица ковариантной производной векторного поля $v$ в точке $P$ в системе
координат экспоненциального отображения есть $(\dfrac{\partial v^i}{\partial x^j})$.

(d) Найдите эту матрицу (с указанием базиса) для плоскости Лобачевского.

(e) Найдите эту матрицу в произвольных координатах.

(f) Докажите эквивалентность приведенного определения определению
ковариантного дифференцирования векторов из [Gr90, 2.2].

(g) Разности $\Gamma^k_{ij}-\tilde\Gamma^k_{ij}$ символов Кристоффеля двух
римановых метрик $g_{ij}$ и $\tilde{g_{ij}}$ на одном и том же $M$ образуют тензор типа~$(1,2)$.

(h) Любой тензор типа~$(1,2)$ может быть представлен таким образом.

\smallskip
{\bf 9.*} Сформулируйте и докажите аналоги всех определений и теорем из \S\ref{0di}.

\smallskip
{\bf 10.*} Пусть $M$ --- поверхность в $\R^m$.

Кривые $\gamma_1,\gamma_2:[-1,1]\to M$ с условием
$\gamma_1(0)=\gamma_2(0)=P$ называются {\it эквивалентными в точке $P\in M$},
если для любого $\varepsilon>0$ существует такое $\delta>0$, что при
$|t|<\delta$ точки $\gamma_1(t)$ и $\gamma_2(t)$ можно соединить дугой (на
поверхности) длины меньше $\varepsilon t$.
{\it Касательным вектором к $M$ в точке $P\in M$} называется класс
эквивалентности таких кривых.
Через $T_P$ обозначим пространство всех касательных векторов к $M$ в точке
$P\in M$.

(a) Существует взаимно однозначное соответствие между $T_PM$ и $\R^n$.

(b) Определите на этом пространстве операции сложения и умножения на число так,
чтобы получилось векторное пространство (над $\R$).

(c) Постройте аналог понятия риманова многообразия для поверхности $M$ в $\R^m$.

\smallskip
{\bf 11.*} Определения {\it аффинной связности} и заданного ей {\it ковариантного дифференцирования} см., например, в [Ra04, MF04].

(a) На~плоскости с~координатами $u,v$ найдите аффинную связность, относительно
которой векторные поля $\xi=(e^u,1)$ и $\eta=(0,e^v)$ ковариантно постоянны
(т.е. их ковариантная производная равна нулю).

(b) Найдите ({\it тензор кручения}) $\Gamma^k_{ij}-\Gamma^k_{ji}$ для этой
связности.

(c) Существует ли риманова метрика, порождающая эту аффинную  связность?

\smallskip
{\it Указание к 2a.} Доказывается при помощи вариационного исчисления.

\smallskip
{\it Ответ к 8e.}
$(\dfrac{\partial v^i}{\partial x^j}+\sum_k\Gamma^i_{kj}v^k)$.

\smallskip
{\it Указание к 10b.}  Докажите, что для любых двух таких кривых $\gamma_1$ и $\gamma_2$
существует (единственная с точностью до эквивалентности) такая кривая $\gamma$,
что $\nabla_\gamma f=\nabla_{\gamma_1}f+\nabla_{\gamma_2}f$ для любой функции
$f$.
Класс эквивалентности кривой $\gamma$ называется {\it суммой} классов
эквивалентности кривых $\gamma_1$ и $\gamma_2$.

\end{document}